\documentclass[11pt]{amsart}
\usepackage{amsmath}
\usepackage{amsxtra}
\usepackage{amscd}
\usepackage{amsthm}
\usepackage{amsfonts}
\usepackage{amssymb}
\usepackage{eucal}

\numberwithin{equation}{subsection}

\begin{document}

\newcommand{\thmref}[1]{Theorem~\ref{#1}}
\newcommand{\secref}[1]{\S~\ref{#1}}
\newcommand{\defref}[1]{Definition~\ref{#1}}
\newcommand{\lemref}[1]{Lemma~\ref{#1}}
\newcommand{\propref}[1]{Proposition~\ref{#1}}
\newcommand{\corref}[1]{Corollary~\ref{#1}}
\newcommand{\remref}[1]{Remark~\ref{#1}}
\newcommand{\nc}{\newcommand}
\nc{\bs}{\backslash}
\nc{\wh}{\widehat}
\nc{\wt}{\widetilde}
\nc{\inv}{^{-1}}
\nc{\sm}{\setminus}
\nc{\g}{{\mathfrak g}}
\nc{\n}{{\mathfrak n}}
\nc{\h}{{\mathfrak h}}
\nc{\bb}{{\mathfrak b}}
\nc{\pp}{{\mathfrak p}}
\nc{\ab}{{\mathfrak a}}
\nc{\Cx}{{\mathbb C}^{\times}}
\nc{\C}{{\mathbb C}}
\nc{\Z}{{\mathbb Z}}
\nc{\E}{{\mathcal E}}
\nc{\on}{\operatorname}
\nc{\B}{{\mathcal B}}
\nc{\Oo}{{\mathcal O}}
\nc{\Dd}{{\mathcal D}}
\nc{\W}{{\mathcal W}}
\nc{\Gr}{{\mathcal Gr}^X}
\nc{\Gro}{{\mathcal Gr}^{\circ}}
\nc{\Fl}{{\mathcal Fl}}
\nc{\Flo}{\Fl^{\circ}}
\nc{\gln}{\mathfrak{gl}_n}
\nc{\sln}{\mathfrak{sl}_n}
\nc{\Pp}{{\mathbb P}^1}
\nc{\Grp}{{\mathcal Gr}^{\Pp}}
\nc{\Ahat}{{A/A_+}}
\nc{\At}{\wh{A}_{-1}}
\nc{\pline}{\wh{A}_p}
\nc{\AO}{{\mathcal AO}}
\nc{\Op}{{\mathcal Op}}
\nc{\D}{{D}}
\nc{\Dx}{{D}^{\times}}
\nc{\K}{{\mathcal K}}
\nc{\T}{{\mathcal T}}
\nc{\V}{{\mathcal V}}
\nc{\Ll}{{\mathcal L}}
\nc{\HGK}{H\bs G/K}
\nc{\Dt}{{\DD}_t}
\nc{\Pic}{\on{Pic}_X}
\nc{\Jac}{\on{Jac}_X}
\nc{\Po}{{\mathcal P}}
\nc{\Ofin}{{\mbf O}}
\nc{\KO}{{\K^\times/\Oo^\times}}
\nc{\Kx}{{\K^{\times}}}
\nc{\Ox}{{\Oo^{\times}}}
\nc{\un}{\underline}
\nc{\ti}{\times}
\nc{\ol}{\overline}
\nc{\M}{{\mathcal M}}
\nc{\arr}{\rightarrow}
\nc{\al}{\alpha}
\nc{\pa}{\partial}
\nc{\kk}{{\mathfrak k}}

\nc{\Ap}{\widehat{A}_p}
\nc{\AMOp}{{\mathcal A}{\mathcal M}{\mathcal O}}
\nc{\MOp}{{\mathcal M}{\mathcal Op}}
\nc{\la}{\lambda}
\nc{\Oaf}{{\mbf O}_{p_{-1}}}
\nc{\F}{{\mathcal F}}

\nc{\AOp}{{\mathcal A}{\mathcal O}}
\nc{\Oaff}{{\mbf O}_{\on{aff}}}
\nc{\nat}{\natural}

\nc{\rr}{{\mathfrak r}}

\nc{\mbf}{\mathbf}

\nc{\DD}{\widehat{D}}

\title{Spectral Curves, Opers and Integrable Systems}
\author{David Ben-Zvi}

\address{Department of Mathematics, University of Chicago, Chicago IL
60637}

\author{Edward Frenkel}\thanks{Partially supported by grants from the
Packard Foundation and the NSF}

\address{Department of Mathematics, University of California,
Berkeley, CA 94720}

\begin{abstract}

We establish a general link between integrable systems in algebraic
geometry (expressed as Jacobian flows on spectral curves) and soliton
equations (expressed as evolution equations on flat connections). Our
main result is a natural isomorphism between a moduli space of
spectral data and a moduli space of differential data, each equipped
with an infinite collection of commuting flows. The spectral data are
principal $G$--bundles on an algebraic curve, equipped with an abelian
reduction near one point. The flows on the spectral side come from the
action of a Heisenberg subgroup of the loop group. The differential
data are flat connections known as opers. The flows on the
differential side come from a generalized Drinfeld--Sokolov
hierarchy. Our isomorphism between the two sides provides a geometric
description of the entire phase space of the hierarchy. It extends the
Krichever construction of special algebro--geometric solutions of the
$n$th KdV hierarchy, corresponding to $G=SL_n$.

An interesting feature is the appearance of formal spectral curves,
replacing the projective spectral curves of the classical
approach. The geometry of these (usually singular) curves reflects the
fine structure of loop groups, in particular the detailed
classification of their Cartan subgroups.  To each such curve
corresponds a homogeneous space of the loop group and a
soliton system. Moreover the flows of the system have interpretations
in terms of Jacobians of formal curves.

\end{abstract}

\maketitle

\tableofcontents
\section{Introduction.}
\subsection{Background.}

The Korteweg--deVries hierarchy is an infinite family of commuting
flows on the space of second--order differential operators
$L=\partial_t^2+q$ in one variable. It has long been known that this
hierarchy has close ties to the geometry of algebraic
curves. The Krichever construction explains how to obtain such an operator
$L$ from a line bundle $\Ll$ on a hyperelliptic curve $Y$, equipped
with some local data near a point $\infty\in Y$.  Changing $\Ll$ by
the action of the Jacobian of $Y$ changes $L$ by the KdV flows.  This
picture was extended to the $n$--th KdV hierarchy, in which the
second--order operator $L$ is replaced by an $n$--th order operator.
By replacing hyperelliptic curves with $n$--fold branched coverings
$Y$ of $\Pp$, one finds a relation between Jacobian flows on line
bundles on $Y$ and KdV flows on associated differential operators.
The resulting special ``algebro--geometric'' solutions to KdV may be
understood in great detail.

This picture illustrates a general phenomenon: an integrable system,
naturally expressed in terms of differential data (differential
operators, flat connections etc.), may be characterized using spectral
data, on which the flows become linear and which have
group--theoretic, and sometimes geometric, significance.  However,
most differential operators $L$ do not arise from the geometry of
curves in this way. Instead, the full phase space may be described
using the beautiful algebraic formalism of the Sato Grassmannian and
pseudodifferential operators (see \cite{DJKM,SW,M1,M3}).

Drinfeld and Sokolov \cite{DS} generalized the differential side of
KdV by replacing $n$--th order differential operators by connections
on rank $n$ vector bundles, that is by translating from $n$--th order
equations to first--order systems. This enabled them to associate a
generalized KdV hierarchy to an arbitrary semisimple Lie group $G$.
These hierarchies live on spaces of connections on principal
$G$--bundles on the line.  Recently, these connections were given new
importance, a coordinate--free formulation, and the name ``opers'' by
Beilinson and Drinfeld \cite{opers} in the course of their work on the
geometric Langlands correspondence \cite{Hecke}.  An oper is a
$G$--bundle on a complex curve with a (flat) holomorphic connection
and a flag, which is not flat but obeys a strict form of Griffiths
transversality with respect to the connection. 

Drinfeld and Sokolov also incorporated a spectral parameter into their
connections -- thereby providing a generalization of the eigenvalue
problem for a differential operator.  They showed that the resulting
connections (a loop group version of opers, which we name affine
opers) may be brought into a canonical gauge, where an
infinite--dimensional abelian group of symmetries becomes apparent.
Using these symmetries, it is easy to write a commuting hierarchy of
flows as ``zero--curvature equations'' -- constraints expressing the
flatness of the connection, when extended to new variables using the
symmetries.  Algebraic generalizations of the Drinfeld--Sokolov
hierarchies have been introduced (see \cite{dGHM,Feher}) in which the
abelian group underlying the Drinfeld--Sokolov equations is replaced
by more general abelian subgroups of loop groups.

The spectral side of KdV has also been greatly developed and
generalized (see \cite{M2,AB,DM,LM}). Classically, one studies line
bundles $\Ll$ on a curve $Y$ which is an $n$--fold branched cover of
$\Pp$, or more generally of some projective curve $X$. Taking the
pushforward of $\Ll$ down to $X$ produces a rank $n$ vector bundle
$\E$. Away from the branch points, the bundle $\E$ decomposes into a
direct sum of lines, while at the branch points this decomposition
degenerates, producing a flag. This additional structure makes $\E$
into a Higgs bundle. Conversely, from this Higgs data on $\E$ we may
recover the ``spectral curve'' $Y$ and the line bundle $\Ll$ on $Y$.
(Usually one defines Higgs fields on $\E$ as one--form valued
endomorphisms of $\E$.  The decomposition of $\E$ is then achieved by
considering the eigenspaces of the endomorphism, and the spectral
curve $Y$ parameterizing the eigenvalues is naturally embedded in
$T^*X$. We will only be interested in ``abstract'' Higgs fields, where
we retain the decomposition structure on $E$ but forget the
endomorphism which induced it.)

By reformulating decompositions into lines as reductions to maximal
tori, one can extend this picture from vector bundles to principal
$G$--bundles, following the general formalism developed by Donagi
\cite{Donagi,DM,DG} (see also \cite{Fa}). One considers reductions of
a $G$--bundle $\E$ to a family of Cartan subgroups of $G$, which is
allowed to degenerate at certain points. This gives rise to the
definition of a (regular) principal Higgs field as a sub-bundle of
regular centralizers in the adjoint bundle of a $G$--bundle. Moduli
spaces of principal Higgs bundles provide natural models for
completely integrable systems in algebraic geometry. This was probably
first realized by Hitchin \cite{Hi}. Similar ideas have been used by
Cherednik \cite{Ch1,Ch2,Ch3} in his study of algebro--geometric
solutions of generalized soliton hierarchies. However, apparently no
attempts have been made to {\em identify} the spectral and
differential sides of soliton equations.

\subsection{The Present Work.}
In the papers \cite{FF0,FF1,NLS,EF1,EF2} a new approach to the study
of KdV equations was introduced by Feigin, Enriquez and one of the
authors (see \cite{Five} for an overview).  This approach is based on
the study of certain homogeneous spaces for (subgroups of) the loop
groups (which also arose in \cite{W} from a different point of
view). These spaces come with an obvious action of an
infinite--dimensional abelian group, and carry simple systems of
coordinates in which the flows are easily understood. Using these
coordinates, it is shown that these spaces are isomorphic to the KdV
phase spaces and that the flows agree with the Drinfeld--Sokolov
hierarchies.

Our original motivation for the current work was to understand
geometrically how opers arise from homogeneous spaces for loop groups.
This involved three main steps:
\begin{enumerate}
\item[(1)] Identifying a moduli space interpretation for the homogeneous
spaces (in particular thinking of them as schemes rather than as
sets).
\item[(2)] Finding a natural morphism between the moduli spaces describing
the spectral side and the differential side (opers), explaining the
explicit construction of \cite{FF1}.
\item[(3)] Establishing an intrinsic reason for this morphism to be an
isomorphism.
\end{enumerate}

It is well known that homogeneous spaces for loop groups correspond to
moduli spaces of bundles on a curve, with some extra structure. In the
present case, we found this extra structure to be a formal
generalization of the Krichever data.  Specifically, the relevant
moduli space is the ``abelianized Grassmannian'' $\Grp_A$. This is the
moduli space of $G$--bundles on $\Pp$, equipped with a spectral curve
description on the formal neighborhood $\D$ of a point $\infty \in
\Pp$. This spectral datum may be formulated as a principal Higgs field
on $\D$ with prescribed branching, or as a reduction of the structure
group of the bundle to a twisted family of Cartan subgroups of
$G$. (In the case of $SL_n$, the appropriate spectral curve is an
$n$--fold cover of $\D$, fully branched over $\infty$.)

Suppose we are given such a $G$--bundle $\E$ on $\Pp$ with a Higgs
field on $\D$. If the Higgs field extends to all of $\Pp$, it will do
so uniquely. Thus our moduli space contains a distinguished subspace,
consisting of Higgs bundles on $\Pp$ whose spectral curve is a global
$n$--fold branched cover $Y$ -- in other words a classical spectral
curve for KdV. Thus we have embedded the Krichever data into a much
bigger space of formal spectral data.  This explains why the
homogeneous spaces in \cite{FF1} have a chance to be isomorphic to the
{\em entire} phase space of KdV hierarchy (the space of all opers on
the disc), while the global spectral data only recover special
``finite--gap'' differential operators. In particular, these formal
spectral curves serve as an algebro--geometric substitute for the
analytic theory of infinite--genus spectral curves, \cite{McKean}.

The most important aspect of the abelianized Grassmannians $\Grp_A$ is
that they come with a canonical action of an infinite--dimensional
abelian Lie algebra. The formal group $\Ahat$ of this Lie algebra can
be interpreted as the formal Jacobian (or Prym) variety of the formal
spectral curve. Its action is the natural generalization (and
extension) of the Jacobian flows appearing in the Krichever
construction.

Now that both sides of the isomorphism from \cite{FF1,Five} have
geometric interpretations, corresponding to the differential and
spectral sides of KdV, the second step is a natural construction of
flat connections from moduli spaces of bundles. In the Krichever
setting of line bundles on global branched covers, such a construction
was explained in the classic works of Drinfeld and Mumford
\cite{Dr,Mum}.  We discovered that the calculations in \cite{Five} can
be interpreted as a generalization of this idea, where we work not
with line bundles on a projective curve, but with $G$--bundles on
$\Pp$ with formal spectral data. The main idea is that the flows on
the moduli space $\Grp_A$ lift to tautological bundles, and that this
lifting leads to the construction of flat connections. These
connections naturally live on the formal group $\Ahat$ itself --
i.e. on the formal Prym variety. Their restrictions to distinguished
one-parameter formal subgroups are identified as opers. The
zero--curvature equations relating the different one--parameter flows
(i.e. the condition of the flatness of the connection on $\Ahat$)
translate precisely into the zero--curvature formulation of the KdV
hierarchy.

We have thus found a natural morphism from $\Grp_A$ to the space
of opers, under which the action of $\Ahat$ on $\Grp_A$ translates
into the KdV flows on the space of opers. The final step is to find
out why this morphism is an isomorphism. As in \cite{DS}, one gains
more insight by replacing opers by their loop group cousins, the
affine opers (in other words, by incorporating the spectral parameter
into the connection). 
Hence we explain how to go back and forth between opers 
and affine opers. Using the Drinfeld--Sokolov gauge for affine opers, we
obtain a simple inverse to our map from $\Grp_A$ to opers, in
particular proving that it is an isomorphism.

Thus, our main result is that when suitably generalized, the Krichever
construction can actually be made into an {\em isomorphism} between a
moduli space of bundles with formal spectral data and the phase space
of a soliton hierarchy. In the abelian setting, the Krichever
construction has been explained by Rothstein in \cite{Ro1,Ro2} and
Nakayashiki \cite{N1,N2} in the language of the (generalized)
Fourier--Mukai transform (\cite{Lau}). Thus our isomorphism should
perhaps be thought of in the context of a non--abelian Fourier
transform.

On closer examination, the construction of connections we use turns
out to be independent of the specifics of the problem, but rather an
application of a very general construction.  The underlying structure
is an isomorphism between any double quotient of an algebraic group
with an appropriate space of flat connections on a subgroup. In fact,
these connections reflect a certain remnant of the connections coming
from the trivial Harish--Chandra structures on homogeneous spaces.
When this construction is applied to the spaces $\Grp_A$, one
naturally obtains affine opers in the Drinfeld--Sokolov gauge. 

We find several interesting contexts in which to apply these
abstractions. Our isomorphism between formal spectral data and
differential data not only specializes to the Krichever construction and
extends it to principal bundles, but
\begin{enumerate}
\item[(1)] We may replace the base curve $\Pp$ by an arbitrary curve.
\item[(2)] We may allow arbitrary monodromy of the spectral curve,
obtaining geometric descriptions of all of the generalized
Drinfeld--Sokolov hierarchies of \cite{dGHM,Feher} (labeled by
conjugacy classes in the Weyl group).
\item[(3)] We may allow arbitrary singularities of the spectral
curve (replacing the smooth spectral data appearing above),
obtaining continuous families of new integrable systems.
\end{enumerate}
In forthcoming work, we apply this approach as follows:
\begin{enumerate}
\item[(1)] The description of (generalized) Drinfeld--Sokolov
  hierarchies as flows on spectral data automatically implies a strong
  compatibility with the Hamiltonian structure of the Hitchin system
  and its meromorphic or formal generalizations (recovering in
  particular results of \cite{DM}). This is closely related to the
  geometry of the affine Springer fibration.

\item[(2)] In the case of the principal Heisenberg $A$, we generalize
the isomorphism of \thmref{tripartite theorem},(3) between the open
subspace of the moduli space $\Grp_A$ and the moduli space of (ordinary)
opers on the disc to the case of an arbitrary curve $X$. Namely, we
obtain a map from a subspace of $\Gr_A$ corresponding to a ``generic''
$G$--bundle on $X$ to the space of opers on the disc
.
\item[(3)] We extend the ideas of this paper to $GL_{\infty}$,
  providing a similar point of view on the KP hierarchy and
  pseudo--differential operators (allowing more ``natural''
  modifications of the Krichever construction in the case of line
  bundles.)  
\end{enumerate}

Some areas for future work include understanding the behavior of the
exotic new integrable systems and their relations with the geometry of
singular spectral curves; the interpretation of tau functions as theta
functions for formal Jacobians, and its application to explicit
formulas for solutions; relations with conformal blocks, vertex
algebras and Virasoro actions; identifying the ``spectral'' meaning of
the Gelfand--Dickey hamiltonian structure; and analogs where we
replace differential operators by $q$--difference operators or
polynomials in Frobenius, relating to the $q$--KdV equations and
elliptic sheaves respectively.

\subsection{Summary of Contents.}

The paper proceeds in the opposite direction from the introduction: we
start with the most general notions, and step by step specialize them,
until we end with the calculations which motivated the work. This
simplifies the exposition because the proofs become elementary in the
appropriate light, and we hope this will help clarify the underlying
ideas. We refer the reader desiring a more concrete and explicit
picture to our descriptions in the most important case (the principal
Heisenberg algebra and the usual KdV hierarchy) and to the survey
\cite{Five} for the origins of our approach.

In \secref{general} we explain a general group--theoretic construction
of isomorphisms between moduli spaces of bundles and moduli spaces of
connections, which is responsible for the spectral--differential
equivalence for KdV.  The connections arise from pulling back
equivariant bundles on a space with a group action, to the group
itself. When the space is a double quotient and the bundle
tautological, one easily characterizes precisely which connections are
obtained. Roughly speaking, we identify double quotients $\HGK$ with
moduli of certain connections on the normalizer of $K$.  This
characterization is phrased in terms of the relative position of a
reduction of a bundle with respect to a connection, a notion we
describe in \secref{types of connections}. The ideas behind this are
that of a period map and the localization for Harish--Chandra
pairs. We also present a formulation in terms of differential schemes,
which is closer in spirit to the theory of elliptic sheaves in
characteristic $p$.

In \secref{Krichever}, we summarize the Krichever construction, which
relates line bundles on an algebraic curve with differential operators
in one variable. The exposition is inspired by \cite{Dr,Mum} and
\cite{Ro1,Ro2}, and informed by our general approach.

In \secref{moduli spaces} we introduce loop groups and some of their
homogeneous and double quotient spaces, which are interpreted in a
standard way as moduli spaces of bundles on a curve. This gives a
context in which to apply the general constructions. To find
interesting connections, however, we need interesting group actions.
So in \secref{Heisenbergs}, we study the Heisenberg subgroups of the
loop group at some length.
In particular, we describe their fine classification 
and explain their relation with the geometry of formal
spectral curves. This section may be read
independently of the remainder of the paper.

This leads us in \secref{Higgs} to the study of the main objects of
interest, the abelianized Grassmannians, the moduli spaces of bundles
equipped with a reduction to a Heisenberg subgroup $A$. Alternatively,
they can be described as the moduli of $G$--bundles with formal Higgs
field and fixed spectral curve. The abelianized Grassmannians come
with the action of an infinite--dimensional abelian formal group
$\Ahat$, which is naturally interpreted as the Jacobian (or Prym)
variety of the formal spectral curve. We apply the abstract
construction of \secref{general} to obtain an isomorphism between the
abelianized Grassmannian and a certain moduli space of flat bundles.

These flat bundles can be recast in the more tangible form of affine
opers, discussed in \secref{opers}. Affine opers are $LG$--bundles
with a flat connection and a reduction having a distinguished relative
position. The concept of affine opers is modeled on that of $G$--opers
introduced by Beilinson and Drinfeld \cite{opers} following
\cite{DS}. For classical groups $G$--opers are identified with special
differential operators. Although the concepts of $G$--opers and affine
opers turn out to be essentially equivalent, the latter is more
suitable in our context. The most important property of affine opers
is a canonical abelian structure identified by the Drinfeld--Sokolov
gauge (\propref{DS gauge}).  The notion of affine opers and the
Drinfeld--Sokolov gauge may be extended to arbitrary Heisenbergs with
good regular elements.

Our main results are presented in \secref{Main}. There we use the
Drinfeld-Sokolov gauge to establish an isomorphism between the moduli
of affine opers and the abelianized Grassmannians, thus establishing a
general differential--spectral equivalence for a wide range of
integrable systems.

{}From the point of view of the theory of integrable systems, our main
result in \secref{Main} is a natural and coordinate independent
construction of an integrable hierarchy of flows on the appropriate
space of affine opers, associated to an arbitrary Heisenberg
subalgebra of the loop algebra $L\g$, and a strongly regular
element. In the special case when the Heisenberg is smooth we recover
the generalized Drinfeld-Sokolov hierarchies introduced in
\cite{dGHM,Feher}.

\subsection{Schemes, Stacks, etc.}
This paper is concerned with moduli spaces of bundles and of flat
connections, from the viewpoint of algebraic geometry. Since these
``spaces'' are rarely varieties, this necessitates the use of some
less familiar objects, namely algebraic stacks and ind--schemes. We
refer the reader to \cite{BL,LS,Tel} for a detailed description of
moduli spaces in this language, and to \cite{LMB,Sor} for a general
treatment of stacks. Our main results of interest to experts in
integrable systems are formulated in \secref{Main} in terms of
varieties, differential polynomials and evolutionary derivations.  All
of our stacks will be algebraic (in the Artin sense, so that
automorphism groups of an object may be infinite). We will sometimes
abuse notations and write $x\in M$ for a stack $M$, signifying an
$S$--point of $M$ for some scheme $S$.

All schemes, groups, sheaves, representations etc. will be defined
over $\C$.

Throughout the paper we will refer to group schemes, for which the
underlying scheme structure is obvious, simply as groups. By a
$G$--torsor over a scheme $X$ we will understand a scheme $E \to X$
equipped with a right action of $G$, such that locally in the flat
topology it is isomorphic to $G \times X$. The term $G$--torsor is of
course synonymous to the term principal $G$--bundle.

An ind--scheme is by definition an inductive limit of schemes (in the
category of {\em spaces}, namely sheaves of sets in the {\it fppf}
topology).  Ind--schemes can be very pathological in general; however,
the ones we will encounter owe their inductive nature primarily to
infinite--dimensionality. In particular, they are unions of closed
subschemes and are formally smooth. We will refer to group ind--schemes
simply as ind--groups.

The typical example of an ind--scheme in our setting is the loop group
${\mathcal L}{\mathcal G}$ of an algebraic group $G$, whose
$R$--points are the points of $G$ over the formal Laurent power series
$R((t))$.  It is an inductive limit of schemes corresponding to
Laurent series with bounded poles.  Perhaps a more familiar class of
ind--schemes is that of formal groups, which are group ind--schemes
having points over rings with nilpotents, but no non--trivial points
over any field.  They arise from formally exponentiating the action of
a Lie algebra.

\section{A Construction of Connections.}\label{general}

In this section we describe a general construction, which allows us to
identify double quotient spaces with moduli spaces of connections. In
our applications, the double quotient will be a moduli space of
bundles, while the connections will be opers and their
generalizations, which appear in soliton theory.  We first describe
the notion of {\em relative position} for a connection and a reduction
of a bundle, and its relation to period maps. We then characterize the
connections arising on homogeneous spaces from the theory of
Harish--Chandra pairs \cite{BB}.  In this way our construction is
related to the localization theory of representations. We identify an
aspect of this picture which can be generalized to double quotient
spaces, namely we obtain a map from the double quotient to a space of
connections of a particular type. Finally using period maps we show
that this map is an isomorphism (\propref{abstract isomorphism}): the
double quotient space itself classifies all connections of the
prescribed type,

\subsection{Types of Connections.}\label{types of connections}
Let $G$ be a group scheme and $K\subset G$ a subgroup, with Lie
algebras $\g$ and $\kk$, respectively.  
Given a $G$--torsor $\E$ on a scheme $X$, and a scheme
$M$ equipped with an action of $G$ (e.g., a
representation of $G$), we define the {\em $\E$--twist of $M$} as $\E
\times_G M$, and denote it by $(M)_{\E}$. This is a bundle over $X$,
whose fibers are isomorphic to $M$.

Suppose $\E$ is a $G$--torsor on a smooth scheme $X$, with a
connection $\nabla$ and a reduction $\E_K$ to $K$.  Then we may
describe the failure of $\nabla$ to preserve $\E_K$ in terms of a
one-form with values in $(\g/\kk)_{\E_K}$. Locally, choose any flat
connection $\nabla'$ on $\E$ preserving $\E_K$, and take the
difference $\nabla' - \nabla$.

\subsubsection{Lemma.} The local $(\g/\kk)_{\E_K}$--valued one-forms 
$[\nabla'-\nabla]$ are independent of $\nabla'$. They 
define a global section $\nabla/\E_K$ of $(\g/\kk)_{\E_K}\otimes\Omega^1$. 

\subsubsection{Remark.} 
More abstractly, this construction can be phrased as follows.  Let
$\mathcal A_{\E}$ be the Lie algebroid of infinitesimal symmetries of
$\E$, and $\mathcal A_{\E_K}$ be the subalgebroid preserving
$\E_K$. Thus $\nabla$ is a splitting of the anchor map of ${\mathcal
A}_{\E}$ to the tangent sheaf, and $\nabla'$ does the same for
${\mathcal A}_{\E_K}$. The difference $\nabla'-\nabla$ is therefore a
map from the tangent sheaf to the quotient ${\mathcal
A}_{\E}/{\mathcal A}_{E_K}$. But the latter sheaf is canonically
isomorphic to $(\g/\kk)_{\E_K}$. Thus we obtain a
$(\g/\kk)_{\E_K}$--valued one-form, which is $\nabla/\E_K$.

\subsubsection{} We can also view $\nabla/\E_K$ as a
map from the tangent bundle $TX$ of $X$ to $(\g/\kk)_{\E_K}$ (denoted
in the same way).
This map can be realized as follows.  Consider the
$G/K$--bundle of $\E$, $(G/K)_{\E}$. The $K$--reduction $\E_K$ gives a
section $s:X\to (G/K)_{\E}$, with differential
$ds: TX \to T((G/K)_{\E})$.
The connection $\nabla$ gives rise
to a horizontal subbundle of the tangent bundle of $(G/K)_{\E}$,
and hence a canonical projection $p_{\nabla}:T((G/K)_{\E}) \to
(\g/\kk)_{\E_K}$. The following simple fact will be useful in the
proof of \propref{abstract isomorphism} below.

\subsubsection{Lemma.}    \label{differential}
The map $\nabla/\E_K:TX\to (\g/\kk)_{\E_K}$ coincides with the
composition 
$$TX \stackrel{\longrightarrow}{ds} T((G/K)_{\E})
\stackrel{\longrightarrow}{p_\nabla}
(\g/\kk)_{\E_K}.$$

\subsubsection{Definition: Relative Position.} \label{relative position}
Let $\mbf O$ be an orbit for the adjoint action of $K$ on $\g/\kk$,
and $\xi$ a vector field on $X$.  The reduction $\E_K$ is said to
have {\em relative position $\mbf O$} with respect to $\nabla_\xi$ if
the image of $\xi$ under the map $\nabla/\E_K: TX \arr (\g/\kk)_{\E_K}$
takes values in $({\mbf O})_{\E_K} \subset (\g/\kk)_{\E_K}$.

\subsubsection{Period Maps.}\label{period map}
The action of $K$ on the quotient $\g/\kk$ may be identified with its
action, as the stabilizer subgroup of the identity coset $[1]\in G/K$,
on the tangent space at that point. Thus there is a one-to-one
correspondence between the $K$--orbits on $\g/\kk$ and the $G$--orbits
in the tangent bundle $T(G/K)$ to $G/K$ ($G$--invariant distributions
on $G/K$).  This leads to a period map interpretation of relative
position. Fix a point $x\in X$ and a trivialization $\E_x\cong G$ of
the fiber at $x$. Then the connection defines a canonical
trivialization of $\E$ on the formal completion of $X$ at $x$.
(Complex analytically, we obtain a trivialization on any simply
connected neighborhood of $x$.)  The section $s$ then provides a map,
the formal period map, from the formal completion of $X$ at $x$ to
$G/K$. The reduction has relative position $\mbf O$ if the vector
field $\xi$ is tangent to the $G$--invariant distribution on $G/K$
corresponding to $\mbf O$.

\subsubsection{Remark.}    \label{trivial}
It is useful to observe that if the bundle $\E$ has a global flat
trivialization, then we can define a global period map $X\to
G/K$. According to \lemref{differential}, the relative position map
$\nabla/\E_K$ is simply its differential.

\subsubsection{Remark.} When the orbit 
$\mbf O$ is $\C^\times$--invariant, we can apply the definition of
relative position to
all vector fields simultaneously. Namely, we say that {\em $\E_K$ has
relative position $\mbf O$ with respect to $\nabla$} if the one-form
$\nabla/\E_K$ takes values in $({\mbf O})_{\E_K} \otimes_{\C^\times}
\Omega^1_X$.

\subsubsection{Difference analog.}\label{difference}
It is instructive to compare the above notion of relative position
with the notion of relative position for equivariant bundles. Suppose
a group $A$ acts on $X$, and let $\E$ be an $A$--equivariant
$G$--torsor on $X$.  Suppose furthermore that $\E$ is equipped with a
reduction $\E_K$ to $K$. The group analogs of $K$--orbits in $\g/\kk$
are $K$--orbits in $G/K$. These correspond bijectively to diagonal
$G$--orbits on $G/K\times G/K$.

Let $a\in A$ and let ${\mbf O}\subset G/K$ be a $K$--orbit. Then we
can say $\E_K$ has relative position $\mbf O$ with respect to $a$ if
$a^*\E_K\subset {\mbf O}_{\E_K}\subset (G/K)_{\E_K}$. The notion of
period map also carries over. Assume for simplicity that we may
$A$--equivariantly trivialize the bundle.  The relative position map
then sends $X\times A\to G/K\times G/K$ via $(x,a)\mapsto
(\E_K|_x,\E_K|_{a\cdot x})\in G/K\times G/K$. Thus $\E_K$ has relative
position $\mbf O$ with respect to $a$ if the period map sends
$X\times\{a\}$ to the diagonal $G$--orbit corresponding to $a$.

This group--theoretic notion extends automatically to the case when
$A$ is an ind--group (or sheaf of groups)\footnote{Replacing groups by
groupoids, we may take $A$ to be the formal neighborhood of the
diagonal and recover the case of flat connections discussed above.}. It
also allows one to define relative position for difference operators
on $G$--torsors, as well as for $G$--torsors on a scheme over a finite
field which are equivariant with respect to Frobenius.  We will
briefly return to this idea in \secref{extended family}.

\subsection{Harish--Chandra structures.} 
We present here an elementary aspect of the theory of $(\g,K)$ (or
Harish--Chandra) structures (see \cite{BB}).  Let $G$ be an
(ind--)group, with Lie algebra $\g$, and $K\subset G$ a subgroup. A
$(\g,K)$--structure on a scheme $M$ is a $K$--torsor $\Po$ over $M$,
together with an action of $\g$ on $\Po$. The restriction of the
$\g$--action on $\Po$ to the Lie subalgebra ${\kk}=\on{Lie}K\subset\g$
is assumed to coincide with the action of the latter coming from the
action of $K$ on $\Po$. Moreover, the action of $\g$ is assumed to be
simply transitive, so that the natural map from $\g$ to the tangent
space to $\Po$ at any point is an isomorphism. It follows then that
the tangent bundle of $M$ is identified with the $\Po$--twist
$(\g/{\kk})_{\Po}$ of $\g/{\kk}$.  The basic example of a
$(\g,K)$--scheme is $G/K$ itself, with $\g$ acting {\em from the
right} on the total space of the $K$--bundle $\Po=G\to G/K$.

\subsubsection{Lemma.}\label{Harish-Chandra}
The $G$--torsor $\Po_G=\Po\times_K G$ induced from $\Po$ carries a
canonical flat connection $\nabla$, such that the map $\nabla/\Po:
TM\cong (\g/\kk)_{\Po}$ is the isomorphism
induced by the $\g$--action on $\Po$.

\subsubsection{Proof for $M=G/K$.}
The $G$--torsor $\Po_G=G\times_K G$ is canonically identified with
$G/K\times G$, by the map $(g_1,_K g_2)\mapsto (g_1K,g_1g_2)$.  To a
point $g_1K\in G/K$ we assign the $K$--reduction $(g_1,_K 1)\in
G\times_K G$, which corresponds to the subset $(g_1K,g_1K)\subset
G/K\times G$. Thus we obtain a global trivialization of $\Po_G$, hence
a flat connection. By construction, the period map $G/K\to G/K$
defined using this trivialization is the identity map, and the lemma
follows.

\subsubsection{General Case.} 
Let $(\g,K)^{{\wedge}}$ denote the group ind--scheme generated by the
group $K$ and the formal group $\g^{{\wedge}}$ of $\g$, in other words
the formal completion of $G$ along $K$. Then for any $x$ in $M$, the
completion of $\Po$ along the fiber $\Po_x$ is a
$(\g,K)^{{\wedge}}$--torsor. The resulting principal bundle $\Po^{\wedge}$
on $M$ carries a flat connection, since it is canonically trivialized
over any local Artinian subscheme of $M$: 
the formal neighborhoods of
infinitesimally nearby fibers of $\Po$ are the {\em same}. 
The induced $G$--torsor, which coincides with $\Po_G$,
inherits this flat connection $\nabla$ as well as the $K$--reduction
$\Po$.  By construction, the connection identifies the bundle
$(\g/\kk)_{\Po_G}=(\g/\kk)_{\Po^{{\wedge}}}$ with the tangent bundle
of $M$.

\subsection{Double Quotients.}\label{double cosets}
Now let $G$ be a group scheme, and $H,K$ subgroup--schemes. Let us
consider the double quotient stack $H\bs G/K$. This means the following:
to each scheme $S$ we
attach the groupoid $H\bs G/K(S)$, whose objects are $G$--torsors on
$S$ together with reductions to $K$ and to $H$. The morphisms are the
isomorphisms of such triples.  The definition of the functor $(H\bs
G/K)(S_2) \arr (H\bs G/K)(S_1)$ corresponding to a morphism $S_1 \arr
S_2$ is straightforward.

For example, when $K$ and $H$ are both equal to the identity subgroup,
the objects of the groupoid $(1\bs G/1)(S)$ are $G$--torsors endowed
with two reductions to the identity, hence two sections. We may use
the first section to trivialize the torsor, and the other section
gives us a map from $S$ to $G$, i.e., an $S$--point of $G$. Therefore
$1\bs G/1 = G$. For general $K$ and $H = \{ 1 \}$, we obtain the
equivalence between reductions of the trivial $G$--torsor on $S$ to
$K$ and maps from the base $S$ to $G/K$. In general we have a
surjective morphism from the scheme $G/K$ to $\HGK$, realizing $\HGK$
as an algebraic stack.

The stack $\HGK$ carries a tautological $G$--torsor $\T$. Its
fiber over an $S$--point of $\HGK$, thought of as a $G$--torsor
$\mathcal P$ on $S$, is identified with $\mathcal P$. Moreover $\T$
comes equipped with tautological reductions $\T_K$ and $\T_H$ to $K$
and $H$, respectively. Explicitly, $\T = G \times_{H}
G/K = H \bs G \times_{K} G$, $\T_H = G/K$, $\T_K =  H\bs
G$.

\subsubsection{Connections on Double Quotients.}
When $H\subset G$ is the identity subgroup, the space $G/K$ has an
obvious $(\g,K)$--structure (\secref{Harish-Chandra}). Thus the
tautological $G$--bundle $\T$ carries a canonical flat connection,
which has a ``tautological'' relative position with respect to $K$.
For general $H$ the construction of \secref{Harish-Chandra} breaks
down, since the action of $\g$ on the $K$--torsor $\Po=H\bs G$ is no
longer simply transitive. Accordingly, there is no natural flat
connection on $\Po$. A well-know way to circumvent this problem is to
replace the flat vector bundles associated to our flat $G$--torsors by
$\mathcal D$--modules, obtained by taking coinvariants by the
stabilizers of the $\g$--action.  Our approach explained below is to
replace an action of all vector fields on $H\bs G$ by those coming
from the action of an appropriate subgroup of $G$, and to construct
connections not on $\HGK$ but on the subgroup itself.
This leads, in \propref{abstract isomorphism}, to an identification of $\HGK$
with a moduli stack of special connections on a subgroup.

\subsubsection{Actions give connections.}\label{action->connection} 
Let $A$ be a group--scheme, and $M$ a scheme equipped with an
$A$--action. Suppose $\T$ is an $A$--equivariant $G$--torsor on $M$,
so in particular we may lift the action of the Lie algebra $\ab$ of
$A$ from $M$ to $\T$. If the $A$ action is not free, it does not
follow that $\T$ obtains a partial connection along the
$A$--orbits. Namely, the action of the stabilizers in $\ab$ on $\T$
presents an obstruction for lifting vector fields consistently to the
bundle.

However, for any $x\in M$, with $A$--orbit $\pi_x:A\to M$, the
$A$--action naturally identifies the pullback bundle $\pi_x^*\T$ on
$A$ with the trivial bundle $A\times \T_x$. Therefore $\pi_x^*\T$ has
a canonical flat connection (albeit isomorphic to a trivial
connection).

\subsubsection{}
Now we apply the construction of \secref{action->connection} to
$M=\HGK$ and $\T$. First we need to identify natural group actions on
$\HGK$. Let $A$ be an arbitrary subgroup of $N(K)$, the normalizer of
$K$ in $G$. Then the right action of $A$ on $H\bs G$ descends to
$\HGK$. In fact, if $A_+=A\cap N(K)$, then the quotient group $A/A_+$
acts on $\HGK$. Furthermore, since $A$ acts on $G/K=\T_H$, we have the
following obvious

\subsubsection{Lemma.}\label{lifting} 
The $A/A_+$--action on $\HGK$ lifts canonically to
$\T$, preserving the reduction $\T_H$ to $H$.

\subsubsection{}    
Therefore the bundle $\T$ over $\M = \HGK$ is $A/A_+$--equivariant.
For any $x\in\HGK$, the construction of \secref{action->connection}
results in a $G$--torsor $\E^x$ on $A/A_+$ with a flat connection
$\nabla$ (induced by a trivialization). The $G$--torsor $\E^x$ also
carries reductions $\E_H^x, \E_K^x$ to $K, H$. Since $\T_H$ is
preserved by the $\Ahat$--action, $\E_H^x$ is automatically flat with
respect to $\nabla$.

The behavior of $\E_K^x$ with respect to $\nabla$ mirrors
\lemref{Harish-Chandra} -- the connection $\nabla$ is simply a part of
the structure of \lemref{Harish-Chandra} which descends to $\HGK$.
Since $\ab$ normalizes $\kk$, the action of $K$ on the Lie algebra
$\ab/\ab_+$ of $\Ahat$ is trivial. Thus for every $a\in\ab$ the
$K$--orbit ${\bf O}_a$ of $a\on{mod}\kk$ in $\g/\kk$ is a point. This
leads us to the following definition:

\subsubsection{Definition.}\label{taut} Let $\E$ be a $G$--torsor with a
$K$--reduction $\E_K$ on $A/A_+$. Then $\E_K$ has {\em tautological
relative position} with respect to a connection $\nabla$ if the image
of the vector field $\xi_a$ coming from the left action of $a$ on
$A/A_+$ under the map $\nabla/\E_K: T(A/A_+) \arr (\g/\kk)_{\E_K}$ is
in ${\bf O}_{-a}$.

\subsubsection{} 
The flat connection on $\E^x$ was obtained from an identification
$\E^x\cong A/A_+\times \E^x|_1$. The additional choice of an
identification $\E^x|_1\cong G$ trivializes the $G$--bundle $\E^x\cong
A/A_+\times G$. We now recall from \remref{trivial} that this global
flat trivialization allows us to define a global period map from
$A/A_+$ to $G/K$, whose differential is the relative position map
$\nabla/{\E^x_K}$.

\subsubsection{Lemma.}\label{tautological position} 
The $K$--reduction $\E^x_K$ is in
tautological relative position with $\nabla^x$.

\subsubsection{Proof.} 
We choose a trivialization $\E_H^x|_1\cong H$, inducing $\E^x|_1\cong
G$ as above.  The resulting trivialization of $\E^x$ preserves the
$H$--reduction $\E^x_H$. This trivialization of $\E_H^x$ gives rise to
a lift of the $A/A_+$--orbit of $x$ on $\HGK$ to an $A/A_+$--orbit on
$G/K$. The period map $A/A_+\to G/K$ induced by the trivialization is
precisely this orbit map. It follows that the relative position of
$\nabla$ is given by the right action of $A/A_+$ on $G/K$, and hence
is tautological.

\subsubsection{Difference version.} The bundle $\E^x$ is 
$A/A_+$--equivariant by construction. Therefore it is natural to
replace the infinitesimal relative position above by its group analog,
\secref{difference}. For $a\in A/A_+$, the $K$--double coset
$KaK=K1Ka$ is a single point. Thus \defref{taut} has an obvious
version, with $a^{-1}$ replacing $-a$. The proof of
\lemref{tautological position} carries over as well.

\subsubsection{}

Denote by $\M_A^{\nabla}$ the stack classifying quadruples
$(\E,\tau,\E_H,\E_K)$, where $\E$ is a $G$--torsor on $A/A_+$, $\tau$
is a trivialization of $\E$, i.e. an identification of $\E$ with
$A/A_+\times \E_1$, where $\E_1$ is the fiber of $\E$ at $1 \in A/A_+$
(this trivialization induces a flat connection $\nabla$ on $\E$),
$\E_H$ is a flat $H$--reduction, and $\E_K$ is a $K$--reduction in
tautological relative position with $\nabla$.  When $A/A_+$ is not
connected, we will always automatically replace the infinitesimal
formulation by its group (difference) version, as above.  When $A/A_+$
is connected, the two are equivalent.

The statement of \lemref{tautological position} holds over any base
$S$, and hence we obtain a natural morphism of stacks $\phi: \HGK \arr
\M_A^{\nabla}$.

\subsubsection{Proposition.}\label{abstract isomorphism}
The morphism $\phi: \HGK \arr
\M_A^{\nabla}$ is an isomorphism of stacks.

\subsubsection{Proof.} There is an obvious forgetful morphism $\psi:
\M_A^\nabla \arr \HGK$, sending $\E\in \M_A^{\nabla}$ to the fiber
$\E|_1$ at the identity of $A/A_+$, considered as a $G$--torsor with
reductions to $K$ and $H$. It is clear that $\psi \circ \phi =
\on{Id}$. It remains to show that $\phi \circ \psi = \on{Id}$.

Given $(\E,\tau,\E_H,\E_K) \in \M_A^\nabla(S)$, where $S$ is an
arbitrary base, we obtain a map $\pi: A/A_+\times S\to\HGK$ classifying
the triple $(\E,\E_H,\E_K)$. Locally on $S'\to S$ (an {\it fppf}
covering), we may trivialize the $H$--torsor $\E_H|_1$, and thus the
$G$--torsor $\E|_1$. Then the map $\tau$ provides a trivialization of
$\E_H$ and $\E$, and hence a lift of $\pi$ to map $\wt{\pi}:A/A_+\times
S'\to G/K$. But this map is precisely the period map, as explained in
\remref{trivial}, and therefore its differential is the relative
position map $\nabla/\E_K$.

Since we know that the relative position of $\E_K$ is tautological, it
follows that the differential of $\wt{\pi}$ coincides with that of the
right $A/A_+$--action on $G/K$.  Hence $\wt{\pi}$ (whence $\pi$) is
$\ab/\ab_+$--equivariant (for $A/A_+$ connected, or
$A/A_+$--equivariant in general). Therefore $\pi(A/A_+\times S')$
equals the $A/A_+$--orbit in $\HGK$ of $\pi(1 \times S')$. This shows
that $\phi \circ \psi = \on{Id}$ and proves the proposition.

\subsubsection{Remark.} 
We note that, due to its difference formulation, the proposition is
applicable in a broader context where we allow $G$ and $A$ to be
ind--groups. In applications below, $A/A_+$ will be an ind--group,
while $K$ will be a group scheme.

\subsection{Alternative Formulation.}\label{differential schemes}
We present a different viewpoint on the above constructions, motivated
by the theory of shtukas in characteristic $p$, and more directly by that
of Krichever sheaves developed by Laumon \cite{Lau2} (see also
\cite{new}). The rough idea is that in order to obtain a
characteristic zero analog of constructions involving Frobenius, one
should consider not schemes $S$ but differential schemes $(S,\pa)$,
where $\pa$ is a distinguished vector field on $S$.\footnote{To obtain
the full parallel of \propref{abstract isomorphism} one simply
replaces differential schemes by schemes with an $A$--action.}

Given a differential scheme $(S,\pa)$ we have the notion of a
differential $G$--torsor $(\E,\pa_{\E})$ on $(S,\pa)$, which is a
$G$--torsor $\E$ on $S$ equipped with an action $\pa_{\E}$ of $\pa$.
There is also a notion of relative position for differential torsors,
following \defref{relative position}: we may require a reduction
$\E_H$ of $\E$ to $H\subset G$ to be in relative position $[-p]$ with
respect to $\pa_{\E}$ (where $p\in\ab/\ab_+, \ab=\on{Lie}(N(H))$).

Let $p\in \ab/\ab_+$ act on $\HGK$ as before. We thus consider the pair
$(\HGK,p)$ as a differential stack.

\subsubsection{Proposition.}\label{differential proof} 
The pair $(\HGK,p)$ represents the functor from differential sche\-mes
to groupoids which assigns to $(S,\pa)$ the category of quadruples
$(\E,\pa_{\E},\E_K,\E_H)$, where $(\E,\pa_{\E})$ is a differential
$G$--torsor on $(S,\pa)$, $\E_K$ is a reduction of $\E$ to $K$
preserved by $\pa_{\E}$, and $\E_H$ is a reduction to $H$ in relative
position $[-p]$ with respect to $\pa_{\E}$ (morphisms being
isomorphisms of such objects).

\subsubsection{Proof.} The proof parallels that of \propref{abstract
isomorphism}. Let $(\E,\pa_{\E},\E_K,\E_H)$ be as above.
Since $\HGK$ classifies triples $(\E,\E_K,\E_H)$ we
obtain a map $S\to\HGK$. The relative position condition implies that
this map gives rise to differential morphism $(S,\pa)\to(\HGK,p)$
classifying $(\E,\pa_{\E},\E_K,\E_H)$ as required. Conversely, given a
differential morphism $(S,\pa)\to(\HGK,p)$, we may pull back the
tautological $G$--torsor $\T$ with reductions to $K$ and $H$, and 
\lemref{tautological position} guarantees that the resulting
quadruple $(\E,\pa_{\E},\E_K,\E_H)$ has the desired properties.

\section{The Abelian Story.}\label{Krichever}

In this section we present the classical construction of Krichever of
algebro--geometric solutions to soliton equations, following the
approach of Drinfeld \cite{Dr} and Mumford \cite{Mum} (see also
Rothstein \cite{Ro1}). We will see that the Krichever construction can
be viewed as a special case of the correspondence between bundles and
flat connections established in the previous section. This is intended
to make the comparison with our generalization in the following
sections more transparent.

\subsection{$GL_n$--opers and differential operators.}\label{GL opers}

\subsubsection{Definition.} A $GL_n$--oper on a smooth curve $Y$ is a
rank $n$ vector bundle $\E$, equipped with a flag
$$0\subset\E_1\subset\cdots\E_{n-1}\subset \E_n=\E,$$ and a connection
$\nabla$, satisfying
\begin{enumerate}
\item[$\bullet$] $\nabla(\E_i) \subset \E_{i+1} \otimes \Omega^1$.
\item[$\bullet$] The induced maps $\E_i/\E_{i-1}\to (\E_{i+1}/\E_i)
\otimes \Omega^1$ are isomorphisms for all $i$.
\end{enumerate}

\subsubsection{} 
In local coordinates a $GL_n$--oper has the form 
$$ \partial_t +\left( \begin{array}{ccccc}
*&*&*&\hdots&*\\
+&*&*&\hdots&*\\
0&+&*&\hdots&*\\
\vdots&\ddots&\ddots&\ddots&\vdots\\ 
0&0&\hdots&+&*
\end{array} \right),$$
where the $*$ are arbitrary and the $+$ are
nonzero.
The oper condition is a strict form
of Griffiths transversality. 

Recall that giving an
$n$--th order differential operator $L$ in one variable
$$\partial_t^n-q_1\partial_t^{n-1}-q_2\partial_t^{n-2}-\cdots-q_n$$ is
equivalent to giving a system of $n$ first--order equations which can
be written in terms of the first--order matrix operator
$$\partial_t-\left( \begin{array}{ccccc}
q_1&q_2&q_3&\cdots&q_n\\
1&0&0&\cdots&0\\
0&1&0&\cdots&0\\
\vdots&\ddots&\ddots&\cdots&\vdots\\
0&0&\cdots&1&0
\end{array}\right).$$
If the $q_i\in \C[[t]]$, then this is a $GL_n$--oper on the formal disc
$\wh{D}=\on{Spf}\C[[t]]$. Conversely it is not hard to see that any
oper may be locally brought into the above form. Thus $GL_n$--opers on
the formal disc are equivalent to $n$--th order differential operators. (A
similar statement holds on global curves if we twist by the appropriate line
bundles.) We thus have

\subsubsection{Lemma.}    \label{opers-to-diff}
$GL_n$--opers on the formal disc $\wh{D}$ are in one--to--one
correspondence with $n$--th order differential operators with
principal symbol $1$.

\subsection{The Krichever Construction.}    \label{alg-geom-setting}
Let $X$ be a smooth, connected, projective curve and $\infty\in X$ a
fixed base point. Denote by $\D$ the ``disc'' around $\infty\in X$,
i.e., $\on{Spec} \Oo$, where $\Oo$ is the completed local ring at
$\infty$.\footnote{Note the difference between the disc $\D =
\on{Spec} \Oo$, which is a scheme, and the formal disc $\wh{D} =
\on{Spf} \Oo$, which is a formal scheme obtained by completing $X$ at
$\infty$. While we can consider a punctured disc $\Dx$, there is no
punctured formal disc.}  If we choose a formal coordinate $z\inv$ on
$\D$ (so that $z$ has a simple pole at $\infty$), we may identify
$\Oo$ with $\C[[z\inv]]$. Let $\Dx$ denote the punctured disc at
$\infty$, i.e., $\on{Spec} \K$, where $\K$ is the field of fractions
of $\Oo$. Choosing a formal coordinate $z\inv$ identifies $\K$ with
$\C((z\inv))$.  However, we note that all of our constructions will be
independent of the choice of formal coordinates.

The field $\K$ has a natural filtration, by orders of poles at
$\infty$: $f\in (\K)_{\geq m}$ if $fz^m\in \Oo$ for any local
coordinate $z\inv$ on $\D$. Thus $\Oo=\K_{\geq 0}$. While the
gradation by order of poles depends on the choice of coordinate $z$,
the filtration is clearly independent of this choice. 
Let $\K^\ti$ denote the group functor of invertible Laurent series: by
definition, the set of $R$--points of $\K^\ti$ is $(R
\widehat{\otimes} \K)^\times \cong R((z\inv))^{\times}$. Note that
$\K^\ti$ is not representable by a scheme, but is a group ind--scheme.
The sub--functor
$\Oo^\ti$ is defined as follows: the set of $R$--points of $\Oo^\ti$
is $(R \widehat{\otimes} \Oo)^\times \cong R[[z\inv]]^{\times}$. This
functor is representable by a group--scheme of infinite type, with Lie
algebra $\Oo$.

The quotient ind--group $\KO$ is isomorphic to a product of $\Z$ and a
formal group. The group of $\C$--points of $\KO$ is naturally
identified with $\Z$. But if $R$ is a ring with nilpotents, then the
group of $R$--points of $\KO$ is much larger: it equals the product of
$\Z$ with the group of all expressions of the form
$$r_{-n}t^{-n}+\cdots+r_{-1}t^{-1}+1$$ where the $r_i$ are
nilpotent. In other words, $\KO$ is isomorphic to the constant group
scheme $\mathbb Z$ times the universal Witt formal group $\wh{\mathbb
W}$ (see \cite{CC}), which is associated with the Lie algebra
$\K/\Oo$.

\subsubsection{}
Let $\Pic$ denote the Picard variety of $X$, i.e., the moduli scheme
of line bundles on $X$ together with a trivialization of the fiber at
a fixed point $0\in X$.  Now consider the moduli scheme $\wt{\Pic}$ of
line bundles $\Ll\in\Pic$ on $X$, together with a trivialization
$\phi$ of $\Ll$ over $\D$. The group $\Oo^\ti$ acts naturally on
$\wt{\Pic}$ by changing trivializations and the quotient
$\wt{\Pic}/\Oo^\ti$ is isomorphic to $\Pic$.

Moreover, the $\Oo^\ti$--action on $\wt{\Pic}$ can be extended to an
action of $\K^\ti$. Informally speaking, given a pair $(\Ll,\phi)
\in\wt{\Pic}$, and an element $k \in \K^\ti$, we define a new line
bundle $\Ll'$ by gluing $\Ll|_{X\sm\infty}$ and $\Oo_{\D}$ over $\Dx$
via $k \phi$; then the bundle $\Ll'$ comes with a natural
trivialization $\phi'$ over $\D$. In other words, we multiply the
transition function of $\Ll$ on $\Dx$ by $k$ (see \cite{LS} for a
discussion of formal gluing of bundles).

%\subsubsection{Proof.}
%The action of $\Oo^\ti$ simply changes the trivialization of $\Ll$ on
%$\D$, and hence does not change $\Ll$. Therefore we obtain an
%injective forgetful morphism $\wt{\Pic}/\Oo^\ti \arr \Pic$. The
%surjectivity of this morphism follows from the fact that for any
%scheme $S$, we can lift any morphism $c:S\to \Pic$ to a morphism
%$\wt{c}: S \to \wt{\Pic}$ locally on $S$.

Since $\K^\ti$ commutes with $\Oo^\ti$, we obtain an action of
$\K^\ti$, and in fact of $\KO$, on $\Pic$. This action is
formally transitive: $\K/\Oo$ surjects onto the tangent space
$H^1(X,\Oo_X)$ to $\Pic$ at any point.  This may be easily seen by
identifying $\K/\Oo\cong H^0(X,i_*\Oo_{X\bs x}/\Oo_X)$, where
$i:X\sm\infty\hookrightarrow X$, and studying the obvious long exact
sequence in cohomology, noting that $H^1(X,i_*\Oo_{X\bs x})=0$ since
$X\sm\infty$ is affine. It follows that we have a surjection from the
connected component of $\KO$ onto the formal group $\wh{\Pic}$ of
$\Pic$ (while the full $\KO$ action changes degrees of bundles as well).

There is a tautological line bundle $\Po$ on $X\times \Pic$ whose
fiber at $x\times\Ll$ is the fiber of $\Ll$ at $x$. The pushforward of
$\Po|_{(X\sm\infty)\times \Pic}$ to $\Pic$ is a quasi-coherent sheaf
$\Po_-=\Po(X\sm\infty)$ on $\Pic$ (its fiber over $\Ll$ is the vector
space of sections of $\Ll$ over $X\sm\infty$).

\subsubsection{Proposition.}\label{abelian lifting} 
The action of $\KO$ on $\Pic$
naturally lifts to $\Po_-$.

\subsubsection{Proof.} This follows immediately from the definition of
the $\K^\ti$ action on $\wt{\Pic}$: in changing the transition
function from $\D$ to $X\sm\infty$ we do not affect the bundle on
$X\sm\infty$. In other words, for $\Ll\in\wt{\Pic}(R)$, there is a
canonical identification between $\Ll|_{X\sm\infty}$ and
$(k\cdot\Ll)|_{X\sm\infty}$, for any $k\in\K^\ti(R)$ and $R$ any
Artinian ring. This identification is $\Oo^\ti$--equivariant, hence
descends to $\Pic$.

%identified with the $\KO$--orbit of $\Ll$. Alternatively we may pull
%$\Po_-$ back to the formal group $\KO$ under the orbit map
%$\KO\to\Pic$ defined by any $\Ll$.

\subsubsection{}
There is a distinguished line $(\K/\Oo)_{\geq -1}$ in the Lie algebra
$\K/\Oo$, consisting of Laurent series with first order pole modulo
regular ones. In a local coordinate $z\inv$ on $\D$, this is the line
$\C z$. Thus we have a distinguished vector field $\pa_1$ on the
Jacobian. The resulting line in the tangent space to the Jacobian
$\Pic^0$ at any point $\Ll$ is naturally identified with the tangent
line to the Abel--Jacobi map based at $\infty$. By \propref{abelian
  lifting}, the vector field $\pa_1$ naturally lifts to the sheaf
$\Po_-$ and provides the latter with a partial connection along this
distinguished direction. Given $\Ll\in \Pic$, we may restrict $\Po_-$
to the formal disc $\wh{D}_t=\on{exp}(t\pa_1)\cdot\Ll$ generated by
$\pa_1$, obtaining a flat vector bundle.

The sheaf $\Po_-$ carries a natural 
increasing filtration, by subsheaves of sections of $\Po$ with
increasing order of pole at $\infty$. These subsheaves are coherent,
but not locally free in general. On the locus of bundles with
vanishing $H^0$ and $H^1$, namely the complement of
the theta--divisor $\Theta \subset \Pic^{g-1}$, these sheaves are
vector bundles whose rank is the order of pole.

The sheaf $\Po_-$ carries one additional structure -- namely, an
action of the ring $\Oo(X\sm\infty)$ of functions away from $\infty$.
This action is compatible with the filtrations on
$\Oo(X\sm\infty)\subset\K$ and $\Po_-$: for $f\in \Oo(X\sm\infty)$
with $n$--th order pole, the sheaf $\Po_-/f\cdot \Po_-$ is a rank $n$
vector bundle over the locus $\on{Pic}^{g-1}\sm \Theta$.

\subsubsection{Proposition.} Let $\Ll\in\on{Pic}^{g-1}\sm\Theta$, 
and $f\in \Oo(X\sm\infty)$ with precisely $n$--th order pole at
$\infty$. Then the rank $n$ bundle $\Po_-/f\cdot\Po_-$ restricted to
the formal disc $\wh{D}_t$, with its natural filtration and
connection, is a $GL_n$--oper.

\subsubsection{} Equivalently, to every 
$\Ll$ and $f$ we assign an $n$--th order operator on $\Dt$. This
extends to a homomorphism $K:\Oo(X\sm\infty)\to {\mathcal D}_t$, such
that $f\cdot \psi=K(f)\cdot \psi$, where $\psi$ (the Baker--Akhiezer
function) is a section of $\Po_-$ with first order pole
(cf. \cite{Dr}).

\subsubsection{Example.}    \label{Fourier example}
Let $X=\Pp$, with $z$ a coordinate on ${\mathbb A}^1$ with first order
pole at $\infty$. We construct a formal one--parameter deformation of
the trivial line bundle $\Oo_X$ by the action of
$-z\in\K/\Oo$. Analytically, this means we are multiplying the
transition function at $\infty$ by $e^{zt}$, where $t$ is a parameter
on $\D$. The resulting line bundle on $\Pp\times\Dt$ has a connection
in the $\Dt$ direction (that is, an action of
$\frac{\partial}{\partial t}$). The connection does not affect the
trivialization of $\Oo_X$ on $\Pp\sm\infty$ (in which coordinate it is
written as $\frac{\partial}{\partial t}$). Using the transition
function to pass to a trivialization on the punctured disc $\Dx$
around $\infty$, the connection becomes $\frac{\partial}{\partial
t}-z$ and the constant section $1$ on $\Pp\sm\infty$ is written as
$\psi(z,t)=e^{zt}$, which has an essential pole at $\infty$ (but is
well defined when $t$ is a formal parameter). The Krichever
homomorphism $K$ induced by $ze^{zt}=\frac{\partial}{\partial
t}e^{zt}$ is simply the Fourier transform, sending $\C[z]$ to
$\C[\frac{\partial}{\partial t}]$.

\subsubsection{} Suppose now that $X$ is endowed with a degree $n$ map
$\phi$ to $\Pp$, such that $\phi\inv(\infty)=\infty$ (so that $\phi$
is completely branched over $\infty$). Giving such a $\phi$ is the
same as specifying a function $f\in\Oo(X\sm\infty)$ with $n$--th
order pole at $\infty$. Applying the above construction to $f$ we
attach an $n$--th order differential operator $L$ with principal
symbol $1$ (equivalently, a $GL_n$--oper, see \lemref{opers-to-diff})
on the $t$--disc to every line bundle $\Ll\in\Pic\sm\Theta$.  The key
fact is the following

\subsubsection{Theorem}\cite{Kr} The action of $\KO$ on $\Pic$
corresponds to the flows of the $n$th KdV hierarchy on the space of
all $GL_n$--opers on $\Dt$.

\subsubsection{Remark} The main result of this paper is an extension of
the above construction of commuting flows on opers on the
formal disc from the case of line bundles on $X$ to the case of
principal $G$--bundles on $X$, where $G$ is a semisimple algebraic
group. The oper connection and the flows will come from an action of
an ind--group $A/A_+$ generalizing $\KO$, the filtrations will come
from a refinement of the order--of--pole filtration, and the action of
$\Oo(X\sm\infty)$ will be replaced by the data of a reduction to the
ind--group $G(X\sm\infty)$.

\subsection{Other Perspectives.}

\subsubsection{The Fourier--Mukai Transform.}    \label{Rothstein}
The Krichever construction may be described as an application of the
Fourier--Mukai transform, as was discovered by Rothstein \cite{Ro2}
(see also \cite{N1,N2}), thus clarifying the meaning of the Krichever
homomorphism $a\cdot \psi =K(a)\cdot \psi$.

The Fourier--Mukai transform is the equivalence between the derived
categories of $\Oo$--modules on an abelian variety $A$ and its dual
$A^{\vee}$, obtained by convolution with the universal line bundle on
the product.  Laumon \cite{Lau} generalized this transform,
establishing in particular an equivalence between the derived category
of $\mathcal D$--modules on $A$ and the derived category of modules
over a sheaf $\Oo^{\nat}$ of commutative $\Oo$--algebras
($\Oo^{\nat}=\Oo(A^\nat)$, where $A^{\nat}$ is the moduli of line
bundles on $A$ equipped with a flat connection).

Now let $A=\on{Jac}_X$, so that the dual variety
$A^{\vee}=\on{Jac}_X$ as well.
Consider the Abel--Jacobi map $a_{\infty}:X\hookrightarrow
\Jac$ based at $\infty$. Let $\Oo_X(*\infty)$ denote the sheaf of
holomorphic functions on $X$ with arbitrary poles at $\infty$
allowed. 
Rothstein \cite{Ro2} proves that the pushforward
of $i_*\Oo_{X\sm x}$ to $\Jac$ under the Abel--Jacobi map has a natural
$\Oo^{\nat}$--module structure. Therefore the
result of applying the Fourier-Mukai transform to
$(a_{\infty})_*\Oo_X(*\infty)$ is a ${\mathcal D}$--module on $\Jac$.
But this transformed sheaf is easily seen to be
precisely the sheaf
$\Po_-$ on the Jacobian. Thus we obtain a
$\mathcal D$--module structure on $\Po_-$. This structure is
consistent with our constructions above, in the sense that
the action on $\Po_-$ of
the subalgebra $\C\pa_1\subset\K/\Oo$ 
comes from its embedding into $\mathcal D$.

\subsubsection{Formal Jacobians.}\label{formal Jacobians}
We wish to comment on the geometric significance of the ind--group
$\KO$ acting on $\Pic$, following \cite{CC} (see \cite{AMP} where the
ideas of \cite{CC} are explained and developed in the context of
conformal field theory). This group represents the moduli functor of
line bundles on the disc, trivialized away from the basepoint.  Thus
$\KO$ may be considered as a substitute for the Picard variety of the
disc. As we mentioned above, it is isomorphic to the constant group
scheme $\mathbb Z$ times the universal Witt formal group $\wh{\mathbb
W}$, thus identifying the latter as the Jacobian of the disc. It
carries a formal Abel--Jacobi map, whose tangent line is
$(\K/\Oo)_{-1}$. There are formal analogues of many of the usual
properties of the Jacobian, including the Fourier--Mukai transform
(following a general construction of Beilinson).

\subsubsection{Concluding Remarks.}\label{Krichever sheaves}
The action of $\KO$ on the Picard scheme of a curve identifies the
formal neighborhood of any $\Ll\in \Pic$ with a double quotient of the
group ind--scheme $\Kx$. It is in this fashion that the Krichever
construction relates to the general ideas of \secref{general}.
However, since this only captures a formal piece of the Picard
variety, it is hard to characterize the connections coming from
arbitrary line bundles all at once\footnote{One may substitute this
  formal picture by an adelic one, realizing the entire Picard as a
  double quotient for the group of id\`eles of $X$, though id\`ele
  bundles with connection seem rather daunting.}.  The solution
adopted in the theory of Krichever sheaves (\cite{Lau2,new}),
paralleling the theory of elliptic sheaves, is to retain the entire
curve $X$, and to consider line bundles $\Ll$ on $X$ times a
differential scheme $(S,\pa)$, with the $\pa$ action lifting to $\Ll$.
One then finds that the scheme $\Pic^{g-1}\setminus\Theta$ with its
$\pa_1$--action classifies Krichever sheaves (of rank $1$) for $X$
(compare \secref{differential schemes} -- the non-abelian version of
this statement will be discussed in \secref{geometric opers},
\secref{geometric interp}).  In this work, we will concentrate on the
moduli of $G$--bundles for $G$ semisimple, which do have a simple
global double quotient description. Once we introduce ``abelianized''
versions of these moduli, we obtain interesting flows and a
construction of differential data extending the above picture for line
bundles.

\section{Loop Groups and Moduli Spaces.}\label{moduli spaces}

In order to develop an analog of the Krichever construction for $G$--bundles
on curves, we wish to apply the general construction of
\secref{general} in the case when $G$ is the (formal) loop group $LG$,
$H$ its subgroup $LG_-^X$ of loops that extend ``outside'' on an
algebraic curve, and $K$ is $A_+$, an abelian subgroup of $LG_+$
(loops that extend ``inside''). This is the subject of the rest of
this paper.

In this section we introduce the loop groups and review the relation
between their quotient spaces and moduli of bundles on curves. We also
write an explicit form for the flat connections in an important
special case.

\subsection{Loop Groups.}\label{loop groups}

In the rest of this paper, unless noted otherwise, $G$ will denote a
connected semisimple algebraic group over $\C$.

Recall the setting of \secref{alg-geom-setting}. Let $LG$ be the group
ind--scheme $G(\K)$, whose $R$--points are the $R((z\inv))$--points of
$G$.  We refer to $LG$ as the loop group.  The subgroup $LG_+\subset
LG$ is defined to be the group scheme (of infinite type) $G(\Oo)$.
The Lie algebra of $LG$ is the loop algebra $L\g=\g(\K)$, with
positive half $L\g_+=\g(\Oo)$. These algebras may be identified, after
choosing a coordinate, with $\g((z\inv))$ and $\g[[z\inv]]$
respectively.  The loop algebra carries a natural filtration,
generalizing the filtration on $\K=\mathfrak{gl}_1(\K)$:

\subsubsection{Definition.}\label{homogeneous filtration} The
homogeneous filtration on the loop algebra
is defined by  
$$L\g_{\geq l} = \{ f \in L\g | fz^{-l} \in L\g_+ \}.$$ 
The induced filtration on the loop group will similarly be denoted by
$LG_{\geq l}$. Both filtrations are independent of the choice of $z$.

\subsubsection{}
Define $LG_-^X=LG_{\leq 0}\subset LG$ to be those
loops which extend holomorphically to maps $X\sm \infty\to G$. We
reserve the notation $LG_-$ for the case when $X = \Pp$. For any
projective $X$, $LG_-^X\cap LG_+\cong G$, as the only global loops are
constants. In the case $X=\Pp$, this leads to a direct sum
decomposition on the level of Lie algebras, $L\g\cong L\g_-\oplus
L\g_{\geq 1}$. (This is the infinitesimal form of the Birkhoff
decomposition, \cite{PS}.)

We now introduce infinite Grassmannians and interpret them as moduli
spaces of bundles.  For a detailed treatment of this material, we
refer the reader to \cite{BL,LS,Tel}.

An important fact about principal $G$--bundles on algebraic curves is
that if $G$ is semisimple, then any $G$--bundle on an affine curve
over $\C$ is trivial \cite{Har}. It follows that (in our previous
notations) a $G$--bundle on $X$ may be trivialized on $X\sm\infty$ and
on $\D$, and is thus determined by a transition function on $\Dx$,
which is an element of the loop group. This provides a description of
the set of isomorphism classes of $G$--bundles on $X$ as a double
quotient of the loop group. However, to obtain a similar statement for
moduli {\em stacks} (that is, to recover the algebraic structure
behind this set) one must appeal to a theorem of Drinfeld and Simpson
\cite{DSi}, which gives a version of the above trivialization
statement for families.

\subsubsection{Definition.}\label{Grassmannian definition}
Let ${\mathcal L}{\mathcal G}$ be the stack that classifies the
$G$--torsors on $X$ (for $G$ semisimple) equipped with
trivializations on $X\sm\infty$ and on $\D$. More precisely, given a
scheme $S$, ${\mathcal L}{\mathcal G}(S)$ is a groupoid whose objects
are $G$--torsors on $S \times X$ with a trivialization on $S \times
(X\sm\infty)$ and $S \times \D$, and morphisms are isomorphisms between
such objects.

The infinite Grassmannian $\Gr$ of $X$ is the moduli stack that
classifies $G$--torsors on $X$, trivialized on $\D$. More precisely,
given a scheme $S$, $\Gr(S)$ is a groupoid whose objects are
$G$--torsors on $S \times X$ with a trivialization on $S \times \D$,
and morphisms are isomorphisms between such objects.

The following remarkable description of the moduli stack of
$G$--bundles is due to Beauville--Laszlo and Drinfeld--Simpson
\cite{BL,DSi} (see \cite{Tel,Sor} for more detailed discussions).

\subsubsection{Uniformization Theorem}\cite{BL,DSi}\label{moduli
as coset}
\begin{enumerate}
\item[(1)] The stack ${\mathcal L}{\mathcal G}$ is representable by
the ind--scheme $LG$.

\item[(2)] For any scheme $S$ and any $G$--torsor $\Po$ on $S \times
X$, the restriction of $\Po$ to $S \times (X\sm\infty)$ becomes
trivial after an \'etale base change $S' \arr S$.

\item[(3)] The moduli stack ${\mathcal M}_G$ of $G$--torsors on $X$ is
canonically isomorphic to the double quotient stack $LG_-^X\bs LG/LG_+$.
It is smooth and of finite type.
\end{enumerate}

\subsubsection{Proposition.}    \label{Gr-iso}
$\Gr$ is canonically isomorphic to $LG_-^X\bs LG$, and is
representable by a scheme of infinite type.

\subsubsection{Proof.}
The fact that $\Gr$ is canonically isomorphic to $LG_-^X\bs LG$ follows
immediately from parts (1) and (2) of \thmref{moduli as coset}. Indeed,
part (2) of \thmref{moduli as coset} shows that the canonical
forgetful morphism $p: {\mathcal L}{\mathcal G} \arr \Gr$ is
surjective. Part (1) shows that ${\mathcal L}{\mathcal G} \simeq
LG$. The group $LG_-^X$ acts simply transitively on trivializations of a
$G$--torsor on $X\sm\infty$, which are the fibers of $p$. Hence $p$
gives us an isomorphism $LG_-^X \bs LG \simeq \Gr$.

It remains to prove that $\Gr$ is a scheme. The following proof was
communicated to us by C.~Teleman. Let $LG_{\geq n}\subset LG_+$
$(n>0)$ be the congruence subgroup with Lie algebra $L\g_{\geq n}$
(see \secref{homogeneous filtration}), consisting of loops regular at
$\infty$ and agreeing with the identity $1\in LG$ there to order
$n$. Thus the double quotient stack $LG_-^X\bs LG/ LG_{\geq n}$ is the
stack of $G$--torsors on $X$ equipped with a trivialization on an
$n$--th order neighborhood of $\infty$.  For any $G$--torsor
${\mathcal P} \in LG_-^X\bs LG/LG_+$ we can find an $n>0$ so that a
choice of level $n$ structure on the bundle fixes all of its
automorphisms. More precisely, there is a fine local moduli scheme for
bundles near $\mathcal P$ with level $n$ structure at $\infty$ (see
\cite{Tel}, Construction 3.12).  It follows that for every $\E\in
LG_-^X\bs LG$ there is an $LG_+$--invariant Zariski neighborhood $U$
and an $N>0$ such that $U/LG_{\geq n}$ is an affine scheme for
$n>N$. Thus $U$ represents the projective limit of $U/LG_{\geq n}$ in
the category of affine schemes, and hence is an affine
scheme. Therefore every $\E\in\Gr$ has a Zariski neighborhood which is
an affine scheme (of infinite type), so that $\Gr$ itself is a scheme
of infinite type.

\subsubsection{Warning.} 
It is important to note that the ``thick'' or ``in'' Grassmannian
$\Gr=LG_-^X\bs LG$ is {\em not} the loop Grassmannian considered,
e.g., in \cite{Ginzburg,MV,LS}, which is the ind--scheme $LG/LG_+$
that classifies the $G$--torsors on $X$ trivialized {\em outside} of
$\infty$ (and is independent of $X$). In particular,
$\Gr$ is an ordinary scheme (of infinite type) which does depend on
the curve $X$.  It is however closely related to the Sato Grassmannian
and its sub--Grassmannians studied in \cite{SW}. In algebraic geometry
the ``in'' and ``out'' Grassmannians are very different, while in the
analytic context of \cite{SW} this distinction is obscured (in genus
$0$) since Fourier series on $S^1$ can be infinite in both directions.

\subsubsection{} The Grassmannian comes equipped with
several universal bundles. As a homogeneous space $LG_-^X \bs LG$, it
comes with a tautological $LG$--bundle $\T(\Dx)$, as in \secref{double
cosets}, with a reduction $\T(X\sm\infty)$ to $LG_-^X$ and a
trivialization. From the moduli space description, there is a
tautological $G$--torsor $\T$, on $\Gr\times X$ whose fiber over
$\E\times x$ is the fiber of $\E$ at $x$.  The bundle $\T(\Dx)$ is
recovered as the sections of $\T$ over $\Dx$. More generally the
sections of $\T$ over an affine subscheme $X' \subset X$ form a
$G(X')=\on{Mor}(X',G)$ torsor. In particular, for every point $x\in X$
there is a $G$--torsor $\T(x)$ on $\Gr$, whose fiber at a point $\E\in
\Gr$ is the fiber $\E_x$ of $\E$ at $x$.  Since $\Gr$ parameterizes
bundles which are trivialized on $\D$, the bundles $\T(X')$ for
$X'\subset \D$ are canonically trivialized.

\subsubsection{Proposition.}\label{Grassmann lifting}
Let $U\subset X\sm \infty$ be a subscheme.
Then the $LG$--action on $\Gr$ lifts to the tautological bundle
$\T(U)$.

\subsubsection{Proof.} 
The total space of the $LG_-^X$--bundle $\T(X\sm\infty)$ is naturally
identified with $LG$, and hence it is clearly $LG$--equivariant. For
general $U \subset X\sm \infty$, the bundle $\T(U)$ is associated to
$\T(X\sm\infty)$ under the restriction $\T(X\sm\infty) \arr \T(U)$. In
other words, $\T(U) = G(U) \times_{LG_-^X} LG$, and hence $T(U)$ is also
$LG$--equivariant (cf. \lemref{lifting}).

Geometrically, the lifting property can be interpreted as saying that
when we change a bundle by deforming the transition function near
$\infty$, fibers away from $\infty$ are unchanged.

\subsubsection{} It follows from \propref{Grassmann lifting} and
\lemref{action->connection} that we may construct connections on
various subgroups of $LG$ by considering their action on $\Gr$ lifted
to the tautological bundles $\T(U)$. In order to obtain interesting
connections, however, we will need to pick out interesting subgroups
of $LG$ and nontrivial structures on $\T(U)$ they preserve. The desire
to obtain commuting families of flows on the resulting spaces of
connections singles out {\em Heisenberg subgroups} of $LG$, and we
will take up this idea in \secref{Heisenbergs}.

\subsection{The Big Cell.}\label{the big cell}

Consider the action on $\Gr$ of the subgroup $LG_+\subset LG$ of loops
that extend to $\D$. Acting on a pair $(\E,\phi) \in \Gr$, an element
$g \in LG_+$ does not change the $G$--torsor $\E$, but changes the
trivialization $\phi$ of $\E|_{\D}$ to $\phi g\inv$.

Note that $\Gr = LG_-^X \bs LG$ has a distinguished point
corresponding to the identity coset. From the point of view of the
moduli description of $\Gr$, this is the pair $(\E_0,\phi_0)$, where
$\E_0$ is a trivial $G$--torsor on $X$, and $\phi_0$ is its
trivialization on $\D$, which extends to a global trivialization on
the whole $X$.

Let $\Gro\subset\Gr$ be the $LG_+$--orbit of $(\E_0,\phi_0)$. This is
a scheme of infinite type that classifies bundles on $X$, trivialized
on $\D$, which admit a global trivialization.  Consider the bundle
$\T(X)$ of global sections of the tautological $G$--torsor $\T$ over
$\Gro \times X$ along $X$. Since $X$ is projective, the only global
sections of a trivial bundle are the constant sections, so that
$\T(X)$ is a $G$--torsor over $\Gro$. Furthermore, for any $x_1,x_2\in
X$, there are canonical isomorphisms $\T(x_1)\cong \T(X)\cong \T(x_2)$
obtained from restricting global sections to the different fibers.
These isomorphisms enable us to transfer extra structures, such as
decompositions or connections, from one fiber to another.

\subsubsection{}    \label{trivialization}
Unfortunately, most interesting group actions do not preserve the
subscheme $\Gro\subset\Gr$, so our main construction cannot be applied
there.  However, if $\Gro$ were open, we could restrict the action of
any Lie subalgebra of $L\g$, and any formal subgroup of $LG$, to
$\Gro$.  The orbit $\Gro$ is open when $H^1(X,\g)=0$, which is
satisfied when $X=\Pp$ is the projective line. Therefore from now on
we reserve the notation $\Gro$ for the case of $\Pp$ and call $\Gro$
the {\em big cell}.

The stabilizer of the $LG_+$--action at $(\E,\phi) \in \Gro$ consists
of elements of $LG_+$ which extend to all of $\Pp$ as automorphisms of
$\E$, namely the global sections of the adjoint group scheme $\E
\times_G \on{Ad} G$ (where $G$ acts on itself by conjugation). In the
realization $\Grp=LG_-\bs LG$, this stabilizer is the intersection
$LG_-\cap LG_+\cong G$. Thus we obtain:

\subsubsection{Lemma.} The big cell $\Gro$ is canonically isomorphic to
$G\bs LG_+$. Furthermore let $LG_{>0}\subset LG_+$ be the congruence
subgroup, consisting of loops which take the value $1\in G$ at
$\infty$. Then we have a canonical factorization $LG_+ = G \cdot
LG_{>0}$, and therefore $\Gro$ is isomorphic to $LG_{>0}$. Thus $\Gro$
may be identified with a pro--unipotent group, and hence it is
isomorphic to a projective limit of affine spaces.

\subsubsection{}    \label{cantriv}
Recall that the total space of the bundle $\T(\Pp\sm\infty)$ over
$\Grp$ is naturally identified with $LG$. The restriction of
$\T(\Pp\sm\infty)$ to $\Gro$ is then identified with an open part of
$LG$ that consists of elements $K$ admitting the factorization $K =
K_-K_+$, with $K_-\in LG_-$ and $K_+\in LG_{>0}$. This
factorization is unique. We will similarly denote by $k=k_- + k_+$ the
direct sum decomposition of
$$L\g=L\g_-\oplus L\g_{\geq 1}\cong \g[z]\oplus z\inv \g[[z\inv]]$$
into negative and positive halves.

It follows that $\T(\Pp\sm\infty)|_{\Gro}$ is
canonically trivialized: the fiber over $K_+ \in LG_{>0} \simeq
\Gro$ is identified with $LG_-$ by sending $K \in
\T(\Pp\sm\infty)|_{K_+}$ to $K_-$.

\subsubsection{}
We are in the setting of \lemref{action->connection}, where $M = \Gro,
\T = \T(\Pp\sm\infty)$. For simplicity, assume that $\ab\cong \C p$ is
a one--dimensional Lie subalgebra of $L\g$ and choose as $A$ its
formal group $\Ap = \{e^{tp}\}$. The group $\Ap$ acts on $\Gro$. Hence
we obtain for each $\E \in \Gro$ a connection on the $LG_-$--bundle
$\pi_{\E}^*(\T(\Pp\sm\infty))$ over $\Ap$ (here $\pi_{\E}: \Ap \to
\Gro$ is the $\Ap$--orbit of $\E$). The above trivialization of
$\T(\Pp\sm\infty)$ induces a trivialization of
$\pi_{\E}^*(\T(\Pp\sm\infty))$, and allows us to write down an
explicit formula for this connection.

\subsubsection{Lemma.}\label{trivial connection}
In the trivialization of $\T(\Pp\sm\infty)$ induced by the
factorization of loops, the connection operator on the $LG_-$--bundle
$\pi_{\E}^*(\T(\Pp\sm\infty))$ takes the form
$$\nabla=\partial_t+(K_+(t) p K_+(t)\inv)_-,$$ where $K_-(t)K_+(t)$ is
the factorization of $K_+ e^{-tp}$, and $K_+=K_+(0)$ is the
representative of $\E \in \Gro$ in $LG_{>0}$.

\subsubsection{Proof.}
The total space of the bundle $\T(\Pp\sm\infty)$ is an open part
$LG^{\circ}$ of $LG$ that consists of elements admitting factorization
$K = K_- K_+$. The group $A$ acts on it as follows:
$$
e^{tp}: K \mapsto K e^{-tp}.
$$
Consider an element $K_+ \in LG_{>0} \simeq \Gro$. The fiber of
$\T(\Pp\sm\infty)$ over $K_+$ consists of all $K' \in LG^{\circ}$,
which can be represented in the form $K' = K_- K_+$. By construction
of \lemref{action->connection}, the flat sections of
$\pi_{K_+}^*(\T(\Pp\sm\infty))$ are precisely the pull-backs of the
$A$--orbits $K'(t)$ of such $K'$ in $LG^{\circ}$.

Under the trivialization of $\T(\Pp\sm\infty)$ introduced in
\secref{trivialization}, the pull-back of $\T(\Pp\sm\infty)$ to $A$ is
identified with the trivial $LG_-$--bundle. The $A$--orbit of $K_+$ in
$\T(\Pp\sm\infty) = LG^{\circ}$ looks as follows: $K(t) = K_+
e^{-tp}$. Hence the corresponding section of the trivial
$LG_-$--bundle over $A$ is $K_-(t)$, where we write $K(t) = K_-(t)
K_+(t)$. This is a flat section with respect to our
connection. Therefore the connection operator reads
$$
\nabla = K_-(t) \pa_t K_-(t)^{-1} = \pa_t - K_-(t)^{-1} K'_-(t).
$$
Now we find:
$$
K'_-(t) K_+(t) + K_-(t) K'_+(t) = - K(t) p,
$$
and so
$$
K_-(t)^{-1} K'_-(t) + K'_+(t) K_+(t)^{-1} = - K_+(t)pK_+(t)^{-1}.
$$
This gives us the formula
$$
K_-(t)^{-1} K'_-(t) = - (K_+(t)pK_+(t)^{-1})_-,
$$
and the lemma follows.

\subsubsection{Remark.}    \label{triv}
In particular, we see that if $A\subset LG_+$ then we obtain a trivial
connection operator $\pa_t$ (i.e., the connection preserves our
trivialization). From this point of view the $LG_+$--action on $\Gr$
is not interesting. The action of $LG_-$, however, and in particular
of Heisenberg subgroups of $LG_-$, is the subject of our interest,
since they can be identified with the KdV flows. We will return to the
above calculation in \secref{Main}, and use it to derive the zero
curvature representation of soliton equations.

\section{Heisenbergs and Spectral Curves.}\label{Heisenbergs}

In this section we describe the geometry of Cartan subgroups of loop
groups, also known as Heisenberg subgroups. The action of these
subgroups on the moduli spaces from \secref{moduli spaces} will
produce interesting integrable systems of KdV type. We also discuss
the theory of spectral curves, introduce filtrations on the loop
algebra associated with a Heisenberg subalgebra, and consider
examples.

\subsection{Basic Properties.}\label{basic properties}
We first recall some facts about Cartan subgroups of $G$, in a form
convenient for generalization.  Next we introduce Heisenberg subgroups
and their spectral curves. The latter are used to explain (following
Kazhdan and Lusztig \cite{KL}) the theorem of Kac and Peterson
\cite{KP} classifying Heisenberg subgroups up to conjugacy.

Let $H$ be a Cartan subgroup of $G$, $N(H)$ its normalizer, and
$W=N(H)/H$ the Weyl group.  The variety of all Cartan subgroups of $G$
(hence of all Cartan subalgebras of $\g$) is naturally identified with
$G/N(H)$. Equivalently, $G/N(H)$ parameterizes $N(H)$--reductions of
the trivial $G$--torsor: the $N(H)$--torsor corresponding to $H' \in
G/N(H)$ is the torsor $\on{Isom}_G(H',H)$ of conjugacies between $H'$
and $H$. The choice of an isomorphism of groups $[\rho]: H \to H'$
reduces this $N(H)$--torsor to an $H$--torsor (conjugacies inducing
the given isomorphism). This gives a point of the variety $G/H$, which
is a $W$--torsor over $G/N(H)$.  This $W$--torsor is usually
identified with the set of Borel subgroups $B'\subset G$ containing
$H'$.

\subsubsection{Definition.}
A {\em Heisenberg subgroup} of $LG$ is a subgroup obtained by
restriction of scalars from a Cartan subgroup of $G(\K)$.
A {\em Heisenberg subalgebra} of $L\g$ is the Lie algebra of a
Heisenberg subgroup.

\subsubsection{Remarks.}    \label{locet}

\begin{enumerate}

\item[(1)] According to this definition, the Heisenberg subgroups are
abelian. The terminology is explained by the fact that the pull--back
of a Heisenberg subgroup to the Kac--Moody central extension of $LG$ is
a Heisenberg group.

\item[(2)] Recall that a Cartan subgroup $A$ of $G(\K)$ is by definition a
subgroup of $G(\K)$, which becomes isomorphic to the maximal torus
(i.e., a product of multiplicative groups of maximal dimension) over
the algebraic closure $\ol{\K}$ of $\K$. Since $A$ is an algebraic
subgroup, it becomes isomorphic to a maximal torus over a finite
extension of $\K$.

\item[(3)] Heisenberg subalgebras are the maximal
commutative subalgebras of $L\g$ consisting of semisimple elements.

\item[(4)] A Heisenberg subgroup $A\subset LG$ of the loop group is
  uniquely determined by a classifying map $C_A:\Dx\to G/N(H)$ (i.e. a
  family of Cartans of $G$).
\end{enumerate}

\subsubsection{}
The simplest example of a Heisenberg subgroup of $LG$ is the {\em
homogeneous Heisenberg} $LH$, consisting of loops $\Dx\to H$ into the
constant Cartan $H\subset G$. It is given by the constant classifying
map $\Dx\to[H]\in G/N(H)$.  A Heisenberg subgroup is said to be of
homogeneous type if it is $LG$--conjugate to $LH$. In other words,
after conjugacy the classifying map $C_A:\Dx\to G/N(H)$ maps $\Dx$
to a constant $H\subset G/N(H)$.

Since the field $\K$ of Laurent series is not algebraically closed
(equivalently, since the punctured disc is not simply-connected), the
Cartan subgroups of $LG$ are not all conjugate.  Intuitively, this
happens because they may experience monodromy around the puncture, and
these monodromies are given by automorphisms of the Cartan subgroup of
$G$, in other words, by the action of the Weyl group $W$ of $G$.  The
monodromy may best be described as a class in the Galois cohomology
group $H^1(\on{Gal}(\ol{\K}/\K), W)$, as in \cite{KL}.

\subsubsection{Definition.}\label{monodromy}
Consider a Heisenberg subgroup $A\subset LG$, given by its classifying map
$C_A:\Dx\to G/N(H)$. The {\em spectral curve}
$\Dx[A]\twoheadrightarrow \Dx$ is defined to be the pull--back under
$C_A$ of the $W$--cover $G/H\twoheadrightarrow G/N(H)$.  The monodromy
of this Galois cover is a well--defined conjugacy class $[w]$ in $W$,
which is called the type of $A$.

\subsubsection{} 
If two Heisenbergs $A$ and $A'$ are $LG$--conjugate, then their
spectral curves are automatically isomorphic. Now the spectral curve
$C_A$ of $A$ may be described as the $W$--torsor associated to the
$N(H)$--torsor on $\Dx$ of all local conjugacies of $A$ to
$LH$. Denote by $\un{A}$ the sheaf of groups (in the \'etale topology) on
$\Dx$ defined by $A$.  Consider the pullback $\un{A_C}$ of $\un{A}$ to
$C$, so that the spectral curve of $A_C$ has a tautological
section. Since every $H$--torsor on $C_A$ (or $\Dx$) is trivial (which
follows from the vanishing of $H^1(\on{Gal}(\ol{\K}/\K),{\mathbb
G}_m)$, see \cite{KL}), it follows that we can lift this tautological
section to a conjugacy of $A_C$ and $LH_C$:

\subsubsection{Proposition.}    \label{pullback}
The pullback $\un{A_C}$ of $\un{A}\subset \un{G}$ to its own spectral
curve $C=\Dx[A]$ is conjugate to the homogeneous Heisenberg $\un{H}_C$
on $C$.

\subsubsection{Corollary \cite{KP}}\label{Heisenberg conjugacy}
The Cartan subgroups of the loop group $LG$ are classified, up to
conjugacy, by the conjugacy classes in the Weyl group of $G$.

\subsubsection{Remark.} 
The spectral curve $\Dx[A]$ is usually disconnected. If we pick
locally an isomorphism $A\to LH$ (in other words a sheet of the
spectral curve) we obtain a reduction of the $W$--torsor $\Dx[A]$ to
the cyclic subgroup $\Z/n\Z$ of $W$ generated by the monodromy,
corresponding to picking a component $\wt{C}$ of the curve
$\Dx[A]$. In terms of a coordinate $z$ on $\Dx$, $\wt{C}$ is
isomorphic to the Galois cover $\zeta^n=z$. Thus the loop group
$LG_{\wt{C}}$ on $\wt{C}$, namely the sections of the constant group
scheme $G$ over $\wt{C}$, is isomorphic to $G((z^{-\frac{1}{n}}))$.
Thus \lemref{pullback} may be paraphrased as saying that if we allow
ourselves to take $n$th roots of $z$, we may conjugate $A$ to $LH$.
The restriction from $\Dx[A]$ to $\wt{C}$ is inessential -- it simply
allows us to think of sections of $G\times \wt{C}$ over $\wt{C}$ as a
loop group and not a product of loop groups.

Sometimes, when speaking about spectral curves, we will restrict
ourselves to a component $\wt{C}$ of $C$. This should be clear from
the context.

\subsection{The Principal Heisenberg.}\label{principal Heisenberg}
In this section we will discuss the most prominent Heisenberg
subgroup, the {\em principal} Heisenberg. The important features of
Heisenberg subgroups can be seen clearly in this case. To make contact
with the material of \secref{moduli spaces}, we wish to view $LG$ as
attached to the disc at a point $\infty$ on a curve $X$. Although all
local results below can be stated for an arbitrary curve $X$, we will
assume, for concreteness, that $X=\Pp$ and $z$ is a global coordinate
on $\Pp$ with simple pole at $\infty$.

Thus $L\g=\g((z\inv))$ is the (formal) loop Lie algebra at $\infty$
with Lie group $LG=G(\C((z\inv)))$. The positive and negative parts
$L\g_+=\g[[z\inv]]$ and $L\g_-=\g[z]$ consist of loops that extend to
the disc $\D$ at $\infty$ and to $\Pp\sm\infty$ respectively.

\subsubsection{}\label{e theta}
Let $\g=\n_+ \oplus \h \oplus \n_-$ be a Cartan decomposition of
$\g$. Here $\h$ is a Cartan subalgebra of $\g$ and $\n_+$, $\n_-$ are
the upper and lower nilpotent subalgebra. Let $\bb_+ = \h \oplus \n_+$
be the (upper) Borel subalgebra of $\g$. Recall that $\g$ has
generators $h_i, e_i, f_i (i=1,\ldots,\ell=\on{rank}\g$), where $h_i
\in \h, e_i \in \n_+, f_i \in \n_-$. Denote by $e_{\theta}$ (resp.,
$f_\theta$) a non-zero element of $\n_+$ (resp., $\n_-$) of weight
(minus) the maximal root $\theta$.

Recall that $L\g$ has Kac-Moody generators $h_i, i=1,\ldots,\ell; e_i,
f_i, i=0,\ldots,\ell$, where $f_i = f_i \otimes 1, e_i = f_i \otimes
1, i=1,\ldots,\ell$; $f_0 = e_{\theta} \otimes z, e_0 = f_\theta
\otimes z^{-1}$.

Introduce the ``cyclic element'' of $L\g$
\begin{equation*}\label{p-1}
p_{-1}=\sum_{i=0}^\ell f_i = f_1+\ldots+f_\ell + e_{\theta} \otimes z.
\end{equation*}
In the case of $\sln$, with conventional choices, we have
\begin{equation}\label{p minus one}
p_{-1}=\left( \begin{array}{ccccc}
0 & 0& 0&\cdots & z\\
1 & 0 &0&\cdots & 0\\
0 & 1 &0&\cdots & 0\\
\vdots&\vdots&\ddots&\cdots&\vdots\\
0& 0 &\cdots &1&0
\end{array} \right)
\end{equation}
This is a regular semisimple element of $L\g$, so that its centralizer
in $LG$ is a Heisenberg subgroup $A$ of $LG$ which is called the
principal Heisenberg.  The Lie algebra $\ab$ of $A$ is the centralizer
of $p_{-1}$ in $L\g$. It contains a unique element of the form
\begin{equation*}\label{p1}
p_1=c_1\cdot e_1+\ldots+c_l\cdot 
e_\ell + f_{\theta} \otimes z\inv
\end{equation*}
(with $c_i \neq 0$ for all $i=1,\ldots\ell$). In the case of $\sln$,
we have
\begin{equation*}
p_1=\left(\begin{array}{ccccc} 0 & 1& 0&\cdots & 0\\ 0&0&1&\cdots&0\\
\vdots&\ddots&\ddots&\ddots&\vdots\\ 0&0&0&\cdots&1\\
z\inv&0&0&\cdots&0
\end{array}\right).
\end{equation*}
Note that $p_1\in L\g_+$ is regular at $\infty$,
and is also a regular semisimple element
of $L\g$. 

Define the principal gradation on $L\g$ by
setting $\deg e_i = - \deg f_i = 1, \deg h_i = 0$.  The principal
Heisenberg Lie algebra $\ab$ is homogeneous with respect to this
gradation, and has a basis $p_i$ with $i$ (modulo the Coxeter number)
an exponent of $\g$ (see \cite{Kac1}). Except in the $D_n$ case, all
the homogeneous components $\ab^i\subset \ab$ are one--dimensional,
and the $(-1)$--component $\ab_{-1}=\C\, p_{-1}$ is always so.  In the
$\sln$ case, we may take $p_{-i}=p_{-1}^i$, which is in $\sln$ if $i$
is not divisible by $n$.

\subsubsection{The monodromy.}
For $z\in\Pp\sm\{0,\infty\}$ (in particular on the punctured disc near
$\infty$) , $p_{1}(z)$ (or $p_{-1}(z)$) is a regular semisimple
element of $\g$, and hence defines a unique Cartan subalgebra
$\ab(z)$. There are two important features of the principal family
$\ab(z)$. The first is that as $z$ undergoes a loop around $\infty$,
the algebra $\ab(z)$ undergoes a monodromy, which is a well--defined
conjugacy class in the Weyl group $W$ of $G$. In fact, for the
principal Heisenberg this is the conjugacy class of {\em Coxeter
  elements}, which are the products of simple reflections in $W$ taken
in an arbitrary order (see \cite{Kostant TDS,Kac1} and
\secref{monodromy}).

In the case of $\sln$ this monodromy has the following explicit
description. For every $z\in\Cx=\Pp\sm\{0,\infty\}$, the fiber of the
trivial bundle $\Pp\times \C^n$ over $z$ decomposes under the Cartan
$\ab(z)$ into a direct sum of $n$ lines.  This defines an $n$--fold
branched cover of $\Pp$, whose fiber over $z$ is given by the set of
eigenvalues of $p_{1}(z)$ on $\C^n$. Recall that we also have a
spectral curve $\Dx[A]$ as in the general setting, which is a
principal $W=S_n$--cover of $\Dx$. The $n$--fold cover above is the
bundle associated to $\Dx[A]$ under the permutation representation of
$S_n$.

Since the Coxeter number (the order of a Coxeter element) of $\sln$
equals $n$, the monodromy statement above translates into the
statement that this $n$--fold cover is fully branched over $\infty$,
i.e., is described by the equation $\zeta^{-n}=z\inv$. (Note that
$p_{1}^n=z\inv\cdot\on{Id}$, so that $p_{1}$ itself plays the role of
$\zeta\inv$.)  The Coxeter class for $S_n$ consists of $n$--cycles, and by
labeling the cover in different ways we obtain the different
$n$--cycles as monodromies. It is also worth noting that the $n$--fold
cover is isomorphic to a component of the reducible curve
$\Dx[A]$ -- the $S_n$--bundle may be reduced to a cyclic subgroup
generated by an $n$--cycle.

For general $G$, the spectral curve $\Dx[A]$ of the principal
Heisenberg is a union of the fully branched cyclic covers of $\Dx$ of
order the Coxeter number of $G$.

\subsubsection{Degeneration at $\infty$.}
The second important feature of $\ab(z)$ is the way it degenerates at
$z=\infty$, as a subspace of $\g$, to an $l=\on{rk}\g$--dimensional
abelian subalgebra $\ab_{\infty}$. This limit is the centralizer of
the regular nilpotent element $\ol{p}_1=p_1(\infty)$ of $\g$.
Specifically, in the case of $\sln$,
\begin{equation}\label{p one}
\ol{p}_1=
\left( \begin{array}{ccccc}
0&1&0&\cdots&0\\
0&0&1&\cdots&0\\
\vdots&\ddots&\ddots&\ddots&\vdots\\
0&0&0&\cdots&1\\
0&0&0&\cdots&0
\end{array}
\right).
\end{equation}
The centralizer $\ab_{\infty}$ of $\ol{p}_1$ consists solely of upper
triangular matrices.  Geometrically, we see that the $n$--sheeted
spectral cover associated with $\ab$ can be completed to the $n$--fold
cover of the disc defined by $\zeta^{-n}=z\inv$, completely branched
at infinity. The elements $p_i\in \ab$ are now identified with the
powers $\zeta^{-i}$ of the coordinate on the spectral cover. Thus in
particular the principal gradation on $\ab$ agrees with the gradation
in powers of the coordinate upstairs.  The upper triangular matrices
are those that preserve the canonical filtration induced on the
$n$--dimensional vector space $\C[\zeta\inv]/z\inv \C[\zeta\inv]$.
Thus for $z\neq \infty$ we obtain a decomposition of the trivial
vector bundle into lines, while at $z=\infty$ we retain the structure
of a flag.

In general there is a unique Borel subalgebra in the fiber
$L\g|_{\infty}\cong \g$ which contains $\ab_{\infty}$, namely $\bb_+$.
Consider the principal filtration, whose $i$--th piece $L\g^{\geq i}$
consists of elements of degree $\geq i$ in the principal gradation.
It defines the above Borel subalgebra as $\bb_+=L\g_-\cap L\g^{\geq
0}$.  The resulting filtration on $\g$ is canonically determined by
the choice of $\bb_+$, and hence $\ab_+$.

It is also useful to note that due to the structure of the principal
filtration, there is a distinguished line $\ab^{-1}\subset \ab/\ab_+$
given by the intersection of $\ab$ with $L\g_{\geq -1}$. Thus the
element $p_{-1}$ is determined intrinsically by the structure of
$\ab$, up to a constant.  This line is the analog of the line in
$\K/\Oo$ of elements with first order pole -- it consists of ``first
order poles in the spectral coordinate $\zeta$''.

Finally, the full branching of the spectral curve translates
algebraically into the statement that $A$ is an anisotropic torus in
$LG$: it does not contain any split torus (i.e., a product of copies
of $\K^\times$) as a subgroup. It follows that the ind--group $A/A_+$
is in fact just the formal group associated to $\ab/\ab_+$.

\subsection{Singularities of Heisenbergs.}\label{singularities}
In this paper we are interested in the classification of Heisenberg
subgroups of $LG$ not only up to $LG$--conjugacy, but up to
$LG_+$--conjugacy. Thus we are interested in the ``integral models''
of Cartan subgroups in $G(\K)$, that is abelian subgroups $\un{A}_+$
of the trivial group scheme $\un{G}$ on the {\em unpunctured}
disc $\D$ (as opposed to the group scheme $\un{A}$ over the punctured
disc $\Dx$).  The reason for that will become apparent in
\secref{Higgs}, where we consider the moduli of pairs
$(\E,\E^{A_+})$, where $\E$ is a $G$--torsor on a curve $X$, and
$\E_{\un{A}}$ is a reduction of $\E|_{\D}$ to such an $\un{A}_+$.

\subsubsection{}
While Heisenberg subalgebras of $L\g$ are classified up to
$LG$--conjugacy by conjugacy classes in the Weyl group, their
determination up to $LG_+$--conjugacy is in fact much more subtle (as
was first explained to us by R. Donagi).  This finer structure
describes how families of Cartan subalgebras on $\Dx$ degenerate at
$\infty\in \D$. It may also be described in terms of singularities of
completed spectral curves. The integrable systems we construct in
subsequent sections reflect this intricate behavior.  

Let $\ab\subset L\g$ be a Heisenberg Lie algebra of homogeneous
type. In other words we can find $g\in LG$ such that $g \ab
g\inv=L\h$.  Such a $g$ is unique up to left multiplication by the
normalizer $N(L\h_C)\subset LG_C$. Note that $N(L\h_C)/LH_C\cong W$ is
the (finite) Weyl group $W$, since for the homogeneous Heisenberg
$N(L\h)/LH=L(N(\h)/H)$ and loops into $W$ are necessarily constant.
It follows that the map $char:\ab\to L\h/W$ induced by $g$ is uniquely
defined. Let $\ab_+=\ab\cap L\g_+$, and let $g(\ab_+)\subset L\h$ be
the image of $\ab_+$ under $g$.

\subsubsection{Lemma.}\label{normalization} The image 
$g(\ab_+)$ satisfies $g(\ab_+)\subset L\h_+ \subset L\h$. Moreover,
$L\h_{\geq N}\subset g(\ab_+)\subset L\h_+$ for $N\gg 0$, where
$L\h_{\geq j}=L\h\cap L\g_{\geq j}$ is the homogeneous filtration.

\subsubsection{Proof.} According to the Chevalley theorem, we have a
``characteristic polynomial'' map $\g\to \on{Spec}\C[\g]^G\cong\h/W.$
Over the field $\K$ of Laurent series, we obtain a map $char:L\g\to
L(\h/W)$ from the pointwise application of the characteristic
polynomial. In particular, $\ab$ is homogeneous when $char$ sends it
to $L\h/W$, loops which may be lifted to $\h$.  In fact, to conjugate
$\ab$ to $L\h$ is equivalent to choosing a lift $\ab\to L\h$ of the
characteristic map.  Since the characteristic map is defined over
$\Oo\subset\K$, it follows that for $a\in\ab_+$, $char(a)\in
L\h_+/W\subset L\h/W$. Since the map $g$ is a lift of the
characteristic map, the first statement follows immediately.

Since $g$ has a finite order of pole at $\infty$, it follows that
$g\inv(h)\in L\g_+$ for $h\in L\h_{\geq N}$ and $N\gg 0$ (as can be
easily seen for example for $SL_n$ and hence in a faithful matrix
representation). This proves the second statement.

%\subsubsection{} In fact, we have the following fine classification
%theorem for Heisenberg subalgebras:

%\subsubsection{Theorem \cite{BZ}}\label{plus conjugacy} 
%Let $\ab$,$\ab'$ be two Heisenberg subalgebras of $L\g$. Then $\ab$
%and $\ab'$ are $LG_+$--conjugate if and only if there is an
%$LG$--conjugacy between them sending $\ab_+$ to
%$\ab'_+$. Equivalently, letting $char(\ab),char(\ab')\subset L(\h/W)$
%denote the images of the characteristic map, $\ab$ and $\ab'$ are
%$LG_+$--conjugate if and only if $char(\ab_+)=char(\ab'_+)$.

\subsubsection{} Now let $\ab\subset L\g$ be a general Heisenberg
subalgebra. Recall that the pullback $\ab_C$ of $\ab$ to its own
spectral curve $C=\Dx[\ab]$ is conjugate to the homogeneous Heisenberg
$L\h_C$. If $g_1, g_2 \in LG_C$ satisfy $g_i \inv L\h_C g_i=\ab_C$,
then the conjugation by $g_1 g_2^{-1}$ preserves the homogeneous
filtration on $L\h_C$. Therefore the following definition makes
sense.

\subsubsection{Definition.}    \label{canonical filtration}
The {\em canonical filtration} $\{\ab^{\geq
i}\}$ on $\ab$ is the filtration induced on $\ab$ from the homogeneous
filtration on $L\h_C$ via the embedding $\ab\subset\ab_C$ and any
conjugacy $\ab_C\cong L\h_C$.

We say $\ab$ is {\em smooth} if the subalgebra $\ab^+:=\ab^{\geq 0}$
and its subalgebra $\ab_+=\ab\cap L\g_+$ coincide.

%\thmref{plus conjugacy} implies:

%\subsubsection{Corollary.} For every conjugacy class $[w]\in W$, there
%exists a unique smooth Heisenberg subalgebra of type $[w]$, up to
%$LG_+$--conjugacy. 

\subsubsection{}\label{difference of A's} 
Let $A_+\subset A^+\subset A$ be the ind--groups 
corresponding to the Lie algebras $\ab_+\subset \ab^+\subset \ab$. It
follows that the group $A^+$ is isomorphic to $LH_+^{[w]}$, where
$LH^{[w]}$ is a smooth Heisenberg of the same type. Hence in particular
  $A^+$ is actually a group scheme, and $A/A^+$ is the product of the
  formal group associated to $L\h^{[w]}/L\h_+^{[w]}$ by a lattice. The
  ``difference'' between $A$ and $LH^{[w]}$ is the finite--dimensional
  group scheme $A^+/A_+$.

\subsubsection{} The Heisenberg algebras which have been studied in the
literature \cite{KP,dGHM,Feher} are graded with respect to an
associated gradation on the loop algebra. These graded Heisenbergs all
satisfy the smoothness condition $\ab_+=\ab^+$ (as well as
$LG^+\subset LG_+$), and provide standard smooth representatives of
all $LG$--types.

We want to stress that there are plenty of $LG_+$--conjugacy classes
of non-smooth Heisenbergs -- there exist continuous families of those
within given type $[w]$.  The collection of $LG_+$--classes of
Heisenbergs of type $[w]$ is naturally parameterized by the
infinite--dimensional double quotient $N(A)\bs LG/LG_+$, where $A$ is
any such subgroup (say, a graded representative) and $N(A)$ its
normalizer.  While the integrable systems that have been studied in
the literature so far are associated to graded Heisenbergs only, in
this paper we construct integrable systems attached to arbitrary
Heisenberg subalgebras possessing strongly regular elements (that is,
for all Heisenbergs of type $[w]$ where $[w]$ varies through many, but
not all, conjugacy classes in $W$).

\subsubsection{Remark.}\label{picture of A} 
The picture of Heisenberg algebras that emerges
closely parallels the structure of branched covers of $\D$. First
there is a topological invariant, the monodromy of the cover over
$\Dx$, which is resolved by passing to an \'etale cover. For
Heisenberg algebras this is the monodromy of the spectral curve and
the passage from $\ab$ to $\ab_C$. The next step is to consider the
behavior at the marked point. The normalization of the completed
spectral curve is a smooth curve isomorphic to $\D$, so that we have a
finite codimension embedding $\Oo_C\subset \Oo\cong\C[[z\inv]]$ of
coordinate rings. This is the meaning of the embedding
$\ab_+\hookrightarrow L\h_+$ of \lemref{normalization}. The filtration
$\{\ab^{\geq i}\}$ is the filtration by order of pole on
$\D=\on{Spec}\Oo$, transferred to the subring $\Oo_C$, and the
smoothness condition $\ab^+=\ab_+$ is the Lie algebra version of the normality
condition on $C$ (so the group $A^+/A_+$ ``measures'' the singularity).

\subsection{Regular Centralizers.}
Heisenberg algebras are classified up to $LG$ conjugacy by the
associated spectral curves over $\Dx$. Their classification up to
$LG_+$ conjugacy reflects the geometry of curves over $\D$. While we
do not have a general theory of completed spectral curves for
arbitrary Heisenberg algebras, such a theory is available (thanks to
\cite{Donagi, DM,DG}) for the large class of regular Heisenberg
algebras.

Since $G/N(H)$ parameterizes the Cartan subalgebras of $\g$, it embeds
into the Grassmannian $Gr^\ell(\g)$ of $\ell$--dimensional subspaces
of $\g$. Let $\wt{G/N(H)}$ be the variety of all abelian subgroups of
$G$, which are centralizers of regular elements. It also embeds into
$Gr^\ell(\g)$ and hence can be thought of as a partial
compactification of $G/N(H)$.

\subsubsection{Definition.} A Heisenberg subgroup $A$ is regular if
the classifying map $C_{A}:\Dx\to G/N(H)$ extends to a map
$C_{A_+}:\D\to \wt{G/N(H)}$.

\subsubsection{} Equivalently, $A$ is regular if the fiber of $A_+$ at
$\infty$ is a regular centralizer in $G$.

In order to extend the spectral curve $\Dx[A]\twoheadrightarrow \Dx$
over $\infty$, we introduce the scheme $\wt{G/H}$ which classifies
pairs $(R,B)$ with $R\in \wt{G/N(H)}$ a regular centralizer, and $B$
Borel subgroup of $G$ containing $R$.  There is an obvious morphism
$\wt{G/H}\to\wt{G/N(H)}$, and the fiber over $R$ may be identified
with the fixed point scheme $\B^R$ of $R$ in the flag variety
$\B=G/B$.

\subsubsection{Definition.}\label{completed curve} 
The {\em (completed) spectral curve} $\D[A]$ of a regular Heisenberg
$A$ is the pullback to $\D$ of the morphism $\wt{G/H}\to\wt{G/N(H)}$
under the classifying map $C_{A_+}$ of $A$.

\subsubsection{Proposition.} Let $A,A'\subset LG$ be regular
Heisenbergs. Then $A,A'$ are $LG_+$--conjugate if and only if 
the spectral curves $\D[A]$ and $\D[A']$ are isomorphic.

\subsubsection{Proof.} It is clear that if $A,A'$ are
$LG_+$--conjugate, then $\D[A] \simeq \D[A']$. To prove the converse,
we note that there is a one-to-one correspondence between the set of
$LG_+$--conjugacy classes of regular Heisenbergs with a fixed spectral
curve $C$ and the set of isomorphism classes of Higgs bundles over $\D$
with spectral cover $C$ (see \secref{Higgs}). A theorem from \cite{DG}
states that the category of $G$--Higgs bundles on $X$ with fixed
spectral cover $\wt{X}$ carries a simply transitive action of the
Picard category of torsors over a specific (abelian) group
scheme $\wt{T}$ on $X$. In the case $X=\D$, all $G$-- and
$\wt{T}$--torsors are trivial, and so we obtain from that statement
that there is a unique up to isomorphism Higgs bundle over $\D$ with
spectral cover $C$. Therefore there is a unique $LG_+$--conjugacy
class of regular Heisenbergs with a fixed spectral curve $C$.

\subsection{Examples.}    \label{non-smooth}
\subsubsection{}
The principal Heisenberg is smooth, reflecting the smoothness of the
$n$--fold branched cover defined by taking the $n$--th root of $z$. 
In fact, the principal gradation on
$\ab$ is a refinement of the canonical filtration. 
In the $\mathfrak{sl}_2$ case, the generators
$$p_i=\left(\begin{array}{cc} 0&z^{1-i}\\z^{-i} &
0\end{array}\right)\in \ab_+, \quad \quad i>0$$ 
of $\ab_+$ are conjugate to elements
$$\left(\begin{array}{cc} z^{-i+1/2}&0\\0&-z^{-i+1/2}
\end{array}\right),$$ which lie in the positive part of the
homogeneous Heisenberg subalgebra of ${\mathfrak s}{\mathfrak
l}_2((z^{-1/2}))$.

The centralizer of 
$$\left(\begin{array}{cc} 0& 1\\z^{-2} &0\end{array}\right)\in
L\mathfrak{sl}_2,$$ is, as in the principal case, a regular nilpotent
centralizer at $\infty$. But this Heisenberg Lie algebra is
$LG$--conjugate to the {\em homogeneous} Heisenberg.

It is easy to see that it is not smooth. Indeed, denote by
$$
\wt{p}_i=\left(\begin{array}{cc} 0&z^{2-i}\\z^{-i} &
0\end{array}\right), \quad \quad i\in\Z,
$$
the generators of this Lie algebra $\ab$. The matrix
\begin{equation}    \label{conjelement}
\left(\begin{array}{cc} -\frac{1}{2} & -z \\ -\frac{1}{2} z^{-1} &
1 \end{array}\right)
\end{equation}
conjugates $\wt{p}_i$ to
$$\left(\begin{array}{cc} z^{-i+1}&0\\0&-z^{-i+1}
\end{array}\right).$$
Hence we see that $\ab_+$ is generated by $\wt{p}_i, i>1$, while $\ab^+$
is generated by $\wt{p}_i, i>0$.

Because $\ab$ is not smooth, it is not $LG_+$--conjugate to $L\h$, as
one can see from the explicit formula \eqref{conjelement} for one of
the conjugating elements.

The spectral curve $\Dx[A]$ is isomorphic to the trivial ${\mathbb
Z}_2$--cover of $\Dx$ (the same as the spectral curve of $LH$). But
the completed spectral curve is singular: it has two irreducible
components, with a simple node over $\infty$.

The above Heisenberg is obtained from the pullback of the principal
Heisenberg to its own spectral curve, which is automatically
homogeneous. More generally, let $\ab$ be a Heisenberg subalgebra with
the limit $\ab_{\infty}\subset L\g_+/L\g_{>1}\cong \g$.  By
considering the pullback of $\ab$ to its spectral curve we obtain an
example of a Heisenberg subalgebra that $LG$--conjugate to $L\h$ but
with a limit $\ab_{\infty}$ at $\infty$.
 
\subsubsection{} Since there are
continuous families of non-isomorphic local singularities, it is easy
to find continuous families of Heisenbergs none of which are $LG_+$
conjugate. For example, one can consider a Heisenberg of homogeneous
type for $SL_4$, whose $4$--sheeted spectral curve is planar and
isomorphic to four copies of $\D$ joined at $\infty$. The tangent
lines to the four components define four points in the projectivized
Zariski tangent space at $\infty$. The cross ratio of the resulting
four points in $\Pp$ is an invariant of the curve (thus of the
associated Heisenberg) which may be varied continuously.
Below we will
assign integrable systems to $LG_+$--conjugacy classes of Heisenberg
subgroups of $LG$. It follows that there are continuous families of
integrable systems obtained by our construction.

\subsection{Filtrations.}\label{filtrations}

In this section we prove some technical results concerning
filtrations, which will be useful when we consider generalizations of
the notion of an oper in \secref{opers} (in particular \propref{DS
gauge}). 

\subsubsection{}
The homogeneous Heisenberg algebra $L\h\subset L\g$ has a strong
compatibility property with the 
homogeneous filtration \secref{homogeneous filtration}.
Denote by $L\h_{\geq i}$ the $i$--th
piece of the induced filtration, $L\h_{\geq i}=L\h\cap L\g_{\geq i}$. 
In particular $L\h_{>0}=L\h_{\geq 1}$ is the Lie algebra
of loops to $\h$ which vanish at $\infty$.
To an element $p_i\in L\h_{\geq i}$ we associate its (principal)
symbol 
$$\ol{p}_i=p_i\on{mod}L\g_{>i}\in L\g_{\geq i}/L\g_{>i}\cong \g.$$

\subsubsection{Definition.}
We say that $p_i$ is {\em strongly regular} if its symbol
$\ol{p}_i\in\g$ is a regular element, i.e., if the centralizer of
$\ol{p}_i$ in $\g$ is precisely $\h$.

\subsubsection{Lemma.} Suppose $p_i\in L\h_{\geq i}$ is strongly regular.
Then:

\begin{enumerate}
\item[(1)] $p_i$ is regular.

\item[(2)] $\on{Ker}(\on{ad}p_i)=L\h$ and
$L\g\cong L\h\oplus\on{Im}(\on{ad}p_i)$.

\item[(3)] The operator $\on{ad}p_i$
induces isomorphisms $L\g_{\geq k}/L\h_{\geq k}\to L\g_{\geq
k+i}/L\h_{\geq k+i}$. 
\end{enumerate}

\subsubsection{Proof.} We may pick a coordinate $z$ on $\D$, thereby
picking a gradation refining the homogeneous filtration. Then we may
write $p_i=\ol{p}_i+\sum_{j>i}\ol{p}_j$. Suppose
$a=\sum_{k=k_0}a_k\in L\g$
is the graded decomposition of an element satisfying $[p_i,a]=0$.
By equating each graded component of the commutator to zero, we find
equations $\sum_{j=0}^n [\ol{p}_{i+j},a_{k_0+n-j}]=0$. By induction
on $n$, and using the regularity of $\ol{p}_i$, we obtain that each
$a_k\in\h$ and hence $a\in L\h$, establishing part (1).

Now since all elements of $L\h$ are semisimple, $p_i$ is regular
semisimple, and we obtain part (2).

Since the filtration $\{L\g_{\geq k}\}$ is a Lie algebra filtration,
and since by part (2) we have an $\on{ad}p_i$--invariant decomposition
of $L\g$, we obtain a well defined operator $L\g_{\geq k}/L\h_{\geq
k}\to L\g_{\geq k+i}/L\h_{\geq k+i}$ as required, only depending on
the symbol $\ol{p}_i$.  The fact that it is an isomorphism is now an
easy consequence of the regularity of the symbol, as may be checked
using the $z$--gradation.

\subsubsection{Definition}
Let $\ab\subset L\g$ be a general Heisenberg subalgebra. By {\em a
filtration associated with $\ab\subset L\g$} we will understand a
filtration on $L\g$ induced by the homogeneous filtration on $L\g_C$
via the homomorphism $\on{Ad} g: L\g \to L\g_C$, where $g$ is {\em
an} element of $LG_C$, such that $g \ab_C g\inv = L\h_C$, where
$\ab_C$ is the pullback to $=ab$ to its own spectral curve $C$.

\subsubsection{}
The restriction of the above filtration to $\ab \subset L\g$ is
canonical, i.e., it does not depend on the choice of $g$. But the
filtration on the whole $L\g$ does depend on the choice of $g$,
because $g$ is specified by the above condition only up to left
$N(LH_C)$ multiplication. We do not know how to endow $L\g$ with a
canonical filtration that restricts to the canonical filtration on
$\ab$ defined in \secref{canonical filtration}.  However, any of these
many filtrations (for varying $g$) have the following nice property.

\subsubsection{Definition}    \label{strongly regular}
An element $p_i \in \ab$ is strongly regular if it corresponds to a
strongly regular element in $L\h_C$.

\subsubsection{Lemma.}\label{strict compatibility} 
Let $\{L\g^{\geq k}\}$ be any filtration on the loop algebra
associated with $\ab\subset L\g$ as above. If $p_i\in \ab^{\geq i}$
is strongly regular, then $\on{ad}p_i$ induces isomorphisms $L\g^{\geq
k}/\ab^{\geq k}\to L\g^{\geq k+i}/\ab^{\geq k+i}$.

\subsubsection{Proof.} 
Note first that the centralizer of $p_i$ in $L\g$ is the intersection
of its centralizer $\ab_C$ in $L\g_C$ with $L\g$, hence $\ab$, so that
$p_i$ is indeed a regular element of $L\g$ and
$L\g\cong\on{Ker}(\on{ad}p_i) \oplus\on{Im}(\on{ad}p_i)$. Now the
statement of the lemma for $L\g_C$ and $\ab_C$ follows immediately
from the corresponding statement for $L\h_C$.

Let us decompose $L\g_C$ into characters for the Galois group of
$C\twoheadrightarrow \Dx$, producing a Lie algebra gradation. 
The subalgebra $L\g\subset L\g_C$ consists
of the invariants of this action, the zeroth graded component. Since
$p_i\in L\g$, $\on{ad}p_i$ preserves the gradation on $L\g_C$. It
follows that the statement of the lemma must be true component by
component in $L\g_C$, and hence in $L\g$ itself, as desired.

\subsubsection{Remark.}\label{strongly regs} 
The most important classes of Heisenbergs, the homogeneous and
principal, have canonical filtrations associated with them and contain
many strongly regular elements. It is not true unfortunately that
every Heisenberg algebra contains strongly regular elements. However
this property depends only on the type of $\ab$, and not on its fine
structure: if $\ab$ is $LG$--conjugate to $\ab'$ then $\ab$ contains
strongly regular elements precisely when $\ab'$ does. Thus one may
inquire for which conjugacy classes $[w]$ in the Weyl group
Heisenbergs of type $[w]$ contain strongly regular elements. It
suffices to answer this question for graded Heisenbergs, and in this
setting it has been shown in \cite{DF} (see also \cite{Feher}), that
the graded Heisenberg of type $[w]$ contains strongly regular elements
precisely when $[w]$ is a regular conjugacy class of the Weyl group.
Those have been previously classified by Springer. For instance, in
the case of $\g={\mathfrak s}{\mathfrak l}_n$ these conjugacy classes
correspond to partitions on $n$ either into equal integers, or into
equal integers plus $1$ (see \cite{FehHarMar}).

The generalized Drinfeld--Sokolov construction of integrable systems
described below is only applicable to Heisenbergs of these
types. However, within each ``topological'' type of Heisenbergs, there
usually exist continuous families of Heisenbergs, which are not
$LG_+$--conjugate, and therefore continuous families of integrable
hierarchies.

\section{Abelianization.}\label{Higgs}

In this section we combine the contents of the previous two to produce
a class of interesting moduli spaces, the abelianized Grassmannians
$\Gr_A$, parameterizing bundles with additional structure on the
disc. This structure is a reduction to the positive part $A_+$ of a
Heisenberg subgroup of $LG$.  These moduli spaces carry actions of
abelian (ind--)groups $\Ahat$, deforming the abelian structure on the
disc.  Following the prescription of \secref{general}, we use these
actions to relate these moduli to moduli of special connections, which
will turn out to be affine opers (see \secref{opers}).

We then identify the abelianized Grassmannians (when $A$ is a {\em
  regular} Heisenberg) with moduli spaces of Higgs bundles. Higgs
  bundles have a natural interpretation in terms of line bundles on
  spectral curves, thus ``abelianizing'' our moduli spaces by
  comparing them to Picard varieties. In particular, we will identify
  the ind--group $A/A_+$ as a generalized Prym variety associated with
  the spectral curve attached to $A$. This allows us to interpret the
  action of $A/A_+$ on $\Gr_A$ as a generalization of the Jacobian
  flows in the Krichever construction \secref{Krichever}. In
  \secref{Main} we identify the action of $A/A_+$ on $\Gr_A$ with a
  hierarchy of generalized KdV flows on the space of affine opers.

\subsection{Abelianized Grassmannians.}
Let $A\subset LG$ be a Heisenberg subgroup, with Lie algebra
$\ab$. Let $A_+=A\cap LG_+$ be its positive half, with Lie
algebra $\ab_+$.

\subsubsection{Definition.}\label{A Grassmannian}
The $A$--Grassmannian $\Gr_A$ is the moduli stack of
$G$--torsors on $X$, equipped with a reduction $\E^{A_+}$ of the
$LG_+$--torsor $\E|_{\D}$ to $A_+$.

\subsubsection{Lemma.} $\Gr_A\cong\Gr/A_+$.

\subsubsection{Proof.} 
Since $\Gr$ parameterizes $G$--torsors on $X$ with an isomorphism
$\un{G}|_{\D} \simeq \E|_{\D}$, there is a natural surjection
$\Gr\twoheadrightarrow \Gr_A$, and the group $A_+$ acts transitively
along the fibers.

\subsubsection{Remark.} The relevance of the double quotient space
$LG_-^X\bs LG/A_+$ and its ad\`elic versions to the study of
integrable systems has been pointed out in \cite{EF1}.

%\subsection{The Abelianized Big Cell.}

\subsubsection{}    \label{Gro}
Since $\Gr_A$ is the quotient of the scheme $\Gr$ by the action of the
group $A_+$, it is an (Artin) algebraic stack. $\Gr_A$ does not,
however, possess a coarse moduli scheme. This is because the
automorphisms of different bundles $\E\in\Gr_A$ may be very
different. These automorphisms are given by the sections of the
adjoint group--scheme $\on{Ad}\E^{A_+}$ of the $A_+$--torsor
associated with $\E$, which extend to global sections of the
group--scheme $\on{Ad}(\E)$.  Thus if $\E$ is trivial as a
$G$--torsor, $\on{Aut}(\E)\cong A_+\cap G$, since $G=LG_-^X\cap LG_+$
are the global sections of the trivial $G$--torsor. For general $\E$
the automorphisms will be isomorphic to the intersection of $A_+$ with
some conjugate of $LG_-^X$ which is not transversal to $LG_+$.

However, recall that (in the case $X=\Pp$) $\Grp$ has a big cell
$\Gro$, which is canonically isomorphic to $G \bs LG_+ \simeq
LG_{>0}$. Define the big cell $\Gro_A$ of $\Grp_A$ as the image of
$\Gro$ under the projection $\Grp \twoheadrightarrow \Grp_A$. While
$\Grp_A$ is the moduli stack of triples $(\V,\V_-,\V^{A_+})$, where $\V$
is an $LG$--torsor with $LG_-$--reduction $\V_-$ and $A_+$--reduction
$\V^{A_+}$, $\Gro_A$ is its open substack that classifies those
triples for which the $LG_-$--reduction $\V_-$ and the induced
$LG_+$--reduction $\V^{A_+} \times_{A_+} LG_+$ are in general
position.

We want to show that
$\Gro_A$ possesses a coarse moduli space, $Gr^{\circ}_A$. Hence we
will be able to translate results concerning the stack $\Gro_A$ into
results concerning the scheme $Gr^{\circ}_A$. We start with a
preliminary result.

\subsubsection{Lemma.}    \label{silly}
\begin{enumerate}
\item[(1)] Let $A\subset LG$ be a Heisenberg subgroup, and $G\subset
LG$ the subgroup of constants. Then either $A$ is $LG_+$--conjugate to
a homogeneous Heisenberg $LH$ (where $H$ is a Cartan subgroup of $G$),
in which case $A\cap G$ is a Cartan subgroup of $G$; or $A\cap G$ is
the center of $G$.

\item[(2)] Let $A_0 = A\cap G$, $A'_+ = A \cap LG_{>0}$, for
homogeneous $A$, and $A'_+ = \exp \ab_+$ for all other $A$. Then
$A_+ = A_0 \cdot A'_+$.
\end{enumerate}

\subsubsection{Proof.} An element $a\in A\cap G$ is a 
constant section of the group scheme $\un{A}$ on $\D$. Its value at
the base point $\infty \in \D$ is then a semi--simple element of
$G$. Thus, the fiber $A_\infty$ of $\un{A}$ is an abelian subgroup of
$G$ of dimension greater than or equal to $\ell$, containing a
semi-simple element $a$. It follows that either $A_\infty$ is a Cartan
subgroup, or $a$ is a central element of $G$. This proves part
(1). Part (2) follows immediately from part (1).

\subsubsection{Lemma.} $\Gro_A$ is isomorphic to the quotient of
the affine scheme of infinite type $Gr^{\circ}_A = G\bs LG_+/A'_+$ by
the trivial action of the group $A_0$.

\subsubsection{Proof.} Since $\Gro \simeq G\bs LG_+$, $\Gro_A$ is isomorphic to $G\bs LG_+/A_+$. By
\lemref{silly}, (2), it is isomorphic to the quotient of $G\bs
LG_+/A'_+$ by the trivial action of $A_0$. Since $g A'_+ g\inv \cap G
= \{ 1 \}$ for any $g \in LG_+$ by \lemref{silly},(1), the
pro--unipotent group $A'_+$ acts freely on $G \bs LG_+$. It is easy to
find locally a transversal slice for this action. Therefore $G\bs
LG_+/A'_+$ is a scheme.

\subsubsection{Actions.} The stack $\Gr_A=LG_-^X\sm LG/A_+$ occupies an
important intermediate position between the Grassmannian
$\Gr=LG_-^X\sm LG$ and the moduli stack of bundles ${\mathcal
M}_G^X=LG_-^X\sm LG/LG_+$.  Unlike the Picard group $\Pic$ (the moduli
scheme of line bundles on $X$), ${\mathcal M}_G^X$ does not have a
group structure and carries no natural group actions.  The
Grassmannian does carry an action of the entire loop group $LG$.  It
is however too big, with ``redundant'' directions coming from the
action of its subgroup $LG_+$, which merely changes the trivialization
on $\D$.

The abelianized Grassmannian $\Gr_A$ carries an interesting remnant of
the $LG$ action on $\Gr$, which is similar to the action of the
Jacobian on itself. Indeed, the quotient $\Gr/A_+$ carries the right
action of the normalizer of $A_+$ in $LG$. In particular, the group
$A$ acts on $\Gr_A$.  Of course, the right action of $A_+$ on $\Gr_A$
is trivial, and so we are left with the action of the quotient
ind--group $A/A_+$, discussed in \secref{difference of A's}. Recall that
this group has three ``parts'' -- the finite--dimensional group scheme
$A^+/A_+$, a finite rank lattice, and the formal group of $\ab/\ab^+$.

\subsubsection{Remark.}
This action is the key property of $\Gr_A$. It
can also be used as a motivation for studying $\Gr_A$. Indeed,
\lemref{action->connection} gives us a general procedure for
constructing flat connections from the actions of a group $A$ on a
homogeneous space. We take as the homogeneous space, the Grassmannian
$\Gr$. The most interesting case of actions to consider is that of a
maximal abelian subgroup. That is why we look at the action of a
Heisenberg subgroup $A$. But as we explained in \remref{triv}, the
$LG_+$--action on $\Gr$ is not interesting. Therefore we mod out $\Gr$
by $A_+ = A \cap LG_+$ and look at the residual action of $\Ahat$.

\subsection{The Principal Grassmannian.}\label{principal Grassmannian}

In this section we concentrate on the abelianized Grassmannian
$\Gr_A=LG_-^X\bs LG /A_+$ in the case when $A \subset LG$ is the
principal Heisenberg and $X=\Pp$. For $G=GL_n$ it may be described as
the moduli stack of rank $n$ vector bundles on $\Pp$, identified over
$\D$ with the pushforward of a line--bundle from the $n$--fold
branched cover $\zeta^n=z$. For general $G$, $\Grp_A$ classifies
$G$--torsors $\E$ on $\Pp$, with a local structure on $\D$ that can be
described as follows: $\E|_{\Dx}$ is reduced to an abelian group
subscheme of the constant group scheme $\un{G}$, whose fibers are
Cartan subgroups of $G$ undergoing a Coxeter class monodromy around
$\infty$ and degenerating at $\infty$ to a regular nilpotent
centralizer.

\subsubsection{}
The tautological $LG$--bundle $\T(\Dx)$ on $\Grp_A$ comes with
reductions $\T(\Pp\sm\infty)$ and $\T^{A_+}(\D)\subset\T(\D)$ to $LG_-$ and
$A_+\subset LG_+$, respectively. Since $A_+$
defines a unique Borel subgroup $B\subset G$ at $\infty$, it follows
that the fiber of any $\E\in\Grp_A$ at $\infty$ has a canonical
flag. In particular, $\T(\Dx)$ has a canonical reduction to the
Iwahori subgroup $LG^+\subset LG_+$, whose sections take values in $B$
at $\infty$.

Recall that since $A$ is smooth and ``maximally twisted'',
i.e. anisotropic, the ind--group $A/A_+$ is actually a formal group.
Thus the formal group $\Ahat$ and its Lie algebra $\ab/\ab_+$ act on
$\Grp_A$.  We are particularly interested in the action of the element
$p_{-1}\in\ab/\ab_+$ and of the formal one-dimensional additive group
$\At=\{e^{tp_{-1}}\}$ that it generates.

\subsubsection{} We are now in the setting of our general result,
\propref{abstract isomorphism}, describing the correspondence between
the double quotients $\HGK$ equipped with an action of a group $A$ and
the moduli of certain flat bundles on $A$. Namely, we take
$\Grp_A=LG_-\bs LG /A_+$ as the double quotient and $\At$ as the group
$A$. The moduli stack of flat bundles on the other side of the
correspondence classifies quadruples $(\V,\nabla,\V_-,\V^{A_+})$,
where $\V$ is an $LG$--torsor on $\At$ with a flat connection
$\nabla$, a flat $LG_-$--reduction, and an $A_+$--reduction
$\V^{A_+}$--reduction in tautological relative position with respect
to $\nabla$.\footnote{Since $\At$ is isomorphic to the formal disc,
  $\nabla$ automatically induces a trivialization of $\V$ and
  $\V_-$.} We denote this moduli stack by $\M^{\Pp}_{p_{-1}}$.

On the other hand, we can consider the action of the whole group
$A/A_+$ on $\Grp_A$. The corresponding moduli space $\M^{\Pp}_A$
classifies quadruples as above defined on all of $A/A_+$ rather than
on $\At$. Then \propref{abstract isomorphism} gives us the following
result.

\subsubsection{Proposition.}\label{pre-pre-oper}
$\Grp_A$ is canonically isomorphic to $\M^{\Pp}_{p_{-1}}$, and to
$\M^{\Pp}_A$.

\subsection{General Case.}    \label{general M}
Now let $X$ be an arbitrary smooth curve, $\infty$ a point of $X$,
and $LG$ the loop group corresponding to the formal neighborhood of
$\infty$. Recall that $LG$ has subgroups $LG_-^X$ and $LG_+$. Let
$A$ be an arbitrary Heisenberg subgroup of $LG$, $A_+ = A \cap LG_+$,
and $\ab$, $\ab_+$ be the Lie algebras of $A$, $A_+$,
respectively. Each non-zero element $p \in \ab/\ab_+$ gives rise to a
one-dimensional formal additive subgroup $\Ap = \{ e^{tp} \}$ of the
group $A/A_+$.

Consider the corresponding abelianized Grassmannian $\Gr_A$. The group
$\Ap$ acts on $\Gr_A$ from the right, and we can again apply
\propref{abstract isomorphism}. Denote by $\M^X_{A,p}$ the moduli
stack that classifies quadruples $(\V,\nabla,\V_-,\V^{A_+})$, where
$\V$ is an $LG$--torsor on $\At$ with a flat connection $\nabla$, a
flat $LG_-$--reduction $\V_-$, and an $A_+$--reduction
$\V^{A_+}$ in tautological relative position with
$\nabla$.

Similarly, $\M^X_A$ will denote quadruples as above, but defined on
all of $A/A_+$ rather than on $\Ap$.  Then we obtain the following
generalization of \propref{pre-pre-oper}.

\subsubsection{Proposition.}    \label{general iso}
$\Gr_A$ is canonically isomorphic to $\M^X_{A,p}$ and to $\M^X_A$.

\subsubsection{Differential Setting.} Rather than consider connections
on $\Ap$, we may follow the prescription of \secref{differential
schemes} and introduce differential schemes. Thus $\V$ will now be a
$LG$--torsor on a differential scheme $(S,\pa)$ with a lifting
$\pa_{\V}$ of $\pa$, and reductions to $LG_-$ and $A_+$ which are
respectively preserved and in relative position $[-p]$ with respect to
$\pa_{\V}$. \propref{differential proof} immediately implies

\subsubsection{Proposition.}\label{differential opers} 
The pair $(\Gr_A,p)$ represents the functor which assigns to a
differential scheme $(S,\pa)$ the groupoid of quadruples
$(\V,\pa_{\V},\V_-,\V^{A_+})$ as above.

\subsection{Higgs Fields and Formal Jacobians.} 

Recall that the Heisenberg algebra $\ab$ is represented by a
classifying map $C_A:\Dx\to G/N(H)$. Since $G/N(H)$ is embedded into
the Grassmannnian $Gr^{\ell}(\g)$, which is proper, we can extend this
map to a map $\D \arr Gr^{\ell}(\g)$. The image of the base point of
$\D$ under this map will be an abelian Lie subalgebra of
$G$. Intuitively, a reduction $\E^{A_+}$ of $\E|_{\D}$ to $A_+$ is a
reduction of the $G$--torsor $\E$ on $\Dx$ to a Cartan subgroup, which
twists in a prescribed way around $\infty$, and degenerates at
$\infty$ to a reduction of $\E_{\infty}$ to the limiting abelian
subgroup of $G$.

\subsubsection{Definition. \cite{DM}}\label{Higgs def}
\begin{enumerate}
\item[(1)] A (regular) principal Higgs field on a $G$--torsor $\E$
over a scheme $Y$ is a sub--bundle $c\subset \on{ad}(\E)$ of regular
centralizers.

\item[(2)] The spectral cover associated with a Higgs bundle $(\E,c)$
on $Y$ is the scheme $Y[c]\twoheadrightarrow Y$ parameterizing Borel
subgroups of $\on{Ad}\E$ containing $c$.  More precisely, $Y[c]$ is
the fixed point scheme $(\B)_{\E}^c$ of $c$ on the relative flag
manifold $(\B)_{\E}=\E\times_G G/B$.
\end{enumerate}

\subsubsection{Remark.} The spectral cover $Y[c]$ may also be defined as
follows. Let $\wt{G/N(H)}$ be the partial compactification of $G/N(H)$
parameterizing regular centralizers in $G$. Let
$\wt{G/H}\twoheadrightarrow \wt{G/N(H)}$ be the scheme parameterizing
pairs $(r,B)$ where $r\in\wt{G/N(H)}$ is a regular centralizer and
$B\subset G$ is a Borel subgroup containing $c$.  Now trivialize the
bundle $\E$ locally on some flat covering $Y'\to Y$, so that $c$
defines a local classifying map $Y'\to\wt{G/N(H)}$. The pullback to
$Y'$ of the morphism $\wt{G/H}\twoheadrightarrow \wt{G/N(H)}$ is
independent of the trivialization up to isomorphism. Thus the
resulting covers $Y'[c]\to Y'$ glue together to give a scheme $Y[c]$
over $Y$. (Without choosing trivializations, we obtain a classifying
morphism $Y\to G\bs\wt{G/N(H)}$, and $Y[c]$ is the pullback of the
representable morphism $G\bs\wt{G/H}\to G\bs\wt{G/N(H)}$.)

\subsubsection{} Recall from \secref{basic properties} the relation
between Cartan subgroups of $G$ and reductions: a reduction of the
trivial $G$--torsor to a Cartan subgroup $H$ (in other words, a point
$s\in G/H$) is equivalent to the data of a subgroup of $G$ conjugate
to $H$, together with a distinguished group isomorphism $\phi: H'\cong
H$. Equivalently, we may replace the isomorphism $\phi$ by the data of
an identification between the set of Borels containing $H'$ and those
containing $H$. This is a consequence of the fact that $N(H)/H$ acts
faithfully on the set $\B^H$ of Borels containing $H$.

\subsubsection{Proposition \cite{DG}}\label{regular centralizers} 
Let $R\subset G$ be a regular centralizer.
\begin{enumerate}
\item[(1)] For any regular element $r\in R$, the fixed point scheme
$\B^R$ is identified with $\pi\inv(char(r))$, where
$\pi:\h\to\h/W$ and $char:\g\to\h/W$ is the adjoint quotient map. In
particular, $\B^R$ carries a natural action of the Weyl group
$W=N(H)/H$.

\item[(2)] There is a morphism $R\to \on{Hom}_W(\B^R,H)$, which is an
isomorphism for $G$ simply connected. Thus $N(R)/R$ acts faithfully on
$\B^R$.

\item[(3)] The data of a reduction of $G$ to $R$ is equivalent to
the data of a regular centralizer $R'\subset G$ and an identification
of $W$--schemes $\B^{R'}\cong \B^R$.
\end{enumerate}

\subsubsection{Corollary.}\label{reduction to Higgs}
Let $A\subset LG$ be a regular Heisenberg, with positive part
$A_+\subset A$.  Then $\Gr_A$ is naturally isomorphic to the
moduli stack of $G$--torsors $\E$ on $X$, equipped with a principal
Higgs field $c$ on $\E|_{\D}$ and an isomorphism $\D[c]\cong\D[A_+]$
of spectral covers over $\D$.

\subsubsection{Proof.} The corollary is a version of \propref{regular
centralizers},(3) over $\D$. Recall from Definitions \ref{completed
curve} and \ref{Higgs def}) that the spectral cover of a Heisenberg
algebra or Higgs field is the global version of the fixed point
schemes $\B^R$. A reduction of $\E|_{\D}$ to $A_+$ is determined by a
section of the associated scheme
$(\un{G}/\un{A}_+)_{\E}=\E\times_{\un{G}} \underline{G}/\un{A}_+$ over
$\D$.  Using \propref{regular centralizers}, we can identify sections
of $\un{G}/\un{A}_+$ with pairs $(c,\phi)$, where
$c\in\underline{G}/N(\un{A}_+)$ is a subgroup of $\un{G}$ and
$\phi:\D[c]\to \D[A_+]$ is an isomorphism of spectral covers.  It
follows that sections of $(\un{G}/\un{A}_+)_{\E}$ are identified with
pairs $(c,\phi)$, where $c\in(\underline{G}/N(\un{A}_+))_{\E}$ is a
principal Higgs field on $\E|_{\D}$ and $\phi:\D[c]\to \D[A_+]$ is an
isomorphism of spectral covers. This implies the Corollary.

\subsubsection{Higgs Bundles and Line Bundles.} 
The most important feature of Higgs bundles on a scheme $X$ is their
relation to line bundles on the associated spectral cover. We offer a
brief review of the theory as developed in 
\cite{Donagi,DM} and completed in \cite{DG}. (See also \cite{LM}.)

Let us first suppose that
$G=GL_n$, and that $\E$ is a rank $n$ vector bundle on $X$ equipped
with a subbundle $c\subset \on{ad}\E$ of Cartan subalgebras. Locally,
the action of a regular section $s$ of $c$ on $\E$ decomposes $\E$
into a direct sum of $n$ eigen--line bundles, parameterized by the $n$
eigenvalues of $s$. In other words $c$ defines an $n$--sheeted \'etale
cover $Y=X[c]_n$ of $X$, whose points over $x\in X$ represent
weights of $c_x$ on $\E_x$.  The cover $Y$ may be recovered from the
spectral cover $X[c]$ associated to the Higgs field $c$ as follows:
$\pi:Y\to X$ is the bundle $X[c]\times_{S_n}\{1,\cdots,n\}$ associated
to the principal $S_n$--bundle $X[c]$ under the permutation
representation of $S_n$ on the set $\{1,\cdots,n\}$.

We have expressed $\E\cong
\pi_*\Ll$, where $\pi:Y\to X$ and $\Ll$ is a line bundle on
$Y$. Conversely, any line bundle $\Ll$ on $Y$ pushes forward
to a rank $n$ vector bundle, equipped with a canonical semisimple
Higgs field $c$ (i.e. bundle of Cartans) preserving the decomposition
$\E_x\cong\oplus_{\pi\inv(x)}\Ll_x$.  In the case $G=SL_n$, we
obtain a correspondence between semisimple Higgs bundles
$(\E,c)$ on $X$, and line bundles $\Ll$ on $Y$ equipped with an
isomorphism $\on{det}\pi_*\Ll\cong \Oo_X$.

Now suppose the Higgs field $c$ is allowed to be an arbitrary bundle
of regular centralizers in $\on{ad}\E$. It is still possible to define
an $n$--sheeted spectral cover $Y=X[c]_n$, which will now be ramified
over the locus where $c$ is not a Cartan. Moreover $\E$ will be
identified with the pushforward of a line bundle on $Y$. Conversely,
the pushforward of a line bundle from an $n$--sheeted cover produces a
rank $n$ vector bundle. The case when $X$ is a curve, with local
coordinate $z$ near $\infty$, and $X[c]$ is locally isomorphic to the
cover $\zeta^n=z$, was discussed in \secref{principal Heisenberg}: the
stalk of $\pi_*\Ll$ at $z=\infty$ will be isomorphic to the stalk of
the bundle of $(n-1)$--jets of $\Ll$ at $\infty$. Thus instead of a
direct sum decomposition, we see a flag on the stalk. The flag is
defined by multiplication by the coordinate $\zeta$, which acts as a
regular nilpotent element on the stalk.

In \cite{Donagi,DM,DG}, the case of arbitrary $G$ is worked out in
detail. The category of Higgs bundles with a given spectral cover
(possibly ramified) is described precisely in terms of Prym varieties
(that is, in terms of $W$--equivariant torus bundles on the spectral
covers).  We only describe the very first step in this work, to
motivate our interpretation of $\Ahat$ as a Prym variety.

Suppose that $G$ is an arbitrary reductive group and $(\E,c)$ a
regular semisimple Higgs bundle on $X$ (so that $c$ is a bundle of
Cartans). The $W$--cover $X[c]\to X$ is then \'etale, and the pullback
$\E_{X[c]}$ of $\E$ to $X[c]$ carries a tautological reduction to a
Borel containing $c$. Thus given an irreducible representation $V$ of
$G$, the associated bundle $V_{\E_{X[c]}}$ on $X[c]$ has a
distinguished line bundle, defined as the highest weight space for the
tautological Borel.  Thus for every (semisimple) Higgs bundle on $X$,
one obtains a $W$--equivariant homomorphism from the weight lattice
$\Lambda$ of $G$ to the Picard variety of $X[c]$, in other words an
element of the cameral Prym variety:

\subsubsection{Definition.}\label{cameral Prym}\cite{DM} 
The cameral Prym variety
associated with the Higgs bundle $(\E,c)$ on $X$ is
$\on{Hom}_W(\Lambda,\on{Pic}(X[c]))$. 

\subsubsection{}
Recall from \secref{formal Jacobians} that the ind--group
$\KO$ plays the role of the Picard of the disc.  \defref{cameral Prym}
then suggests an analogous role for $\Ahat$.  Let us denote by
$\Kx[A]$ the invertible functions on $\Dx[A]$, considered as a group
ind--scheme. Similarly $\Ox[A]$ will be the group scheme of invertible
functions on $\D[A]$ and $\KO[A]$ the quotient group ind--scheme. The
following lemma is an easy consequence of \propref{regular
centralizers}.

\subsubsection{Lemma.} Let $G$ be simply connected, and $A\subset LG$ 
a regular Heisenberg subgroup. Then $A\cong \on{Hom}_W(\Dx[A],H)\cong
\on{Hom}_W(\Lambda,\Kx[A])$,
$A_+\cong \on{Hom}_W(\Lambda,\Ox[A])$ and $A/A_+\cong\KO[A]$.

\subsubsection{Definition.} The Prym variety $\on{Prym}(A)$ associated
with a Heisenberg $A$ is the ind--group $\Ahat$.

\subsubsection{Corollary.} 
\begin{enumerate}
\item[(1)] For any Heisenberg subgroup $A$ of homogeneous type, there
is a map $\on{Prym}(A)\twoheadrightarrow \on{Prym}(LH)$, uniquely
specified up to $W$.

\item[(2)] For $A$ regular, there is an injection
$\on{Prym}(A)\hookrightarrow \on{Prym}(\wt{A})$, where
$\un{\wt{A}}\subset \un{G}_{\D[A]}$ is the Heisenberg subgroup of
homogeneous type obtained from pulling $A$ back to $\D[A]$.
\end{enumerate}

\subsubsection{Proof.} The first statement is a consequence of 
\lemref{normalization}, which describes $\ab$ in terms of $L\h$.  The
second follows from the inclusion $i:A\hookrightarrow \wt{A}$ and the
equality $i(A_+)=i(A)\cap \wt{A}_+$.

\subsubsection{} The first part of the corollary is the analog, for
an arbitrary Heisenberg subgroup, of the morphism between the Prym
varieties corresponding to the pullback of line bundles from a cover
($\D[A]$) to its normalization ($\D[LH]$). The second is the analog of
the morphism corresponding to pulling back line bundles from a curve
($\D$) to its branched cover ($\D[A]$).

\subsection{Krichever Construction, Revisited.}
Let $G=GL_n$ and $A\subset LGL_n$ a regular Heisenberg, with
$n$--sheeted spectral cover $\D[A]_n$. Let $\Sigma$ be a projective
curve, $i:\D[A]_n\hookrightarrow \Sigma$ an inclusion and
$\pi:\Sigma\twoheadrightarrow X$ a morphism extending
$\D[A]_n\twoheadrightarrow \D$. Let $\KO[\Sigma]$ denote the
ind--group representing the quotient of $\Oo^{\ti}(\pi\inv(\Dx))$ by
$\Oo^{\ti}(\pi\inv(\D))$.

\subsubsection{Definition.} 
The $\Sigma$--Grassmannian $\Gr_A[\Sigma]$ is the
moduli stack of $GL_n$ Higgs bundles $(\E,c)$ on $X$, equipped with
an isomorphism $X[c]_n\to\Sigma$ of the $n$--sheeted spectral cover
of $c$ with $\Sigma$.

\subsubsection{Proposition.}\label{Krichever revisited} 
\begin{enumerate}
\item[(1)] The $\Sigma$--Grassmannian $\Gr_A[\Sigma]$ is a substack
of $\Gr_A$, which is preserved by the action of $\Ahat$.
\item[(2)] $\Gr_A[\Sigma]\cong \on{Pic}_{\Sigma}$.
\item[(3)] There is a natural isomorphism $\Ahat\cong
\KO[\Sigma]$, which identifies their actions on $\on{Pic}_{\Sigma}$.
\end{enumerate}

\subsubsection{Proof.}
The isomorphism of \corref{reduction to Higgs} may be reformulated in
the case $G=GL_n$, by replacing spectral covers $\D[c]$ with
$n$--sheeted spectral covers $\D[c]_n$. 
It follows that there is a natural morphism
$\Gr_A[\Sigma]\to\Gr_A$, obtained by restricting $c$ to $\D$
and composing the identification $\D[c]_n\to\pi\inv(\D)\subset
\Sigma$ and the inverse of the inclusion $i:\D[A]_n\to\Sigma$. This
morphism is a monomorphism, since $c|_{\D}$ determines the extension
$c$ on $X$ (when it exists) and $\pi:\Sigma\to X$ has no automorphisms
preserving $i$.

The action of $\Ahat$ on $\Gr_A$ affects neither the $G$--torsor
$\E|_{X\sm\infty}$, nor the reduction of $\E|_{\Dx}$ to $A$ (since the
action is deduced from an $A$ action, which clearly has this
property). It follows that the action preserves the substack
$\Gr_A[\Sigma]$, whose definition depends only on $\E|_{X\sm\infty}$
and its reduction to $A$ on $\Dx$. This completes the proof of part
$(1)$.

Part $(2)$ of the proposition is the standard identification between
the moduli stack of $GL_n$ Higgs bundles with $n$--sheeted spectral
cover isomorphic to $\Sigma$ and that of line bundles on $\Sigma$
(\cite{DG}).  Applying this statement to $\Dx$, we obtain a canonical
isomorphism between the automorphism group of a $GL_n$ Higgs bundle on
$\Dx$ with spectral curve $\Dx[A]$ and the automorphism group of a
(trivial) line bundle on $\Dx[A]$. This gives us an identification
$A\cong \Kx[\Sigma]$. The equivalent statement $A_+\cong\Ox[\Sigma]$
over $\D$ then implies the desired identification $A/A_+\cong
\KO[\Sigma]$. Moreover, it is clear that the actions of $\Ahat$ on
$\Gr_A[\Sigma]$ and of $\KO[\Sigma]$ on $\on{Pic}_\Sigma$
coincide. This concludes the proof.

\subsubsection{}
We may now replace $GL_n$ by an arbitrary reductive group $G$.
Let $A\subset LG$ be a regular Heisenberg, with
spectral curve $\D[A]$. Let $\Sigma$ be a projective
curve, $i:\D[A]\hookrightarrow \Sigma$ an inclusion and
$\pi:\Sigma\twoheadrightarrow X$ a morphism extending
$\D[A]\twoheadrightarrow \D$. 

\subsubsection{Definition.} 
The $\Sigma$--Grassmannian $\Gr_A[\Sigma]$ is the
moduli stack of $G$--Higgs bundles $(\E,c)$ on $X$, equipped with
an isomorphism $X[c]\to\Sigma$ of the spectral cover
of $c$ with $\Sigma$.

\subsubsection{Conclusion.} As in \propref{Krichever revisited}, the
stack $\Gr_A[\Sigma]$ is naturally a substack of $\Gr_A$, and is
preserved by $\Ahat$.  By the results of \cite{Donagi,DM,DG},
$\Gr_A[\Sigma]$ has a description in terms of line bundles on
$\Sigma$.  As we observed, $\Ahat$ serves as a
generalized Prym variety for the spectral curve $\D[A]$. Thus, the
$\Sigma$--Grassmannian $\Gr_A[\Sigma]$ together with the action of
$A/A_+$ is a natural generalization of the setting of the Krichever
construction to arbitrary semisimple groups $G$. The points of
$\Gr_A[\Sigma]$ correspond to an interesting class of global
``algebro--geometric'' solutions of generalized soliton
hierarchies. On the other hand, the total abelianized Grassmannian
$\Gr_A$ captures all formal solutions of those hierarchies, as we will
see in \secref{Main}.

\subsubsection{Remark.} Cherednik \cite{Ch1,Ch2,Ch3} also considers
versions of the stacks $\Gr_A[\Sigma]$ as parameterizing generalized
algebro--geometric solutions to soliton equations. He uses the
language of reductions to subgroups of the group scheme $G\times X$
over $X$, which are generically Cartans. As in \corref{reduction to
Higgs}, this is essentially equivalent to the notion of principal
Higgs bundles.

\section{Opers.}\label{opers}

In this section we study special connections,
introduced by Drinfeld-Sokolov \cite{DS} and Beilinson-Drinfeld
\cite{opers} under the name opers, and their ``spectral''
generalizations, affine
opers. Opers may be defined on any curve, and have an interpretation
in terms of differential operators.  Opers are also equivalent to a
special class of affine opers.  The Drinfeld--Sokolov gauge
(\propref{DS gauge}) relates affine opers with the connections
produced in \secref{Higgs}. Tying the loose ends together, we obtain
in the next section a broad generalization of the relation between
line bundles on a curve and differential operators reviewed in
\secref{Krichever}.

\subsection{Introducing Opers.}\label{ordinary opers}
Let $G$ be a reductive Lie group and $B$ a Borel subgroup, with Lie
algebras $\bb\subset \g$.
There is a distinguished $B$--orbit ${\mbf O}\subset \g/\bb$,
which is open in the subspace of vectors stabilized by the radical
$N\subset B$. 
$\mbf O$ consists of vectors which have precisely their negative
simple root components nonzero -- that is, the $B$-- (or $H$--)orbit
of the sum of the Chevalley generators $f_i$ in $\g/\bb$.
In the case of $GL_n$, ${\mbf O}$ can be identified with the
set of matrices of the form
\begin{equation*}
\left( \begin{array}{ccccc}
*&*&*&\hdots&*\\
+&*&*&\hdots&*\\
0&+&*&\hdots&*\\
\vdots&\ddots&\ddots&\ddots&\vdots\\ 
0&0&\hdots&+&*
\end{array} \right)
\end{equation*}
where the $+$ represent arbitrary nonzero entries.

Now recall the notion of relative position of a reduction of a bundle
with respect to a connection given in \secref{types of connections},
and the notion of $GL_n$--oper from \secref{GL opers}.

\subsubsection{Definition.}\cite{opers} Let $Y$ be a smooth curve. A
$G$--{\em oper} on $Y$ is a $G$--torsor $\E$ on $Y$ with a connection
and a reduction $\E_B$ to $B$, which has relative position $\mbf O$
with respect to $\nabla$.

In other words, the one-form $\nabla/\E_B$ takes values in ${\mbf
O}_{\E_B} \otimes \Omega^1 \subset (\g/\bb)_{\E_B} \otimes \Omega^1$
(the orbit $\mbf O$ is $\Cx$--invariant). 

The $G$--opers on $Y$ form a stack denoted by $\Op_G(Y)$, or $\Op(Y)$
when there is no ambiguity. A $\g$--oper is by definition an oper for the adjoint group
of $\g$.

\subsubsection{} Identifying $GL_n$--principal bundles equipped 
with $B$--reductions with rank $n$ vector bundles equipped with flags,
we see that this notion agrees with the $GL_n$--opers as introduced in
\secref{GL opers}.  Thus in the case of $GL_n$, opers are identified
with $n$--th order differential operators with principal symbol $1$.
The $SL_n$--opers correspond (in the notation of \secref{GL opers}) to
the case $q_1=0$, that is differential operators with vanishing
subprincipal symbol.  For the classical series $B_n$ and $C_n$, we
obtain differential operators which are either self-- or
skew--adjoint.

\subsubsection{Lemma.}\cite{opers}    \label{automor} 
If $G$ is a group of adjoint type, then the moduli stack $\Op_G(Y)$ of
$G$--opers on $Y$ is an affine ind--scheme denoted by $Op_{\g}(Y)$
(which is in fact a scheme for $Y$ projective or $Y=\wh{D}$). For
general $G$, the stack $\Op_G(Y)$ is the quotient of the (ind--)scheme
$Op_{\g}(Y)$ of $\g$--opers by the trivial
action of the center $Z(G)$.

\subsubsection{}\label{twidle opers}
Let $\DD$ be the formal disc. In \secref{Main} we will
need an explicit description of the scheme $Op_\g(\DD)=Op(\DD)$.

Consider $(\E,\nabla,\E_B) \in Op_\g(\DD)$. If we choose a
trivialization of $\E_B$, then the connection operator $\nabla$ can be
written as
\begin{equation}    \label{genoper}
\pa_t + \sum_{i=1}^\ell \phi_i(t) \cdot f_i + {\mathbf b}(t).
\end{equation}
Here $t$ is a coordinate on $\DD$, $\phi_i(t)$ are invertible elements
of $\C[[t]]$, and ${\mathbf b}(t) \in \bb_+[[t]]$.

Let $\wt{Op}(\DD)$ be the affine space of operators of the form
\eqref{genoper}. The group $B_+[[t]]$ acts on $\wt{Op}(\DD)$ by gauge
transformations corresponding to the changes of trivialization. We
have: $Op(\DD) \simeq \wt{Op}(\DD)/B_+[[t]]$. This allows us to
identify $Op(\DD)$ with a projective limit of affine spaces.

Recall that $\ol{p}_{-1} = \sum_{i=1}^\ell f_i$. We have a direct sum
decomposition $\bb_+ = \oplus_{i\geq 0} \bb_{+,i}$ with respect to the
principal gradation. The operator $\on{ad} \ol{p}_{-1}$ acts from
$\bb_{+,i+1}$ to $\bb_{+,i}$ injectively for all $i>1$. Hence we can
find for each $j>0$ a vector subspace $V_j \subset \bb_{+,j}$, such
that $\bb_{+,j} = [\ol{p}_{-1},\bb_{+,{j+1}}] \oplus V_j$. Note that
$V_j \neq 0$ if and only if $j$ is an exponent of $\g$, and in that
case $\dim V_j$ is the multiplicity of the exponent $j$. In
particular, $V_0=0$. Let $V = \oplus_{j\in E} V_j \subset \n_+$, where
$E$ is the set of exponents of $\g$. We call such $V$ a transversal
subspace.

\subsubsection{Lemma.}\cite{DS}    \label{canreps}
The action of $B_+[[t]]$ on $\wt{Op}(\DD)$ is free. Moreover, each
$B_+[[t]]$--orbit contains a unique representative of the form
\begin{equation}    \label{genoper2}
\pa_t + \ol{p}_{-1} + {\mathbf v}(t), \quad \quad {\mathbf v}(t) \in
V[[t]].
\end{equation}

%\subsubsection{Proof.}
%Given a differential operator of the form \eqref{genoper}, we can
%bring it to the form
%\begin{equation}    \label{genoper1}
%\pa_t + \ol{p}_{-1} + {\mathbf b}(t), \quad \quad {\mathbf b}(t) \in
%B_+[[t]]
%\end{equation}
%by the gauge action of an element $h(t)$ of $H[[t]]$, which is
%uniquely determined by this condition. Hence it remains to show that
%the action of $N_+[[t]]$ (where $N_+ = [B_+,B_+]$) on the operators of
%the form \eqref{genoper1} is free, and each orbit contains a unique
%element of the form \eqref{genoper2}.

%This means that each element \eqref{genoper1} can be uniquely
%represented in the form
%\begin{equation}    \label{gauge}
%\pa_t + \ol{p}_{-1} + {\mathbf b}(t) = \exp \left( \on{ad} U \right) \cdot
%\left( \pa_t + \ol{p}_{-1} + {\mathbf v}(t) \right),
%\end{equation}
%where $U \in \n_+[[t]]$ and ${\mathbf v}(t) \in V[[t]]$. We can
%decompose with respect to the principal gradation: $U=\sum_{j\geq 0}
%U_j$, ${\mathbf b}(t)=\sum_{j\geq 0} {\mathbf b}_j(t)$, ${\mathbf
%v}(t)=\sum_{j\geq 0} {\mathbf v}_j(t)$. Equating the homogeneous
%components of degree $j$ in both sides of \eqref{gauge}, we obtain
%that ${\mathbf b}_i(t) + [U_{i+1},\ol{p}_{-1}]$ is expressed in terms of
%${\mathbf v}_i(t),{\mathbf v}_1(t),\ldots,{\mathbf
%v}_{i-1}(t),U_1,\ldots,U_i$. The injectivity of $\on{ad} \ol{p}_{-1}$ then
%allows us to determine uniquely ${\mathbf b}_i(t)$ and
%$U_{i+1}$. Hence $U$ and ${\mathbf b}(t)$ satisfying equation
%\eqref{gauge} can be found uniquely by induction, and lemma follows.

\subsubsection{Corollary.} The ring of functions on $Op(\DD)$ is
isomorphic to the ring of functions on $V[[t]]$ for any choice of
transversal subspace $V \subset \n_+$.

\subsection{Affine Opers.}\label{affine opers}
While opers provide a Lie--theoretic viewpoint on differential
equations $L\psi=0$ (where $L=\partial^n-\ldots$), the affine opers
incorporate the eigenvalue problem $L\psi=z\psi$ for $L$: by rewriting
this as $(L-z)\psi=0$, we find we are dealing with a natural
one--parameter deformation of the operator $L$. Similarly, an oper for
a general group $G$ comes with a natural one--parameter ``spectral''
deformation. As Drinfeld and Sokolov discovered, it is
beneficial to consider this entire family as a single connection, by
replacing the group $G$ by its loop group in $z$.

Affine opers are thus the analogs of opers for loop groups (i.e. opers
with spectral parameter). Let $p_{-1}$ be as in \secref{principal
Heisenberg}. Recall that $p_{-1}$ is the ``loop analog'' of the
principal nilpotent $\ol{p}_{-1}$ of $\g$. On a differential scheme
(such as the line $\At$) we may thus define an affine oper as a
$LG$--bundle with connection and reductions, where the distinguished
vector field acts in relative position $p_{-1}$ (see \defref{t-line}).
However, in contrast to the case of $G$, the $LG^+$--orbit of
$[p_{-1}]\in L\g/L\g^+$ is not $\C^\times$--invariant. We will modify
this approach, so as to obtain a definition of affine opers on an
arbitrary curve, without using a distinguished vector field. However,
as we will see, affine opers (unlike opers) will only exist on curves
whose canonical bundle is trivial.

Consider the multiplicative group $\Cx$ of automorphisms of
$\Pp$ preserving the points $0$ and $\infty$. It acts naturally on
$LG\supset LG^+$ and $L\g\supset L\g^+$ by rescaling the coordinate
$z$. We can therefore form the semidirect products $\wt{LG}\supset
\wt{LG}^+$ of $LG$ and $LG^+$ with $\Cx$.  Denote by $\Oaff$ the
$\wt{LG}^+$--orbit of $[p_{-1}]\in L\g/L\g^+$. This orbit is the
$\Cx$--span of the $LG^+$--orbit of $[p_{-1}]$. It consists of all
elements of the form $\sum_{i=0}^\ell \lambda_i f_i$, where $\la_i$'s
are arbitrary non-zero complex numbers.

\subsubsection{Definition.}    \label{affoper1}
Let $Y$ be a smooth curve. An {\em affine oper} on $Y$ is a quadruple
$(\wt{\V},\nabla,\wt{\V}_-,\wt{\V}^+)$, where $\wt{\V}$ is a
$\wt{LG}$--torsor on $Y$ with a connection $\nabla$, a flat reduction
$\wt{\V}_-$ to $\wt{LG}_-$ and a reduction $\wt{\V}^+$ to $\wt{LG}^+$
in relative position $\Oaff$.

The affine opers on $Y$ form a stack which is denoted by $\AOp(Y)$.

\subsubsection{Geometric Reformulation} \label{geometric opers}
An affine oper on $Y$ may also be described as follows.  Let ${\Ll}$
be a principal $\Cx$--bundle on $Y$ (equivalently, a line bundle). The
group $\Cx$ acts naturally on the projective line $\Pp$ (preserving
two points, $0$ and $\infty$) and on the loop group $LG$ of maps from
the punctured formal neighborhood of $\infty \in \Pp$ to $G$. Let
$\Pp_\Ll = \Ll \underset{\C^\times}\times \Pp$ and $(LG)_\Ll = \Ll
\underset{\C^\times}\times LG$ be the induced bundles on $Y$. Note
that $\Pp_\Ll$ is equipped with two disjoint sections $0$ and
$\infty$.

Let $\Po$ be a $G$--bundle on $\Pp_\Ll$, equipped with a flag along
the section $\infty$ (i.e. reduction of the $G$--bundle $\Po|_\infty$
to $B$).  The sections of $\Po$ on the punctured formal neighborhood
of the section $\infty$ give rise to a principal bundle $\V$ on $Y$
for the $\Ll$--twist $(LG)_\Ll$ of $LG$ (clearly, the data of $\V$ and
$\Ll$ are equivalent to those of an $\wt{LG}$--torsor $\wt{\V}$ on
$Y$.)  Moreover the sections of $\Po$ over $\Pp_\Ll\sm\infty$ give
rise to a reduction $\V_-$ of $\V$ to the twist $(LG_-)_\Ll$, and
sections of $\Po$ near $\infty$ respecting the flag give rise to a
reduction $\V^+$ to $(LG^+)_\Ll$.  An affine oper is given by a
connection on $\V$ respecting $\V_-$ and having relative position
$[p_{-1}]$ with respect to $\V^+$.

Conversely, an affine oper on $Y$ gives rise to a $\C^\times$--bundle
$\Ll$ on $Y$ and a $G$--bundle $\Po$ on $\Pp_\Ll$ satisfying the above
conditions. In fact, $\Ll$ is nothing but the $\C^\times$--bundle
$\wt{\V} \underset{\wt{LG}}\times \C^\times$, where we use the natural
projection $\wt{LG} \to \C^\times$. Therefore $\Ll$ inherits a
connection from $\wt{\V}$. However, the relative position condition
puts strong constraints on the vector bundle $\V$ (as in
\cite{opers}). Namely, the components $f_i$ and $e_\theta\otimes z$ of
$p_{-1}=\ol{p}_{-1}+e_\theta\otimes z$ must transform as sections of
the tangent bundle of $Y$, hence $e_\theta$ transforms as a section of
$\Omega_Y^{\otimes h}$, where $h$ is the Coxeter number of
$G$. Therefore $z$ itself must transform as a section of
$\Omega^{\otimes(-1-h)}$. This implies the following

\subsubsection{Lemma.} The line bundle associated to the
$\C^\times$--bundle $\Ll$ is isomorphic to $\Omega_Y^{\otimes(-1-h)}$,
where $h$ is the Coxeter number of $G$.

\bigskip

Thus, there are no affine opers on a projective curve of genus other
than one, since $\Ll$ will then have nonzero degree.

\subsubsection{Definition.}    \label{genaffoper}
A {\em generic affine oper} is an affine oper, whose $\wt{LG}_-$ and
$\wt{LG}^+$ reductions are in general position (i.e., they correspond
to the open $\wt{LG}^+$--orbit of $LG_-\bs \wt{LG}$). Let
$\AOp^{\circ}(Y)$ be the open substack of generic affine opers in
$\AOp(Y)$.

\bigskip

The affine oper is generic when the bundle $\Po$ is trivial along the
fibers of $\Pp_\Ll\to Y$. In this case we can identify the
$G$--bundles $\Po|_\infty$ and $\Po|_0$ canonically.  The former comes
equipped with a flag, while the latter inherits a connection from the
affine oper connection on $\V_-$ (via the evaluation at $0$,
$(LG_-)_\Ll\to G$).

Thus, a generic affine oper gives rise to a $G$--bundle with a flag
and a connection, which are the data required for a
$G$--oper. Moreover, since the evaluation of $p_{-1}$ at $z=0$ gives
us $\ol{p}_{-1}$, we obtain a morphism $\AOp^{\circ}(Y) \to
\Op(Y)$. Namely, to an affine oper we assign the $G$--torsor $\E=\V_0$
with connection $\nabla$ and the $B$--reduction $\E_+$ coming from the
identification of $\V_0$ with $\V_\infty$. To see that the triple
$(\E,\nabla,\E_+)$ is an oper we must verify that the connection has
relative position $\Ofin$ with respect to $\E_+$. This follows from
the relation between the $B$--orbit $\Ofin \subset \g/\bb$ and the
$\wt{LG}^+$--orbit $\Oaff \subset L\g/L\g^+$. Under the standard
inclusion $\g\hookrightarrow L\g$, we have $\bb\hookrightarrow
L\g^+$. It follows from the explicit form of the orbits that
$\Oaff\subset \Ofin+(L\bb)_{\on{out}}$, where $(L\bb)_{\on{out}}$
consists of loops into $\bb$ which extend to $\Pp\sm\infty$.  Hence in
any trivialization of $\E_+$, evaluation at $0$ sends $\nabla$ to
$\Ofin$, so that $(\E,\E_+,\nabla)$ is indeed an oper.

\subsubsection{Proposition.}\label{affine oper to oper}
Let $Y$ be a curve with trivial canonical line bundle. Then the
canonical morphism $\AOp^{\circ}(Y)\to \Op(Y)$ is an isomorphism.

\subsubsection{Proof.} 
We construct the inverse morphism using a trivialization of the
tangent bundle, by reinterpreting the formula
$$p_{-1}=\ol{p}_{-1}+e_\theta\otimes z.$$ Given a $G$--oper
$(\E,\E_+,\nabla)$, let $\E^{\theta}$ denote the line subbundle of the
adjoint bundle $(\g)_\E$ corresponding to the highest weight line $\C
e_\theta\subset \g$ with respect to the Borel reduction
$\E_+$. Consider the $\Cx$--bundle $\Ll$ whose associated line bundle
is dual to the line bundle $\Omega_Y \otimes \E^{\theta} \simeq
\Omega_Y^{\otimes(1+h)}$. Let us pick a trivialization of the
canonical line bundle $\Omega$. Then we obtain a connection on $\Ll$.

As in \secref{geometric opers}, we associate to $\Ll$ a $\Pp$--bundle
$\Pp_\Ll$ on $Y$ and the corresponding twisted loop group scheme
$(LG)_\Ll$, which now inherits a connection from $\Ll$, and the
subgroups $(LG_-)_\Ll,(LG^+)_\Ll$. Let $\Po$ denote the $G$--bundle
on $\Pp_\Ll$ which is the pullback of $\E$ from $Y$, and
$(\V,\V_-,\V^+)$ the corresponding
$((LG)_\Ll,(LG_-)_\Ll,(LG^+)_\Ll)$--torsors, where the
$(LG^+)_\Ll$--reduction $V^+$ is defined so that the values of its
sections at $\infty$ lie in $\E_+$. Finally, the torsor $\V$ for the
group scheme with connection $(LG)_\Ll$ itself carries a connection
$\wt{\nabla}$ induced from $\nabla$. This connection needs to be
corrected to make $\V$ into an affine oper, because while $\V^+$ has
relative position $\Ofin\subset \g/\bb=L\g_+/L\g^+$ with respect to
$\wt{\nabla}$, we need relative position $\Oaff$.  But by
construction, there is a canonical section $e_\theta\otimes z$ of the
line bundle $\Ll\otimes\Omega_Y\otimes \E^{\theta}$, which is a
one--form on $Y$ with values in $\on{ad}\V_-$.  The connection
$\nabla$ obtained by adding this section to $\wt{\nabla}$ still
preserves $\V_-$, and it clearly has the correct relative position.

The above constructions work over an arbitrary base and are
functorial. Hence we obtain a morphism of stacks $\Op(Y) \to
\AOp^{\circ}(Y)$. It is easy to see that it is inverse to the morphism
$\AOp^{\circ}(Y)\to \Op(Y)$ constructed in \secref{genaffoper}. This
completes the proof.

\medskip

Now \lemref{automor} implies:

\subsubsection{Corollary.}    \label{automor1}
Let $Y$ be a curve with trivial canonical line bundle. Then the stack
$\AOp^\circ(Y)$ is the quotient of a scheme $AO^\circ(Y) \simeq Op(Y)$
by the trivial action of the center $Z(G)$ of $G$.

%\subsubsection{Proof.} Let $(\V,\nabla,\V_-,\V^+)$ be an affine
%oper. We need to show that the group of its automorphisms is
%isomorphic to $Z(G)$. These are the automorphisms of $\V$ that
%preserve the reductions $\V_-,\V^+$ and are flat with respect to
%$\nabla$. The last condition means that such an automorphism is
%uniquely determined by its action on the fiber of $\E$ at any point of
%$Y$. This action should preserve the intersection $\V_- \cap
%\V^+$. For generic affine opers, this intersection is a torsor over
%$LG_- \cap A_+ = Z(G)$, the center of $G$. Hence the group of
%automorphisms of an oper is isomorphic to $Z(G)$.
 
\subsection{Drinfeld--Sokolov Gauge.}

We now turn to our main objects of study: the affine opers on the
formal group $\At$. Recall that $\At$ is the formal additive group
with Lie algebra $\C p_{-1}$, where $p_{-1}$ is a fixed element of
$L\g$ given by formula \eqref{p-1} (i.e., from now on we fix our
choice of the generators $f_i, i=0,\ldots,\ell$). Therefore we have a
canonical coordinate $t$ on $\At$, such that $\At = \{ e^{tp_{-1}}
\}$, and a canonical vector field on $\At$ corresponding to $p_{-1}$
(also denoted by $p_{-1}$). Because of that, we can avoid twisting by the
group $\Cx$ that was used in the
general definition \ref{affoper1}.

Recall the notion of tautological relative position from
\defref{taut}. The following definition is equivalent to the
definitions \ref{affoper1}, \ref{genaffoper} in the case when $Y$ is
replaced by $\At$.

\subsubsection{Definition.} \label{t-line} An affine
oper on $\At$ is a quadruple $(\V,\nabla,\V_-,\V^+)$, where $\V$ is an
$LG$--torsor on $\At$ with a connection $\nabla$, a flat reduction
$\V_-$ to $LG_-$ and a reduction $\V^+$ to $LG^+$ in tautological
relative position with $\nabla$.

A generic affine oper on $\At$ is a quadruple $(\V,\nabla,\V_-,\V^+)$
as above such that the reductions $\V_-$ and $\V^+$ are in general
position.

The affine opers on $\At$ form a stack that is denoted by $\AOp(\At)$,
and generic affine opers form an open substack $\AOp^\circ(\At)$.

\subsubsection{} Suppose we are given a quadruple
$(\V,\nabla,\V_-,\V^{A_+}) \in \M_{p_{-1}}$, where $\V$ is an
$LG$--torsor on $\At$ with a flat connection $\nabla$, a flat
reduction $\V_-$ to $LG_-$ and a reduction $\V^{A_+}$ to $A_+$ in
tautological relative position with respect to $\nabla$. Then the
induced reduction of $\V$ to $LG^+$, $\V^+ = \V^{A_+} \times_{A_+}
LG^+$, is clearly in relative position $\Oaf \subset L\g/L\g^+$ with
respect to $\nabla$. Hence $(\V,\nabla,\V_-,\V^+)$ is an affine oper
on $\At$. This can certainly be done over an arbitrary base. Thus, we
obtain:

\subsubsection{Lemma.}    \label{prinmor}
There is a natural morphism of stacks $\M^{\Pp}_{p_{-1}} \to
\AOp(\At)$.

\subsubsection{Drinfeld--Sokolov Gauge.}
Recall that in \propref{pre-pre-oper} we established an isomorphism
between the abelianized Grassmannian $\Grp_A$, associated to the
principal Heisenberg subalgebra $\ab \subset L\g$, and
$\M^{\Pp}_{p_{-1}}$.  Our goal now is to prove that the above morphism
$\M^{\Pp}_{p_{-1}} \arr \AOp(\At)$, and hence the composition
$\Grp_A \to \AOp(\At)$, are in fact isomorphisms. In order to do
that we prove in this section a technical result, which shows the
existence of a canonical gauge for all affine opers on $\At$. This
result will enable us to show that any affine oper has a canonical
reduction to $A_+$, thus giving us an inverse morphism $\AOp(\At) \to
\M^{\Pp}_{p_{-1}}$. This gauge goes back to the original works on the
inverse scattering method in soliton equations. It was first
formulated in the language of connections by Drinfeld and Sokolov
\cite{DS}.

We will state our result in a much more general situation of a
strongly regular element $p$ of an arbitrary Heisenberg subalgebra
$\ab$ of $L\g$. (Recall \defref{strongly regular} of strongly regular
element $p \in \ab^{\geq l}$.) This statement is a direct
generalization of the Drinfeld-Sokolov lemma \cite{DS} (see also
\cite{dGHM}), and our proof essentially follows their argument.

\subsubsection{} Let us fix notation: $\ab\subset L\g$ is a Heisenberg
algebra with the canonical filtration $\{ \ab^{\geq j} \}$;
$p\in\ab^{\geq l}, l>0$; $\{L\g^{\geq j}\}$ is a filtration
compatible with $\ab$, in the sense of \secref{filtrations}.

Since $p \in \ab$, it is semisimple. If $p$ is regular, then we have
the decomposition $L\g\cong \on{Ker}(\on{ad}p)\oplus
\on{Im}(\on{ad}p)$. Furthermore, if $p$ is strongly regular, then
\lemref{strict compatibility} gives us the following
\subsubsection{Lemma.}    \label{compat}
$L\g^{\geq j} \equiv \ab^{\geq j}\oplus \on{ad}p\cdot
L\g^{\geq j-l}$ for every $j$.

\subsubsection{Proposition.}\label{DS gauge}
Let $p \in \ab^{\geq l}, l<0$, be a strongly regular element, and
$\nabla_t=\pa_t+p+q(t)$, with $q(t)\in L\g^+[[t]]$, be a connection on the
trivial $LG$--bundle on the formal disc $\Dt$. Then there is a gauge
transformation $M(t)\in LG^{>0}[[t]]$, such that
\begin{equation}\label{gauging}
M\inv(t) (\pa_t+p+q)M(t)=\pa_t+p+p_+(t),
\end{equation}
with $p_+(t)\in \ab^+[[t]]$. Furthermore, this equation determines $M(t)$
uniquely up to right multiplication by $A^+[[t]]$.

\subsubsection{Proof.}
To simplify notation, we will omit reference to the parameter $t$ on which 
all variables depend
(we write $\pa$ for $\pa_t$ and $\nabla$ for $\nabla_t$ etc).

Let us solve \eqref{gauging} for some $M\in LG^{>0}$ and $p_+\in
\ab^+$ . Since $LG^{>0}$ is a prounipotent group,
every $M\in LG^{>0}$ may be written as $\exp{m}$ for some $m\in
L\g^{>0}$. We construct $m$ and $p_+$ by induction on the filtration
as $m=\sum_{i=l}^{\infty} m_i$ and $p_+=\sum_{i=1}^{\infty} p_i$ with
$m_i,p_i\in L\g^{\geq i}$.  Write $$M\inv \nabla M=\sum_{k\geq
0}\frac{1}{k!}\on{ad}^k m \cdot \nabla.$$ At the first step, $\nabla
\equiv \pa +p (\on{mod}L\g^+)$ and so $m\equiv 0\,(\on{mod} L\g^{\geq
-l})$ indeed works modulo $L\g^+$.

For the inductive step, suppose that $k\geq -l$, and
$M_k=\on{exp}(\sum_{i=-l}^k m_i)$ satisfies

$$M_k\inv\nabla M_k\equiv \pa+p+\sum_{i=1}^{k-l-1} p_i\hskip.3in
(\on{mod} L\g^{\geq k+l}).$$

Denote $\delta_k = M_k\inv\nabla M_k-(\pa+p+\sum_{i=1}^{k+l-1}
p_i)$, so that $\delta_k \in L\g^{\geq k+l}$. By \lemref{compat},
we may decompose $\delta_k=p_{k+l}+[p,m_k]$ for some
$p_{k+l}\in \ab^{\geq k+l}$ and $m_k\in L\g^{\geq k}$.  If we set
$M_{k+1}=\on{exp}(\sum_{i=-l}^{k+1} m_i)$, then
$$M_{k+1}\inv \nabla M_{k+1}\equiv \pa+p+\sum_{i=1}^{k+l}
p_k\hskip.3in (\on{mod}L\g^{\geq k+l+1}),$$
completing the inductive step.

Note that at each step, $M_k$ is unique up to the right multiplication
by an element of $A^+$. This completes the proof.

\subsubsection{} \propref{DS gauge} implies the following statement
for the principal Heisenberg,
whose proof we postpone to \propref{general DS reduction}, where 
the case of general Heisenbergs is taken up.

\subsubsection{Corollary.}\label{DS reduction}
Let $(\V,\nabla,\V_-,\V^+)$ be an affine oper on $\At$. Then $\V$ has
a unique reduction $\V^{A_+}\subset\V^+$ to $A_+$, such that
\begin{enumerate}
\item[(1)] the induced $A$--torsor $\V^A\subset\V$ is flat;
\item[(2)] $\V^{A_+}$ has the tautological relative position with
respect to 
$\nabla$.
\end{enumerate}

\subsubsection{Remark.} The above corollary may be applied to the
variants of affine opers defined on arbitrary differential
schemes. Thus we obtain a morphism from affine opers on a differential
scheme to the bundles considered in \propref{differential opers}, and
thus to $\Gr_A$.

\subsubsection{Remarks: The Extended Oper Family.}\label{extended family}
The direct analog of the notion of an oper for the loop group would
not require a flat reduction to $LG_-$. The resulting objects relate
to pseudo--differential operators. We only note here that the
Drinfeld--Sokolov lemma applies to them as well.

We may also replace the curve $\Pp$ in the above definition by an
arbitrary smooth curve $X$ (in particular, replacing the
$LG_-$--reduction by an $LG_-^X$--reduction). The resulting
$X$--affine opers are analogous to elliptic sheaves for $X$. Moreover
we have a notion of ``generic'' $X$--affine oper, akin to the
genericity conditions on line bundles in the Krichever
construction. Choosing a basepoint $0\in X\sm \infty$ we obtain a map
from generic affine opers for any $X$ to opers on the disc.  Other
generalizations of affine opers will be discussed in \secref{general
affine opers}.

We briefly mention some other relatives of opers. The difference opers
and Frobenius opers arise from replacing differential operators by
difference operators and polynomials in the Frobenius, respectively.
They may be expected to be useful in the theory of the $q$--KdV
hierarchies and elliptic sheaves, respectively.  Recall from
\secref{difference} that there is a notion of relative position for
equivariant bundles with reductions.  Let $(Y,l)$ be a pair consisting
of a scheme $Y$ and an automorphism $l$ of $Y$. Let $\chi$ be an
$LG^+$--orbit in $LG/LG^+$. Recall that these orbits are in
one--to--one correspondence with the affine Weyl group. In
applications one usually takes $\chi$ to be a Coxeter element.

\subsubsection{Definition.} A difference oper on $(Y,l)$ is a
quadruple $(\V,l_{\V},\V_-,\V^+)$, where $\V$ is an $LG$--torsor on
$Y$, $l_{\V}$ a lifting of $l$ to $\V$, $\V_-$ an $LG_-$--reduction of
$\V$ preserved by $l_{\V}$, and $\V^+$ an $LG^+$--reduction of $\V$ in
relative position $\chi$ with respect to $l_{\V}$.

\subsubsection{} If $Y$ is a scheme defined over a finite field, we
take $l$ to be the Frobenius of $Y$. The resulting objects, the
Frobenius opers, are generalizations for semisimple groups $G$ of the
notion of elliptic sheaves (see e.g. \cite{new}). As was the case for
their analogs, Krichever sheaves from \secref{Krichever sheaves}, the
definition of elliptic sheaves for $X$ as sheaves on $X\times Y$ is
simplified in the semisimple case: by the Drinfeld--Simpson theorem
\cite{DSi} we may replace $G$--torsors on $X\times Y$ by $LG$--torsors
on $Y$, with reductions as above.

\section{Integrable Systems.}    \label{Main}

In this section we combine the results of sections \secref{Higgs} and
\secref{opers} to obtain isomorphisms between abelianized
Grassmannians and the moduli of affine opers. Under this isomorphism,
the action of the group $\Ahat$ gives rise to a collection of
infinitely many commuting flows on the space of affine opers. They
form a generalized Drinfeld-Sokolov hierarchy. We first work out in
detail the case of the principal Heisenberg subalgebra with $\Pp$ as
the underlying curve $X$. After that we generalize the construction to
the case of an arbitrary Heisenberg subalgebra.

\subsection{The Principal Case: KdV}    \label{prinh}

Let us recall the notation: $\infty$ is a point on $\Pp$, $LG$ is the
loop group associated to the disc at $\infty$, and $L\g$ its Lie
algebra. In this subsection we denote by $\ab$ the principal
Heisenberg subalgebra of $L\g$, and by $p_{-1}$ a generator of the
one-dimensional space $\ab_{\geq -1}/\ab_+$. Finally, $\At$ is the
formal additive subgroup of $A/A_+$ generated by $p_{-1}$.

\subsubsection{Theorem.}\label{tripartite theorem}
\begin{enumerate}
\item[(1)] There is a canonical isomorphism between the abelianized
Grassmannian $\Grp_A$ and the moduli stack $\AOp(\At)$ of affine
opers on $\At$.

\item[(2)] The isomorphism of part (1) identifies the big cell
$\Gro_A$ of $\Grp_A$ and the moduli stack $\AOp^\circ(\At)$ of
generic affine opers on $\At$.

\item[(3)] For each point $0 \in \Pp\sm\infty$, there is a canonical
isomorphism between $\Gro_A$ and the moduli stack $\Op(\At)$ of
$G$--opers on $\At$.
\end{enumerate}

\subsubsection{Proof.} 
The construction of \corref{DS reduction} gives us a morphism
$\AOp(\At) \arr \M^{\Pp}_{p_{-1}}$, which sends $(\V,\nabla,\V_-,\V^+)
\in \AOp(\At)$ to $(\V,\nabla,\V_-,\V^{A_+}) \in
\M^{\Pp}_{p_{-1}}$. On the other hand, in \lemref{prinmor} we
constructed a morphism $\M^{\Pp}_{p_{-1}} \arr \AOp(\At)$.  It is
clear from the construction that these morphisms are inverse to each
other. Hence we obtain an isomorphism $\M^{\Pp}_{p_{-1}} \simeq
\AOp(\At)$. But $\Grp_A$ is isomorphic to $\M^{\Pp}_{p_{-1}}$ by
\propref{pre-pre-oper}. Therefore $\Grp_A \simeq \AOp(\At)$. This
proves part (1).

By definition, the big cell $\Gro_A$ of $\Grp_A$ classifies
$LG$--torsors $\V$ with reductions $\V_-$ to $LG_-$ and $\V^{A_+}$ to
$A_+$, such that $\V_-$ and the induced $LG^+$--reduction $\V^{A_+}
\times_{A_+} LG^+$ are in general position (see \secref{Gro}). On the
other hand, $\AOp^\circ(\At)$ classifies the quadruples
$(\V,\nabla,\V_-,\V^+)$, such that the reductions $\V_-$ and $\V^+$
are in general position. Hence under the above isomorphism $\Gro_A$ is
mapped to $\AOp^\circ(\At)$ and we obtain part (2).

Finally, part (3) follows from part (2) and \propref{affine oper to
oper}.

\subsubsection{Remark.} Let us once again spell out the definition of
the morphisms $\Grp_A \arr \AOp(\At)$ and $\AOp(\At) \arr \Grp_A$.

Let $(\V,\V_-,\V^{A_+})$ be a point of $\Grp_A$. Here $\V$ is an
$LG$--torsor, and $\V_-$, $\V^{A_+}$ are its reductions to $LG_-$, $A_+$,
respectively. Let $\T$ be the universal bundle over $\Grp_A$, whose
fiber at $(\V,\V_-,\V^{A_+})$ is $\V$; it has canonical reductions
$\T_-$ and $\T^{A_+}$. Denote by $\pi: \At \arr \Grp_A$ the map
corresponding to the action of $\At$ on $(\V,\V_-,\V^{A_+})$. Note that
the action of $\At$ lifts to $\T$. Now as in
\lemref{action->connection}, pulling back $\T$ to $\At$ by $\pi$, we
obtain an $LG$--bundle $\V_t$ on $\At$ with a flat connection $\nabla$
and the reductions $\V_{-,t}$ and $\V^{A_+}_t$ to $LG_-$ and
$A_+$. Moreover, by \lemref{tautological position}, $\V_{-,t}$ is
preserved by $\nabla$, while $\V^{A_+}_t$ is in tautological relative
position with respect to $\nabla$. Denote by $\V^+_t$ the induced
$LG^+$--reduction $\V^{A_+}_t \times_{A_+} LG^+$. Then
$(\V_t,\nabla,\V_{-,t},\V^+_t)$ is the affine oper corresponding to
$(\V,\V_-,\V^{A_+}) \in \Grp_A$.

The inverse map is constructed as follows. Given an affine oper
$(\V_t,\V_{-,t},\V^+_t)$, we obtain a canonical $A_+$--reduction
$\V^{A_+}_t$ of $\V_t$ using \corref{DS reduction}. Then we look at
the fiber $\V_0$ of $\V_t$ at $0 \in \At$. It comes with the
reductions $\V_{-,0}$ and $\V^{A_+}_0$ to $LG_-$ and $A_+$,
respectively. Hence we attach to our affine oper the point
$(\V_0,\V_{-,0},\V^{A_+}_0)$ of $\Grp_A$.

The above construction is almost tautological. However, there is one
place where in addition to the general functorial correspondence of
\propref{abstract isomorphism} we really have to use the specifics of
our situation: namely, when we switch between reductions to $A_+$ and
$LG^+$. It is easy to pass from an $A_+$--reduction to an
$LG^+$--reduction by using induction, and this allows us to pass from
the unwieldy moduli space $\M^{\Pp}_{p_{-1}}$ to the much nicer moduli
of affine opers. But {\em a priori} one can not go back from an
$LG^+$--reduction to an $A_+$--reduction. For this we need to rely on
the rather technical Drinfeld-Sokolov gauge (see \corref{DS
reduction}). As we will see below, this construction works if we
replace $p_{-1}$ (the generator of the formal group $\At$) by any
strongly regular element $p$ of a general Heisenberg subalgebra $\ab$.

\subsubsection{Geometric interpretation.}\label{geometric interp} The
morphism $\Grp_A \arr \AOp(\At)$ has a simple geometric
interpretation, which is close in spirit to the Krichever construction
as in \cite{Mum} (see also \secref{geometric opers}, \secref{Krichever
sheaves}). To simplify the picture, we explain it using analytic
rather than algebro--geometric language.

Let $\E$ be a $G$--bundle on $\Pp$, equipped with a reduction to $A_+$
on $\D$. We wish to deform $\E$ by the action of $p_{-1}$, as we did
in the trivial abelian case \secref{Fourier example}. Thus we
construct a $G$--bundle $\wt{\E}$ on $\Pp\times\At$ by multiplying the
transition function of $\E$ on $\Dx$ by $e^{-tp_{-1}}$. This change in
the transition function does not affect $\E$ away from $\infty$, so
that we can canonically identify sections $\E(\Pp\sm\infty)$ at time
$t=0$ with the sections $\E_t(\Pp\sm\infty)$ of its deformation at any
time $t\in\At$. In other words, we have a canonical (flat) partial
connection $\nabla$ over $\At$ on
$\wt{\E}|_{(\Pp\sm\infty)\times\At)}$ in the direction of $\At$.

Thus the restriction $\wt{\E}_0$ of $\wt{\E}$ to $0 \times \At$,
where $0 \in \Pp\sm\infty$, is a $G$--bundle on $\At$ with a flat
connection. On the other hand, the restriction $\wt{E}_\infty$ of
$\wt{\E}$ to $\infty \times \At$ is a $G$--torsor on $\At$ with a
$B$--reduction induced by the $A_+$--reduction of $\E|_{\D}$. The
point is that when $\E$ is a trivial bundle on $\Pp$, we may identify
$\wt{\E}_0$ with $\wt{E}_\infty$. Hence we obtain a $G$--bundle on
$\At$ with a flag and a flat connection. Explicit calculation (see
below) shows that it is a $G$--oper. Thus, we attach to $\E$ a
$G$--oper on $\At$.

On the other hand, we can consider the $LG$--bundle on $\At$
corresponding to taking the sections of $\wt{E}$ over $\Dx$. It
carries a reduction to $A_+$, and hence the induced reduction to
$LG^+$. But this reduction is not preserved by the connection
$\nabla$. The particular form of $p_{-1}$ shows that in fact we obtain
an affine oper.

Conversely, an affine oper gives rise to a period map $\At \arr
LG_-\bs LG/LG^+$, which by the (Griffiths) transversality of the
connection is tangent to a certain completely non-integrable
distribution (compare \cite{Mum}). The Drinfeld--Sokolov gauge picks
out a canonical lifting of this period map to $LG_-\bs LG/A_+$ which
is tangent to the vector field $p_{-1}$, and from this data we recover
our original $\E$.

\subsubsection{} According to \lemref{silly}, $\Gro_A$ is the quotient
of the scheme $Gr^\circ_A= G\bs LG_+/A'_+$ by the trivial action of
$A_0=Z(G)$. On the other hand, $\AOp^\circ(\At)$ (resp., $\Op(\At)$)
is the quotient of a scheme $AO^\circ(\At)$ (resp.,
$Op(\At)=Op_{\g}(\At)$) by the trivial action of $Z(G)$, see
\lemref{automor1} and \lemref{automor}. Hence we obtain:

\subsubsection{Corollary.}    \label{iso-for-opers}

\begin{enumerate}
\item[(1)] There is a canonical isomorphism of schemes $Gr^\circ_A \simeq
AO^\circ(\At)$.

\item[(2)] For each point $0 \in \Pp\sm\infty$, there is a canonical
isomorphism $Gr^\circ_A \simeq Op(\At)$.
\end{enumerate}

\subsubsection{The Universal Oper.} 
Let $\T$ be the tautological $G$--bundle over
$Op(\At) \times \At$, whose restriction to
$(\E,\nabla,\E_B) \times \At$ is $\E$. The two additional structures
of an oper translate into a $B$--reduction $\T_B$ of $\T$ and a
partial connection $\nabla_{\T}$ along $\At$. This is the ``universal
oper'' on $Op(\At) \times \At$.

On the other hand, consider the homogeneous space $LG_+/A'_+$, which is
a $G$--bundle over $Gr^\circ_A = G\bs LG_+/A'_+$. The subgroup $LG^+
\subset LG_+$, which is canonically associated to $A_+$ gives us a
$B$--reduction $LG^+/A'_+$ of $LG_+$. Now fix a point $0 \in \Pp$
(i.e., choose a generator $z$ of $\C[\Pp\sm\infty]$ up to a scalar
multiple) and the corresponding subgroup $LG_{<0}$ of $LG$. Then we
can view $LG_+$ as an open part of $LG_{<0}\bs LG$. Hence we have an
action of $p_{-1}$ from the right on $LG_+/A'_+$. Define a partial
connection $\nabla_{p_{-1}}$ on $(LG_+/A'_+) \times \At$ along $\At$ by
the formula $\pa_t + p_{-1}$. \corref{iso-for-opers},(2) can be
interpreted as follows:

\subsubsection{Corollary.}
The universal oper $(\T,\nabla_{\T},\T_B)$ on $Op(\At) \times \At$ is
canonically isomorphic to the triple $(LG_+/A'_+,\pa_t +
p_{-1},LG^+/A_+)$ on $(G\bs LG_+/A'_+) \times \At$.

\subsubsection{Differential Polynomials.} 
According to \lemref{canreps}, $Op(\At)$
is isomorphic to the pro-vector space $V[[t]]$, where $V$ is a
subspace of $\n_+$ satisfying $\n_+ = V \oplus \on{Im} \on{ad}
p_{-1}$. Let us fix such a subspace $V$, and a homogeneous basis of
$V$ with respect to the principal gradation. Then we can identify
$V[[t]]$, and hence $Op(\At)$, with the pro-vector space of
$\ell$--tuples $(v_i(t))_{i=1,\ldots,\ell}$ of formal Taylor power
series.  For instance, in the case $\g={\mathfrak s}{\mathfrak l}_2$
we identify $Op(\At)$ with the space of operators of the form
\begin{equation}    \label{sl2oper}
\pa_t + \begin{pmatrix} 0 & v(t) \\ 1 & 0 \end{pmatrix}.
\end{equation}

Let $v_i^{(n)}$ be the linear functional on $V[[t]]$, whose value at
$(v_i(t))_{i=1,\ldots,\ell}$ is $\pa_t^n v_i(t)|_{t=0}$. We can now
identify the ring of regular functions $\C[Op(\At)]$ on $Op(\At)$ with
the ring of differential polynomials
$\C[v_i^{(n)}]_{i=1,\ldots,\ell;n\geq 0}$. Furthermore, the natural
action of $\pa_t$ on $\C[Op(\At)]$ is given by $\pa_t \cdot v_i^{(n)}
= v_i^{(n+1)}$. We will use the notation $v_i$ for $v_i^{(0)}$.

\subsubsection{} Now
\corref{iso-for-opers} gives us a new proof of the following result,
which is equivalent to Theorem 4.1 from \cite{Five} (to make the
connection with \cite{Five} clear, note that $Gr^\circ_A \simeq N_+\bs
LG^{>0}/A'_+$).

\subsubsection{Theorem}    \label{five}
The ring of functions $\C[Gr^\circ_A]$ on $Gr^\circ_A$ is
isomorphic to the ring of differential polynomials
$\C[v_i^{(n)}]_{i=1,\ldots,\ell;n\geq 0}$. Under this isomorphism, the
right infinitesimal action of $p_{-1}$ on $Gr^\circ_A$ is given by
$p_{-1} \cdot v_i^{(n)} = v_i^{(n+1)}$.

\subsubsection{KdV hierarchy.}    \label{kdv}
The infinite-dimensional formal group $A/A_+$ and its Lie algebra
$\ab/\ab_+$ act from the right on $Gr^\circ_A$.  Hence by
\thmref{five} we obtain an action of $\ab/\ab_+$ on the space
$Op(\At)$ of $G$--opers on $\At$. The action of the element $p_{-1}$
of $\ab/\ab_+$ coincides with the flow generated by $\pa_t$. Other
elements of $\ab/\ab_+$ act on $Op(\At)$ by vector fields commuting
with $\pa_t$. It is known that $\ab/\ab_+$ has a basis $p_i, i \in
-I$, where $I$ is the set of all positive integers equal to the
exponents of $\g$ modulo the Coxeter number. The degree of $p_i$
equals $i$ with respect to the principal gradation of $L\g$ (see
\cite{Kac1}).

Given an element $p_{-m} \in \ab/\ab_+, m \in I$, let $\wt{p}_{-m}$ be
the corresponding derivation of $\C[v_j^{(n)}]_{j=1,\ldots,\ell;n\geq
0}$. We know that each $\wt{p}_{-m}$ commutes with $\wt{p}_{-1} =
\pa_t$, and hence is an {\em evolutionary} derivation. Because of the
Leibnitz rule, the action of an evolutionary derivation on
$\C[v_i^{(n)}]_{i=1,\ldots,\ell;n\geq 0}$ is uniquely determined by
its action on $v_i, i=1,\ldots,\ell$.

We know that $\wt{p}_{-m} \cdot v_i$ is a differential polynomial in
$v_i$'s. The system of partial differential equations
\begin{equation}    \label{KdVeq}
\pa_{t_m} v_i = \wt{p}_{-m} \cdot v_i, \quad \quad
i=1,\ldots,\ell,
\end{equation}
(considered as equations on the functions $v_i(t),
i=1,\ldots,\ell,$ belonging to some reasonable class of functions), is called
the $m$th equation of the {\em generalized KdV hierarchy corresponding
to $\g$}, and the time $t_m$ is called the $m$th time of the
hierarchy. The totality of these equations as $m$ runs over $I$ is
called the generalized KdV hierarchy corresponding to $\g$.

For instance, for $\g={\mathfrak s}{\mathfrak l}_2$, we obtain the KdV
hierarchy. In this case, $I$ consists of all positive odd integers. We
already know the action of $p_{-1}$: it is given by $\pa_t$. In
particular, the equations \eqref{KdVeq} read in this case: $\pa_{t_1}
v = \pa_t v$ (here we use $v$ for $v_1$), which means that $t_1=t$
(this is true for all $\g$). The next element is $p_{-3}$. The
corresponding derivation $\wt{p}_{-3}$ has been computed explicitly in
\cite{Five}. The resulting equation is
\begin{equation}    \label{kdv3}
\pa_{t_3} v = \frac{3}{2} v \pa_t v - \frac{1}{4} \pa_t^3 v,
\end{equation}
which is the KdV equation (up to a slight redefinition of variables).

\subsection{Zero curvature representation.}
It is well-known that the KdV equations can be written in the zero
curvature, or Zakharov--Shabat, form (see \cite{DS}). The zero
curvature formalism is one of the standard methods to write down
equations of completely integrable systems, and it is convenient for
explicit description of the associated Hamiltonian structures. This
form of the equations arises very naturally in our approach as the
equations expressing the flatness of the connection on the formal
group $\Ahat$. Recall that we have {\em identified} $Gr^\circ_A$ with
the moduli space of flat connections on the entire formal group
$A/A_+$ (see \propref{pre-pre-oper} for the precise statement). The
flatness condition, written in explicit coordinates, takes the
familiar form of zero curvature equations as we will now demonstrate.

\subsubsection{}    \label{zcf}
The isomorphism of \propref{pre-pre-oper} assigns to each point $K_+
\in G\bs LG_+/A'_+$ $= Gr^\circ_A$ a quadruple
$(\V,\nabla,\V_-,\V^{A_+})$, where $\V$ is an $LG$--torsor on $\Ahat$
with a flat connection $\nabla$, a flat $LG_-$--reduction $\V_-$, and
an $A_+$--reduction in tautological relative position with
$\nabla$. We want to trivialize $\V$ and calculate the contraction
$\nabla_p$ of $\nabla$ with the vector field on $\Ahat$ coming from
the left action of $p \in \ab/\ab_+$.

It is convenient to trivialize $\V$ in such a way that the
$LG_-$--reduction $\V_-$ and the induced $LG^+$--reduction $\V^{A_+}
\times_{A_+} LG^+$ are preserved. Since we started with a point on the
big cell $Gr^\circ_A$, these reductions are in general
position. Therefore such a trivialization is unique up to the gauge
action of the group $LG_- \cap LG^+ = B_+$. 
Let us first choose one
such trivialization using the factorization of loops as in \secref{the
big cell}. The computation of the connection operator $\nabla_p$ in
this trivialization is the content of \lemref{trivial
connection}, and we obtain:
$$
\nabla_p = \partial_t + (K_+(t) p K_+(t)\inv)_-, \quad \quad p \in
\ab/\ab_+.
$$
But the connection $\nabla$ is flat by our construction. Therefore
$[\nabla_p,\nabla_{p'}]=\nabla_{[p,p']}=0$ for all $p,p' \in
\ab/\ab_+$. This gives us the zero curvature equations
\begin{equation}    \label{zckdv}
\left[ \partial_{t_m} + (K_+(t) p_{-m} K_+(t)\inv)_-,\pa_{t_n} +
(K_+(t) p_{-n} K_+(t)\inv)_- \right] = 0,
\end{equation}
where $m,n \in I$.

In the special case $n=1$ we obtain the equation
\begin{equation}    \label{zckdvm}
\left[ \partial_{t_m} + (K_+(t) p_{-m} K_+(t)\inv)_-,\pa_{t} +
(K_+(t) p_{-1} K_+(t)\inv)_- \right] = 0.
\end{equation}
But the component of our connection in the direction of $p_{-1}$ is
an affine oper:
\begin{equation}    \label{affoper}
\pa_t + (K_+(t) p_{-1} K_+(t)\inv)_- = \pa_t + p_{-1} + {\mathbf
b}(t),
\end{equation}
where ${\mathbf b}(t) \in \bb_+[[t]]$. Denote $L_m = (K_+(t) p_{-m}
K_+(t)\inv)_-$. It is easy to see that any solution for $L_m$ as an
element of $\g[[t]]$ is a differential polynomial in the matrix
elements ${\mathbf b}_\al(t), \al \in \Delta_+$ of ${\mathbf
b}(t)$. Hence formula \eqref{zcm} expresses $\pa_{t_m} {\mathbf
b}_\al(t)$ as a differential polynomial in ${\mathbf b}_\al(t)$'s.

However, we should remember that the formulas for the connection
operators that we obtained refer to a particular trivialization of
$\V$, which was unique up to the gauge action of $B_+[[t]]$. 
We may require the trivialization to preserve not only $\V_-$ and
$\V^+$, but also the form \eqref{affoper} -- in other words to have
symbol $p_{-1}$ -- thus further reducing the gauge freedom to
$N_+[[t]]$. 
The
equations \eqref{zckdvm} are invariant under this action, and we
should consider them as equations not on the space of connections of
the form \eqref{affoper}, but on the space of affine opers, which is
its quotient by the free action of $N_+[[t]]$. Under this action, we
can bring \eqref{affoper} to the form
$$
\pa_t + p_{-1} + {\mathbf v}(t), \quad \quad {\mathbf v}(t) \in
V[[t]],
$$
where $V$ is the transversal subspace to $\on{ad} \ol{p}_{-1}$ in
$\n_+$ (see \lemref{canreps}). If we choose a basis $\{ v_i
\}_{i=1,\ldots,\ell}$, we can interpret the equation \eqref{zckdvm} as
an equation expressing $\pa_{t_m} v_i(t)$ as a differential
polynomial in $v_j$'s. Thus we obtain the {\em zero curvature
representation} of the $m$th generalized KdV equation
\eqref{KdVeq}. This form of the generalized KdV equations was
first introduced in \cite{DS}.

\subsection{Flag Manifolds and the mKdV Hierarchy.}    \label{flags}

\subsubsection{}

Recall from \secref{homogeneous filtration} that for the homogeneous
filtration of the loop group, we were able to define one piece
$LG_-=LG_{\leq 0}$ of an opposite filtration, just using the global
curve $\Pp$. In order to find a similar partial splitting of the
principal filtration, let us fix a Borel subgroup $B_-\subset G$ which
is transverse to the Borel subgroup $B$ defined by $A_+$ at $\infty$
(in other words, $[B_-]\in G/B$ lies in the open $B$--orbit).

Let us fix a point $0 \in \Pp\sm\infty$. Define $LG^-=LG^{\leq
0}\subset LG$ as the subgroup that consists of loops $x$ which extend
to all of $\Pp\sm\infty$ (that is, $x\in LG_-$) whose value at $0$
lies in $B_-$. Thus $LG^-$ is a ``lower Iwahori subgroup'' of
$LG$. (Replacing $B_-$ by a parabolic, one obtains ``parahori subgroups'' of
$LG$.)

\subsubsection{Definition.} The {\em affine flag manifold $\Fl$} is the
scheme of infinite type representing the moduli functor of
$G$--torsors on $\Pp$, equipped with a trivialization on $\D$ and a
reduction to $B_-$ at $0$.

The following proposition is proved in the same way as
\propref{Gr-iso}.

\subsubsection{Proposition.} $\Fl\cong LG^-\bs LG$.

\subsubsection{} The analog of abelianization of $\Gr$ 
for the affine flag manifold is the stack $\Fl_A$ classifying
$G$--torsors $\V$ on $\Pp$, equipped with a reduction $\E^{A_+}$ of
the $LG_+$--torsor $\E|_{\D}$ to $A_+$, and a reduction to $B_-$ at
$0$.

The group $LG^{>0}$ acts on $\Fl$ from the right.  The infinitesimal
decomposition $L\g\cong L\g^-\oplus L\g^{>0}$ implies that the orbit
of the identity coset is open, and since $LG^-\cap LG^{>0}=1$, the
orbit is in fact isomorphic to $LG^{>0}$. This orbit is denoted by
$\Flo$. There is an obvious map $\Fl\twoheadrightarrow \Grp$,
forgetting the flag at $0$. The fibers of this map are isomorphic to
the flag manifold $G/B$ of $G$, and $\Flo$ is the inverse image of the
big cell $\Gro$ under this map. The restriction $\Flo \arr \Gro$ is an
$N$--bundle, whose fibers are identified with the big cell of
$G/B$.

The image of $\Flo$ in $\Fl_A$ is an open substack of $\Fl_A$,
which we denote by $\Flo_A$. In the same way as in the case of
abelianized Grassmannians, one shows that $\Flo_A$ is the quotient
of a scheme of infinite type $Fl^\circ_A = LG^{>0}/A'_+$ by
the trivial action of $Z(G)$. The natural morphism $Fl^\circ_A
\arr Gr^\circ_A$ is again an $N$--bundle.

The following are generalizations of the notions of opers and affine
opers for the flag manifold. Note that though one can define affine
Miura opers on an arbitrary curve $Y$ as in \defref{affoper1}, we only
need them when $Y$ is $\At$. Therefore the definition
below will suffice for our current purposes.

\subsubsection{Definition.} A {\em Miura $G$--oper} on a curve $Y$ is a
quadruple $(\E,\nabla,\E_B,\E'_B)$, where $\E$ is a $G$--torsor on $Y$
with a connection $\nabla$ and two reductions: $\E_B, \E'_B$, to the
Borel subgroup $B$ of $G$, which are in generic position. The
reduction $\E'_B$ is preserved by $\nabla$, while $\E_B$ is in
relative position ${\mbf O}$ with respect to $\nabla$.

An {\em affine Miura oper} is a quadruple $(\V,\nabla,\V^-,\V^+)$,
where $\V$ is an $LG$--torsor on $\At$ with a connection $\nabla$, a
flat reduction $\V^-$ to $LG^-$ and a reduction $\V^+$ to $LG^+$ in
tautological relative position with $\nabla$.

A {\em generic affine Miura oper} on $\At$ is a quadruple
$(\V,\nabla,\V^-,\V^+)$ as above, such the reductions $\V^-$ and
$\V^+$ are in generic position.

The Miura $G$--opers on $\At$ form a stack $\MOp(\At)$.  The affine
opers on $\At$ form a stack that is denoted by $\AMOp(\At)$, and
generic affine opers form its open substack $\AMOp^\circ(\At)$.

We have obvious surjective morphisms $\MOp(\At) \arr \Op(\At)$ and 
$\AMOp(\At) \arr \AOp(\At)$.

The following lemma is proved in the same way as \propref{affine oper
to oper}.

\subsubsection{Lemma.} There is a canonical isomorphism
$\AMOp^\circ(\At) \simeq \MOp(\At)$, which makes the following
diagram commutative.
\setlength{\unitlength}{1mm}
\begin{equation}
\begin{CD}
\AMOp^\circ(\At)  @>{\sim}>> \MOp(\At)
\end{CD}
\end{equation}
\begin{picture}(60,0)(0,0)
\put(52,0){$\Big\downarrow$}
\put(52,1){$\Big\downarrow$}
\put(88,0){$\Big\downarrow$}
\put(88,1){$\Big\downarrow$}
\end{picture}
$$\begin{CD}
\AOp^\circ(\At)\;\;   @>{\sim}>>\;\; \Op(\At)
\end{CD}$$

\noindent Here the isomorphism $\AOp^\circ(\At) \simeq \Op(\At)$ from
\lemref{affine oper to oper} corresponds to the point $0$ used in the
definition of the subgroup $LG^-$.

\subsubsection{}
The stack $\AMOp^\circ(\At)$ (resp., $ \MOp(\At)$) is the quotient of
an affine scheme of infinite type $AMO^\circ(\At)$ (resp., $MOp(\At)$)
by the trivial action of $Z(G)$.

The following statement is proved in the same way as
\corref{iso-for-opers}.

\subsubsection{Theorem.}
$Fl^\circ_A \simeq AMO^\circ_{-1}$ and $Fl^\circ_A \simeq MOp(\At)$.

\subsubsection{}    \label{coord}
Let us describe explicitly the ring of functions on $MOp(\At)$. Let
$(\E,\nabla,\E_B,\E'_B)$ be a Miura oper on $\At$. Since the
$B$--reductions $\E_B$ and $\E'_B$ are in generic position, they
produce a 
unique compatible reduction $\E_H$ of $\E$ to $H$.
Let us trivialize $\E_H$. The operator of connection $\nabla$ reads
relative to this trivialization as follows:
$$
\pa_t + \sum_{i=1}^\ell \phi_i(t) \cdot f_i + {\mathbf h}(t), \quad
\quad {\mathbf h}(t) \in \h[[t]].
$$
By changing trivialization of $\E_H$ (i.e., applying a gauge
transformation from $H[[t]]$), it can be brought to the form
\begin{equation} \label{Miura oper}
\pa_t + \ol{p}_{-1} + {\mathbf h}(t), \quad \quad {\mathbf h}(t) \in
\h[[t]].
\end{equation}
%For $G=SL_n$ this has the form
%\begin{equation*}
%\partial_t+\left( \begin{array}{ccccc}
%u_1 &0& 0&cdots&0\\
%1& u_2&0&cdots&0\\
%0&1&u_3&\cdots&0\\
%\vdots&\ddots&\ddots&\ddots&\vdots\\
%0&0&\cdots&1&u_n
%\end{array}\right)
%\end{equation*}
%with the sum of $u_i$ zero.

Denote $u_i(t) = \al_i({\mathbf h}(t))$. We can now identify
$AOp(\At)$ with the pro-vector space of $\ell$--tuples
$(u_i(t))_{i=1,\ldots,\ell}$ of formal Taylor series. For instance, in
the case $\g={\mathfrak s}{\mathfrak l}_2$, we identify $AOp(\At)$
with the space of operators of the form
$$
\pa_t + \begin{pmatrix} \frac{1}{2} u(t) & 0 \\ 1 & - \frac{1}{2} u(t)
\end{pmatrix}.
$$

Let $u_i^{(n)}, i=1,\ldots,\ell; n\geq 0$, be the function on
$AOp(\At)$, whose value at $(u_i(t))_{1,\ldots,\ell}$ equals
$\pa_t^n u_i(t)|_{t=0}$. The ring of functions $\C[AOp(\At)]$ on
$AOp(\At)$ is isomorphic to the polynomial ring
$\C[u_i^{(n)}]_{i=1,\ldots,\ell;n\geq 0}$, on which $\pa_t$ acts by
$\pa_t \cdot u_i^{(n)} = u_i^{(n+1)}$. Thus we obtain a new proof of
Proposition 4 from \cite{FF1} (see also \cite{Five}).

\subsubsection{Theorem} The ring of functions $\C[Fl^\circ_A]$ on
$Fl^\circ_A=LG^{>0}/A'_+$ is isomorphic to the ring of
differential polynomials $\C[u_i^{(n)}]_{i=1,\ldots,\ell;n\geq 0}$, on
which the action of $p_{-1}$ is given by $p_{-1} \cdot u_i^{(n)} =
u_i^{(n+1)}$.

\subsubsection{mKdV hierarchy.}
The infinite-dimensional abelian Lie algebra $\ab/\ab_+$ acts from the
right on $Fl^\circ_A$. Hence we obtain an infinite set of
commuting flows on $AOp(\At)$, and an infinite set of commuting
evolutionary derivations of $\C[u_i^{(n)}]_{i=1,\ldots,\ell;n\geq
0}$. Denote by $\ol{p}_{-m}$ the derivation corresponding to $p_{-m}
\in \ab/\ab_+, m \in I$. In particular, we have: $\ol{p}_{-1} =
\pa_t$. The equation
\begin{equation}    \label{mKdVeq}
\pa_{t_m} u_i = \ol{p}_{-m} \cdot u_i, \quad \quad
i=1,\ldots,\ell,
\end{equation}
is called the $m$th equation of the {\em generalized modified KdV
hierarchy} (or mKdV hierarchy)
associated to $\g$.

\subsubsection{Miura transformation.}
We have the following commutative diagram of differential rings

\begin{equation}
\begin{CD}
(\C[Gr^\circ_A],p_{-1}) @>{\sim}>> (\C[v_i^{(n)}],\pa_t)
\end{CD}
\end{equation}
\begin{picture}(60,0)(0,0)
\put(54,-1){$\Big\downarrow$}
\put(56,3.5){\oval(1.5,2)[t]}
\put(91,-1){$\Big\downarrow$}
\put(93,3.5){\oval(1.5,2)[t]}
\end{picture}
$$\begin{CD}
(\C[Fl^\circ_A],p_{-1})\; @>{\sim}>> (\C[u_i^{(n)}],\pa_t)
\end{CD}$$

\noindent where the vertical arrows are embeddings and the horizontal
arrows are isomorphisms. Furthermore, the above diagram is compatible
with the action of $\ab/\ab_+$ on all four rings.

The embedding $\C[v_i^{(n)}]_{i=1,\ldots,\ell;n\geq 0} \arr
\C[u_i^{(n)}]_{i=1,\ldots,\ell;n\geq 0}$ is called the {\em Miura
transformation}. The corresponding map of spectra $AOp(\At) \arr
Op(\At)$ is simply the forgetting of the flat $B$--reduction
$\E'_B$.

Explicitly, given a Miura oper \eqref{Miura oper}, we view it as an
element of $\wt{Op}(\At)$ (\secref{twidle opers}) 
and take its projection onto $Op(\At)$,
i.e. apply a gauge transformation by an appropriate element of
$N[[t]]$ to bring it to the form \eqref{genoper2}. For instance, in
the case $\g={\mathfrak s}{\mathfrak l}_2$ we have the following
transformation,:
$$
\begin{pmatrix} 1 & - \frac{u}{2} \\ 0 & 1 \end{pmatrix} \left( \pa_t +
\begin{pmatrix} \frac{u}{2} & 0 \\ 1 & - \frac{u}{2} \end{pmatrix} \right)
\begin{pmatrix} 1 & - \frac{u}{2} \\ 0 & 1 \end{pmatrix}^{-1} = 
\pa_t + \begin{pmatrix} 0 & v \\ 1 & 0 \end{pmatrix}.
$$
Therefore the Miura transformation is in this case
\begin{equation}    \label{miura}
v \arr \frac{1}{4} u^2 + \frac{1}{2} \pa_t u.
\end{equation}

\subsubsection{Zero curvature representation.}

We start out with an analog of \lemref{trivial connection} in the
context of the flag manifold.

Recall that the big cell $\Fl^{\circ}\subset \Fl$ is isomorphic to the
group $LG^{>0}$. Consider the $LG^-$--bundle $\T$: $LG \arr \Fl =
LG^-\bs LG$. The fiber of $\T$ over $K^+ \in LG^{>0} = \Flo$ consists
of all elements $K$ of $LG$ that can be written in the form $K = K^-
K^+$, for some $K^-\in LG^-$. Thus the restriction of $\T$ to
$\Flo$ is canonically trivialized.

We are now in the setting of \lemref{action->connection}, where $M =
\Flo$, and $A = \Ap \subset \Ahat$, where $p \in \ab/\ab_+$. The group
$\Ap$ acts on $\Flo$ from the right. Hence we obtain for each $K^+ \in
LG^{>0} = \Flo$ a connection on the $LG^-$--bundle $\pi_{K^+}^*(\T)$
over $\Ap$ (here $\pi_{K^+}: \Ap \to \Flo$ is the $\Ap$--orbit of
$K^+$). The above trivialization of $\T$ induces a trivialization of
$\pi_{K^+}^*(\T)$, and allows us to write down an explicit formula for
this connection in the same way as in \lemref{trivial connection}.

\subsubsection{Lemma.}    \label{trivconn}
In the above trivialization of $\T$ the connection operator on
$\pi_{K^+}^*(\T)$ equals
\begin{equation}    \label{nablap}
\partial_t+(K^+(t) p K^+(t)\inv)^-,
\end{equation}
where $K^-(t)K^+(t)$ is the factorization of $K^+ e^{-tp}$.

\subsubsection{}
The isomorphism of \propref{pre-pre-oper}, reformulated for
$Fl^\circ_A$, assigns to each point $K^+ \in Fl^\circ_A$ a quadruple
$(\V,\nabla,\V^-,\V^+)$, where $\V$ is an $LG$--torsor on $\Ahat$ with
a flat connection $\nabla$, a flat $LG^-$--reduction $\V^-$, and an
$LG^+$--reduction in tautological relative position with $\nabla$. We
want to calculate the contractions $\nabla_p, p \in \ab/\ab_+$, of
$\nabla$ in the trivialization of $\V$, which preserves both $\V^-$
and $\V^+$, and has symbol $p$. Such a trivialization is unique,
because $LG^- \cap LG^{>0} = \on{Id}$. From \lemref{trivconn}, we
obtain the following formula:
$$
\nabla_p = \partial_t+(K^+(t) p K^+(t)\inv)^-.
$$
The flatness of the connection $\nabla$ leads to the zero curvature
equations
\begin{equation}    \label{zc}
\left[ \partial_{t_m} + (K^+(t) p_{-m} K^+(t)\inv)^-,\pa_{t_n} +
(K^+(t) p_{-m} K^+(t)\inv)^- \right] = 0.
\end{equation}

The $p_{-1}$ component of the connection is a modified affine oper
$$
\pa_t + (K^+(t) p_{-1} K^+(t)\inv)^- = \pa_t + p_{-1} + {\mathbf
h}(t),
$$
where ${\mathbf h}(t) \in \h[[t]]$ (cf. \secref{coord}). Therefore in
the special case $n=1$ we obtain the equation
\begin{equation}    \label{zcm}
\left[ \partial_{t_m} + (K^+(t) p_{-m} K^+(t)\inv)^-,\pa_t + p_{-1} +
{\mathbf h}(t) \right] = 0.
\end{equation}
This is the {\em zero curvature representation} of the $m$th mKdV
equation \eqref{mKdVeq}. Any solution for $(K^+(t) p_{-m}
K^+(t)\inv)^-$ as an element of $\g[[t]]$ is a differential
polynomial in $u_i(t) = \al_i({\mathbf h}(t))$, i.e., $L_m \in
\C[u_i^{(n)}] \otimes \g$. Hence formula \eqref{zcm} expresses
$\pa_{t_m} u_i$ as a differential polynomial in $u_j$'s (see
\cite{EF1,Five} for more detail).

Note that in contrast to the KdV equations, we do not have any
residual gauge freedom in equations \eqref{zcm}, because the
trivialization of $\V$ compatible with $\V^-$ and $\V_+$ is
unique.

\subsubsection{Example.}
Let us derive the mKdV equation, which is \eqref{zcm} in the case
$\g = {\mathfrak s}{\mathfrak l}_2, m=3$, following \cite{Five}.

We have:
\begin{equation}    \label{Lop}
\pa_t + p_{-1} + {\mathbf h} = \pa_t + \begin{pmatrix} \frac{1}{2} u
& z \\ 1 & - \frac{1}{2} u \end{pmatrix}.
\end{equation}

Now we have to compute $(K^+(t) p_{-3} K^+(t)^{-1})^-$. This can be
done recursively using the equation
$$
[p_{-1} + {\mathbf h},(K^+(t) p_{-3} K^+(t)^{-1})^-]=0
$$
(see \cite{Five}). It gives:
$$
(K^+(t) p_{-3} K^+(t)^{-1})^- = \begin{pmatrix} \frac{1}{2} u z^2
- \left(
\frac{1}{16} u^3 - \frac{1}{8} \pa_t^2 u \right) & z^2 + \left(
- \frac{1}{8} u^2 + \frac{1}{4} \pa_t u \right) z \\ z - \left(
\frac{1}{8} u^2 + \frac{1}{4} \pa_t u \right) & - \frac{1}{2} u
z^2 + \left( \frac{1}{16} u^3 - \frac{1}{8} \pa_t^2 u \right)
\end{pmatrix}.
$$
Substituting into formula \eqref{zcm}, we obtain:
\begin{equation}    \label{mkdv1}
\pa_{t_3} u = \frac{3}{8} u^2 \pa_t u - \frac{1}{4} \pa_t^3 u.
\end{equation}
This is the mKdV equation up to a slight redefinition of variables.
One can check that the corresponding equation on $v = \displaystyle
\frac{1}{4} u^2 + \frac{1}{2} \pa_t u$ (applying the
Miura transformation) is the
KdV equation \eqref{kdv3}.

\subsection{Generalized Affine Opers.}\label{general affine opers}

We now introduce generalized affine opers, which appear to provide the
broadest setting in which we can obtain integrable systems of KdV type
using the Drinfeld--Sokolov construction, \propref{DS gauge}. In
particular, their moduli spaces will turn out to be isomorphic to the
appropriate abelianized Grassmannians $\Gr_A$.

\subsubsection{} Let $\ab\subset L\g$ be a Heisenberg algebra with the
canonical filtration $\{ \ab^{\geq j} \}$, and let $p\in\ab^{\geq l},
l<0$, be an element with regular symbol $\ol{p}$ (then $p$ is
automatically regular as well). Let $\{L\g^{\geq j}\}$ be a filtration
compatible with $\ab$, in the sense of \secref{filtrations}. We
consider this as chosen once and for all, so that all superscripts
refer to this filtration, and we otherwise suppress it in the
notation.

Recall the notion of tautological relative position from
\defref{taut}, and the notation $\Ap = \{ e^{tp} \}$ for the
one-dimensional formal additive subgroup of $A/A_+$ generated by
$p$. We then have the following generalization of the notion of an affine
oper.

\subsubsection{Definition.}\label{all affine opers} 
An {\em affine $(A,p)$--oper} is a quintuple
$(\V,\nabla,\V_-,\V_+,\V^+)$, where $\V$ is an $LG$--torsor on $\Ap$
with a flat connection $\nabla$, a flat reduction $\V_-$ to $LG_-^X$,
and compatible reductions $\V_+$ and $\V^+$ to $LG_+$ and $LG^+$,
respectively. We require that $\V^+$ is in tautological relative
position with $\nabla_p$.

In the case when $X=\Pp$, we call an affine $(A,p)$--oper {\em
  generic} if the reductions $\V_-$ and $\V_+$ are in general position.

The moduli stack classifying affine $(A,p)$--opers is denoted by
$\AOp^X_{A,p}$, and the open substack of $\AOp^{\Pp}_{A,p}$
classifying generic affine $(A,p)$--opers is denoted by
$\AOp^\circ_{A,p}$.

\subsubsection{Remarks.}
\begin{enumerate}
\item[(1)] When $\ab$ is a smooth Heisenberg, so that $\ab^+\subset
\ab_+$, (\defref{canonical filtration}), we may require $LG^+ \subset
LG_+$ (e.g. for a filtration coming from a compatible gradation), and
hence the $LG_+$--reduction is redundant. But in general we do not
have such an inclusion, making the above definition somewhat
cumbersome.

\item[(2)] Since $LG^+$ preserves the $\ab$--filtration, the
$\V^+$--twist $(L\g)_{\V^+}$ carries a canonical filtration and hence
the relative position condition makes sense. Also note that compatible
reductions to $LG_+$ and $LG^+$ amount to a reduction $\V_+^+$ to
$LG_+\cap LG^+$.

\item[(3)] It is possible to define affine $(A,p)$ opers on an
arbitrary differential scheme $(S,\pa)$. \thmref{very general} will
then identify $\Gr_A$ (with the action of $p$)
with the moduli stack of such objects.

\item[(4)] For general Heisenbergs, there seems to be no simple
  generalization of the notion of $G$--oper. It is however possible to
  identify generic versions of affine opers with flat connections on
  appropriate finite--dimensional bundles. If $A$ is of Coxeter type
  (i.e. $LG$--conjugate to the principal Heisenberg), then we can always
  define a map from an open subset of the moduli of ``generic'' affine
  opers to the moduli of ordinary opers.

\end{enumerate}

The following key result follows from the generalized Drinfeld-Sokolov
gauge (\propref{DS gauge}):

\subsubsection{Proposition.}\label{general DS reduction}
Let $(\V,\nabla,\V_-,\V^+)$ be an affine oper on $\At$. Then $\V$ has
a unique reduction $\V^{A_+}\subset\V^+$ to $A_+$, such that
\begin{enumerate}
\item[(1)] the induced $A$--torsor $\V^A\subset\V$ is flat;
\item[(2)] $\V^{A_+}$ has the tautological relative position with
respect to 
$\nabla$.
\end{enumerate}

\subsubsection{Proof.}
First recall that the above relative position condition means that the
action of the vector field $p$ on $\Ap$ lies in the $\V^{A_+}$--twist
of the $A_+$--orbit $[-p]\in \ab/\ab_+ \subset L\g/\ab_+$.

We now use \propref{DS gauge} to reduce $\V^+$ to $A^+$.

Pick an arbitrary trivialization $t^+:LG^+\to\V^+$ of the
$LG^+$--bundle $\V^+$.  By the assumption on the connection $\nabla$,
we may pick this trivialization (after possibly using the
$LG^+$--action) so that $t(\nabla_{p})t\inv=\partial_p+p+q$, where $q$
is a section of the trivial $LG^+$--bundle on $\Dt$.  By \lemref{DS
gauge}, we may find an $M\in LG^{>0}$ so that $M\inv t\inv\nabla_p
tM=\partial+p+p_+$, with $p_+\in\ab^{>l}$.  Define $\V^{A^+}=tM(A^+)$,
the image of $A^+\subset LG$ under the isomorphism $tM:LG\to\V$. This
is an $A^+$--torsor in $\V$, since if $x=tM(a)\in\V^{A^+}$ and $b\in
A^+$, then $b\cdot x=tM(ab\inv)\in\V^{A^+}$.

Now suppose we pick a different trivialization $t':LG\to\V$, and
define $M'$ as above. It follows that the composition $N=M\inv t\inv
t'M'$ conjugates the connection $\partial+p+p'_+$ to the form
$\partial+p+p_+$, where both $p_+,p'_+\in\ab^{>l}$. By the uniqueness
statement of \propref{DS gauge}, it follows that $N\in A^+$. Hence
$x=tMa=t'M'a'$, where $a=Na'$, so that the two definitions of
$\V^{A^+}$ agree. Furthermore, the $A^+$--torsor structures agree as
well: $b\cdot x=t'M'(a'b\inv)=tM(Na'b\inv)=tM(ab\inv)$.

Now the intersection $\V^{A^+} \cap \V_+$ is an $A_+$--reduction
$\V^{A_+}$, since by definition the $LG_+$-- and $LG^+$--reductions
are compatible. It is clear that $\V^{A_+}$ satisfies the above
conditions.

\subsection{Generalized Drinfeld-Sokolov Hierarchies.}

In this subsection we extend the construction of integrable
systems to general Heisenberg subalgebras.

\subsubsection{}
Recall the definition of $\M^X_{A,p}$ from \secref{general M} and the
isomorphism between the abelianized Grassmannian $\Gr_A$ and
$\M^X_{A,p}$ from \propref{general iso}. The following results are
proved in the same way as in the case of the principal
Heisenberg. (See \remref{strongly regs} concerning strongly regular
elements.)

\subsubsection{Theorem.}\label{very general} 
Let $p \in \ab^{\geq l}, l<0$, be a strongly
regular element.
\begin{enumerate}
\item[(1)] There is a natural morphism from $\M^X_{A,p}$
to $\AOp^X_{A,p}$, and hence from $\Gr_A$ to $\AOp^X_{A,p}$.

\item[(2)] The morphism between the abelianized
Grassmannian $\Gr_A$ and the moduli stack $\AOp^X_{A,p}$ of affine
$(A,p)$--opers on $\Ap$ is an isomorphism.

\item[(3)] The above isomorphism identifies the big cell
$\Gro_A$ of $\Grp_A$ and the moduli stack $\AOp_{A,p}^\circ$
of generic affine opers on $\Ap$.
\end{enumerate}

\subsubsection{} Recall that $\Gro_A$ is the quotient of the
scheme $Gr^\circ_A= G \bs LG_+/A'_+$ by the trivial action of
a group $A_0$ (see \secref{Gro}). On the other hand, the argument of
\lemref{silly} applied to $\AOp_{A,p}^\circ$ shows that it is the
quotient of a scheme $AO_{A,p}^\circ$ by the trivial action of
$A_0$. Therefore we obtain

\subsubsection{Corollary.}
For each strongly regular $p \in \ab^{\geq l}, l<0$, there is a
canonical isomorphism $Gr^\circ_A \simeq AO^\circ_{A,p}$.

\subsubsection{Canonical form.}
Choosing a canonical form of an affine
$(A,p)$--oper allows us to identify the ring of functions on
$AO^\circ_{A,p}$ (and hence on $Gr^\circ_A$) with the ring of
differential polynomials. This can be done as follows.

Given an $(A,p)$--oper, we can choose a trivialization of the
underlying $LG$--bundle, which preserves the reductions $\V^+,
\V_-$. Such a trivialization is unique up to the action of the
finite-dimensional group $LG^{+} \cap LG_-$. Requiring the
trivialized connection
to have symbol $p$ reduces gauge group to the unipotent
group $R=LG^{>0}\cap LG_-$. The
connection operator of the affine oper then reads:
\begin{equation}    \label{canform1}
\pa_t + p + b(t), \quad \quad b(t) \in (L\g^{>l} \cap L\g_-)[[t]].
\end{equation}
Denote the space of all such operators by $\wt{AO}^\circ_{A,p}$. Then
$AO^\circ_{A,p}$ is the quotient of $\wt{AO}^\circ_{A,p}$ by the gauge
action of $R[[t]]$. Denote by $\rr$ the (nilpotent) Lie algebra of the
group $R$. The following lemma, which is a generalization of
\lemref{canreps}, is proved along the lines of \propref{DS gauge}.

\subsubsection{Lemma.} The action of $R[[t]]$ on $\wt{AO}^\circ_{A,p}$
is free. Moreover, each $R[[t]]$--orbit in $\wt{AO}^\circ_{A,p}$
contains a unique operator of the form
\begin{equation}    \label{canform2}
\pa_t + p + v(t), \quad \quad v(t) \in V[[t]],
\end{equation}
where $V \subset \rr$ is such that $\rr = V \oplus \on{Im} \on{ad}
p$. Thus, we can identify $AO^\circ_{A,p}$ with $V[[t]]$.

\subsubsection{}
Choosing a basis $\{ v_i \}$ of $V$, we identify $\C[AO^\circ_{A,p}]$,
and hence $\C[Gr^\circ_A]$, with the ring of differential polynomials
$\C[v_i^{(n)}]$. Moreover, the action of $p$ on $Gr^\circ_A$
corresponds to the standard action of $\pa_t$ on $\C[v_i^{(n)}]$.

The action of $A/A_+$ on $Gr^\circ_A$ gives rise to an infinite
hierarchy of commuting flows on $AO^\circ_{A,p}$, and hence commuting
evolutionary derivations on $\C[v_i^{(n)}]$. They form the {\em
generalized Drinfeld-Sokolov hierarchy} associated to $A$ and the
strongly regular element $p \in \ab$. One can write these equations
down explicitly in the {\em zero curvature form} following
\secref{zcf}. Namely, for each $q \in \ab/\ab_+$ we have the equation
\begin{equation}    \label{genDS}
[\pa_s + (K_+(t) q K_+(t)\inv)_-,\pa_t + p + v(t)] = 0.
\end{equation}
Here $K_+(t)$ comes from the factorization $K_-(t) K_+(t) = K_+
e^{-tp}$ along the $\Ap$--orbit of $K_+ \in G \bs LG_+/A'_+ \simeq
Gr^\circ_A$. Note that
$$
\pa_t + p + v(t) = \pa_t + (K_+(t) p K_+(t)\inv)_-.
$$
The equation \eqref{genDS} is invariant under the residual gauge group
$R[[t]]$. Hence the above equations with $q$ running over $\ab/\ab_+$
really define commuting evolutionary derivations on the ring of
differential polynomials $\C[v_i^{(n)}]$.

\subsubsection{Remark: How Many Hierarchies?} 
Suppose $A$ and $A'$ are $LG_+$--conjugate
Heisenbergs, so that $gA_+g\inv=A'_+$ for some fixed $g\in LG_+$. Then
there is an isomorphism $\Gr_{A'} \simeq \Gr_A$ intertwining the
actions of $A'/A'_+$ and $A/A_+$. Namely, we have a well--defined map
on double quotients $i_K: LG_- \cdot M \cdot A'_+ \mapsto LG_- \cdot M
g \cdot A_+$. It follows that the integrable systems associated with
$A$ and $p \in \ab$ and $A'$ and $p'=gpg\inv \in \ab'$ are equivalent.

Thus our construction associates an integrable system to each
$LG_+$--conjugacy class of pairs $(\ab,p)$, where $\ab$ is an
arbitrary Heisenberg subalgebra of $L\g$, and $p$ is a strongly
regular element of $\ab/\ab_+$. (Different strongly regular $p$ of the
same $\ab$ give different presentations of the same underlying
integrable system.)

Let $LH^{[w]}$ denote a graded Heisenberg of the same type $[w]$ as
$A$. As we remarked in \secref{difference of A's} and \secref{picture
  of A}, the ``difference'' between $A$ and $LH^{[w]}$ is measured by
the finite--dimensional group scheme $A^+/A_+$. In particular, this
group acts on $\Gr_A$ (commuting with the flows) and the quotient is
isomorphic to $\Gr_{LH^{[w]}}$ (with its natural flows). However this
  isomorphism does {\em not} lead to an isomorphism of the big cells,
  and so the resulting integrable systems on $\Gro_A$ and
  $\Gro_{LH^{[w]}}$ can be quite different.

\subsubsection{} There is an obvious version of the above
construction, in which the abelianized Grassmannian $Gr^\circ_A$ is
replaced by its flag manifold version $Fl^\circ_A$, as in
\secref{flags}. In particular, we obtain an identification between
$Fl^\circ_A$ and an appropriate moduli space of affine Miura
opers. The corresponding flows form the generalized mKdV hierarchy.

One can also introduce ``partially modified'' hierarchies by
considering moduli spaces that are intermediate between $\Grp_A$ and
$\Fl_A$, namely, the moduli space of $G$--bundles $\E$ on $X$ with a
reduction of $\E|_{\Dx}$ to $A_+$ and a reduction of $\E|_{0 \in X}$
to a parabolic subgroup $P$ of $G$. (The partial flag at $0$ is chosen
so as to provide a partial splitting of the $A$--filtration.)

\subsubsection{Examples.}

In the case when $\ab$ is a graded Heisenberg subalgebra, the
corresponding
generalized Drinfeld-Sokolov hierarchy was introduced in
\cite{dGHM}. In this case, one can promote the canonical filtration on
$L\g$ into a $\Z$--gradation, which simplifies the study of the
equations. Delduc and Feher \cite{DF,Feher} have described explicitly
the strongly regular elements $p$ of the graded Heisenberg
subalgebras.

The most widely known example is of course the case when $\ab$ is the
principal Heisenberg subalgebra and $p=p_{-1}$, which corresponds to
the generalized KdV hierarchies discussed above. In \cite{Ba}, Balan
studies, along the lines of \cite{FF0,FF1,EF1}, the hierarchy
corresponding to $p=p_{-3}$ in the case of the principal Heisenberg
subalgebra of $L{\mathfrak s}{\mathfrak l}_2$. The other well-known
example is the generalized AKNS (or non-linear Schr\"odinger)
hierarchy, which corresponds to the case of the homogeneous Heisenberg
subalgebra (see \cite{NLS}).

Finally, we have worked out explicitly the case of the simplest
non--smooth Heisenberg subalgebra, which was introduced in
\secref{non-smooth}. We plan to present this hierarchy elsewhere.

\section*{Acknowledgments.} We would like to acknowledge very helpful
conversations with A. Beilinson, P. Burchard, D. Gaitsgory,
A. Genestier, M. Rothstein, and especially R. Donagi, who explained to
us the theory of spectral curves, alerted us to the role of singular
curves, and provided useful comments on an early draft.


\begin{thebibliography}{Dri}

\bibitem[AB]{AB} M. Adams and M. Bergvelt, The Krichever map, vector
bundles over algebraic curves, and Heisenberg algebras,
Comm. Math. Phys. {\bf 154} (1993) 265-305.

\bibitem[AMP]{AMP} A. \'{A}lvarez, J. Mu\~{n}oz and F. Plaza, The
algebraic formalism of soliton equations over arbitrary base fields,
Workshop on Abelian Varieties and Theta Functions (Morelia 1996),
Aportaciones Mat. Investig. 13, Soc. Mat. Mexicana (3-40)
1998 (alg-geom/9606009).

\bibitem[Ba]{Ba} A. Balan, Les \'equations mKdV g\'en\'eralis\'ees,
Preprint Ecole Polytechnique, No 98-5, April 1998.

\bibitem[BL]{BL} A. Beauville and Y. Laszlo, Conformal blocks and
generalized theta functions, Comm. Math. Phys. 164 (1993) 385-419.

\bibitem[BB]{BB} A. Beilinson and J. Bernstein. A Proof of Jantzen
Conjectures, Advances in Soviet Mathematics, Vol. 16 Part 1 (1-50) 1993. 

\bibitem[BD1]{opers} A. Beilinson and V. Drinfeld, Opers, Preprint,
1994.

\bibitem[BD2]{Hecke} A. Beilinson and V. Drinfeld, Quantization of
Hitchin's Integrable System and Hecke Eigensheaves, in
preparation.

%\bibitem[BZ]{BZ} D. Ben--Zvi. Cartan subgroups of loop groups and
%affine Springer fibers, In preparation.

\bibitem[BS]{new} A. Blum and U. Stuhler. Drinfeld modules and
elliptic sheaves, Vector bundles on curves -- new directions
(eds. S. Kumar et al.) LNM 1649, 1997.

\bibitem[Ch1]{Ch1} I. Cherednik, Group interpretation of Baker
  functions and $\tau$--function, Uspekhi Mat.  Nauk 38, No. 6 (1983)
  133--134.  

\bibitem[Ch2]{Ch2} I. Cherednik, Determination of
  $\tau$--functions for generalized affine Lie algebras, Funct. Anal.
  Appl.i Prilozhen. 17, No. 3 (1983) 93--95.

\bibitem[Ch3]{Ch3} I. Cherednik, Functional realizations of basic
representations of factorizable Lie groups and Lie
algebras, Funct. Anal.  Appl. {\bf 19}, 193-206 (1985)

\bibitem[CC]{CC} C. Contou--Carr\'{e}re, Jacobienne locale, groupe de
bivecteurs de Witt universel, et symbole mod\'{e}r\'{e},
C. R. Acad. Sci. Paris, s\'{e}rie I 318 (1994) 743-746.

\bibitem[DJKM]{DJKM} E. Date, M. Jimbo, M. Kashiwara, T. Miwa,
Transformation groups for soliton equations, In: M. Jimbo, T. Miwa
(eds.) Non-linear Integrable Systems -- Classical Theory and Quantum
Theory, pp. 39-120, Singapore: World Scientific 1983.

\bibitem[dGHM]{dGHM} M. de Groot, T. Hollowood and L. Miramontes,
Generalized Drinfeld-Sokolov hierarchies, Comm. Math. Phys. 145 (1992)
57-84.

\bibitem[DF]{DF} F. Delduc and L. Feh\'{e}r, Regular conjugacy classes
in the Weyl group and integrable hierarchies, J. Phys A28 (1995)
5843-5882.

\bibitem[D]{Donagi} R. Donagi, Decomposition of spectral
covers, Ast\'{e}risque 218 (1993) 145-175.

\bibitem[DG]{DG} R. Donagi and D. Gaitsgory, The gerbe of Higgs bundles, 
Preprint math.AG/0005132.

\bibitem[DM]{DM} R. Donagi and E. Markman, Spectral covers,
algebraically completely integrable Hamiltonian systems, and moduli of
bundles, Integrable systems and quantum groups, Lecture Notes in
Math. 1620, Springer 1995.

\bibitem[Dr]{Dr} V. Drinfeld, Commutative subrings of certain
noncommutative rings, Funct. Anal Appl. 11, No. 1 (1977) 11-14.

\bibitem[DSi]{DSi} V. Drinfeld and C. Simpson, $B$--structures on
$G$--bundles and local triviality, Math. Res. Lett. 2, No. 6 (1995)
823--829.

\bibitem[DS]{DS} V. Drinfeld and V. Sokolov, Lie algebras and
equations of Korteweg-deVries type, Journal of Soviet Mathematics,
vol. 30 (1985) 1975-2035.

\bibitem[EF1]{EF1} B. Enriquez and E. Frenkel, Equivalence of two
approaches to integrable equations of KdV type, Comm. Math. Phys. 185
(1997) 211--230.

\bibitem[EF2]{EF2} B. Enriquez and E. Frenkel, Geometric
interpretation of the Poisson structure in affine Toda field
theories, Duke Math. J. 92 (1998) 459--495.

\bibitem[Fa]{Fa} G. Faltings, Stable $G$-bundles and projective
connections, J. Alg. Geom. 2 (1993) 507--568.

\bibitem[Fe]{Feher} L. Feh\'{e}r, KdV type systems and $\W$--algebras in
the Drinfeld--Sokolov approach, Preprint hep-th/9510001.

\bibitem[FHM]{FehHarMar} L. Feh\'{e}r, J. Harnad and I. Marshall,
Generalized Drinfeld--Sokolov reductions and KdV type hierarchies,
Comm. Math. Phys. 154 (1993) 181-214.

\bibitem[FF1]{FF0} B.~Feigin, E.~Frenkel, Integrals of motion and
quantum groups, in Lect. Notes in Math. 1620, pp. 349-418, Springer
Verlag, 1995.

\bibitem[FF2]{FF1} B. Feigin and E. Frenkel, Kac-Moody groups and
integrability of soliton equations, Invent. Math. 120 (1995) 379-408.

\bibitem[FF3]{NLS} B. Feigin and E. Frenkel, Integrable hierarchies
and Wakimoto modules, in Differential Topology, Infinite-Dimensional
Lie Algebras, and Applications: D. B. Fuchs' 60th Anniversary
Collection, A. Astashkevich and S. Tabachnikov (eds.), pp. 27-60, AMS
1999.
  
\bibitem[F]{Five} E. Frenkel, Five lectures on soliton equations, in
Integrable Systems, C.-L. Terng and K. Uhlenbeck (eds.), Surveys in
Differential Geometry, vol. 4, pp. 131-180. International Press, 1998.

\bibitem[Gin]{Ginzburg} V. Ginzburg, Perverse Sheaves on a Loop Group
and Langlands' Duality, Preprint alg--geom/9511007.

\bibitem[Ha]{Har} G. Harder, Halbeinfache Gruppenschemata \"uber
Dedekindringen, Invent. Math. 4 (1967) 165-191.

\bibitem[Hi]{Hi} N. Hitchin, Stable bundles and integrable
systems, Duke Math. J. 54 (1990) 91--114.

\bibitem[Kac1]{Kac1} V.G. Kac, Infinite-dimensional algebras, Dedekind
$\eta$--function, classical M\"obius function and the very strange
formula, Adv. Math. 30 (1978) 85-136.

\bibitem[Kac2]{Kac} V.G. Kac. Infinite-dimensional Lie algebras. Third
Edition. Cambridge University Press 1990.

\bibitem[KP]{KP} V. Kac and D. Peterson, 112 constructions of the
basic representation of the loop group of E8, Proceedings of
``Anomalies, geometry, topology'' (Argonne, 1985), pp. 276-298, World
Scientific 1985. 

\bibitem[KSU]{KSU} T. Katsura, Y. Shimizu and K. Ueno, Formal Groups
and conformal field theory over $\Z$, In Integrable systems in quantum
field theory and statistical mechanics. Adv. Stud. Pure Math
19. Academic Press, Boston (1989) 347-366.

\bibitem[KL]{KL} D. Kazhdan and G. Lusztig, Fixed point varieties on
affine flag manifolds, Israel J. of Math. 62, No. 1 (1988) 129-168.

\bibitem[Kos]{Kostant TDS} B. Kostant, The principal
three--dimensional subgroup and the Betti numbers of a complex simple
Lie group, Amer. J. Math. 81 (1959) 973--1032.

\bibitem[Kr]{Kr} I. Krichever, Algebro-geometric construction of the
Zakharov-Shabat equations and their periodic
solutions, Dokl. Akad. Nauk SSSR 227, No. 2 (1976) 291-294.

\bibitem[LS]{LS} Y. Laszlo and C. Sorger, The line bundles on the
stack of parabolic $G$--torsors over curves and their
sections, Ann. Sci. Ec. Norm. Sup. 30, No. 4 (1997) 499--525. 

\bibitem[Lau1]{Lau} G. Laumon, Transformation de Fourier
g\'{e}n\'{e}ralis\'{e}e, Preprint alg-geom/9603004.

\bibitem[Lau2]{Lau2} G. Laumon, Sur les modules de Krichever, Preprint
Universit\'e de Paris-Sud 86 T 20.

\bibitem[LMB]{LMB} G. Laumon and L. Moret--Bailly, Champs
Alg\'ebriques. Ergebnisse der Mathematik und ihrer Grenzgebiete, 3. Folge 39.
Springer-Verlag Berlin 2000.

\bibitem[LM]{LM} Y.~Li and M.~Mulase. Prym varieties and integrable
systems, Comm. Anal. Geom. 5 No. 2 (1997) 279-332.

\bibitem[McK]{McKean} H.~McKean, Is there an infinite--dimensional
algebraic geometry? Hints from KdV. Proc. Symp. Pure Math 49 Vol.I,
pp. 27--37, AMS 1989.
 
\bibitem[MV]{MV} I.~Mirkovi\'c and K.~Vilonen, Perverse sheaves on
loop Grassmannians and Langlands duality, Math Res. Lett. 7 (2000) 
13--24.

\bibitem[M1]{M1} M.~Mulase. Cohomological structure in soliton
equations and Jacobian varieties, J. Diff. Geom. 19 (1984) 403-430.

\bibitem[M2]{M2} M.~Mulase. Category of vector bundles on algebraic
curves and infinite--dimensional Grassmannians, Int. J. Math. 1 (1990)
293-342.

\bibitem[M3]{M3} M. ~Mulase. Algebraic theory of the KP equations.
Perspectives in mathematical physics, Conf. Proc. Lect. Notes Math.
Phys. III. Intenat. Press, Cambrdige MA (1994) 151-217.

\bibitem[Mum]{Mum} D. Mumford, An algebro-geometric construction of
commuting operators and of solutions to the Toda lattice equation,
Korteweg-deVries equation and related non-linear equations,
Int. Symposium on Algebraic Geometry, Kyoto, pp. 115-153, Kinokuniya
Book Store 1977.

\bibitem[N1]{N1} A. Nakayashiki, Structure of Baker-Akhiezer modules of
pricipally polarized abelian varieties, commuting partial differential
operators and associated integrable systems, Duke Math J. 62 (1991)
315--358.

\bibitem[N2]{N2} A. Nakayashiki, 
Commuting partial differential operators and vector bundles
over abelian varieties, Amer. J. Math. 116 (1994), 65--100.

\bibitem[PS]{PS} A. Pressley and G. Segal, Loop Groups. Oxford
University Press 1986.

\bibitem[Ro1]{Ro1} M. Rothstein, Connections on the total Picard sheaf
and the KP hierarchy, Acta Applicandae Mathematicae 42 (1996) 297-308.

\bibitem[Ro2]{Ro2} M. Rothstein, Sheaves with connection on abelian
varieties, Duke Math J. 84, No. 3 (1996) 565--598.

\bibitem[Sor]{Sor} C. Sorger. Lectures on moduli of principal
$G$--bundles over algebraic curves. School on Algebraic Geometry
(Trieste 1999) ICTP Lect. Notes 1, Abdus Salam ICTP, Trieste (2000)
1--57.

\bibitem[Tel]{Tel} C. Teleman, Borel--Weil--Bott theory on the moduli
stack of $G$--torsors over a curve, Invent. Math. 134 (1998) 1--57.

\bibitem[SW]{SW} G. Segal and G. Wilson, Loop groups and equations of
KdV type, Publ. Math. IHES 61 (1985) 5-65.

\bibitem[W]{W} G. Wilson. Habillage et fonctions $\tau$, C. R. Acad.
Sci. Paris 299 S\'er. I no. 13 (1984) 587-590.

\end{thebibliography}
\end{document}